\documentclass{article}%
\usepackage[letterpaper,top=1.5cm,bottom=2cm,left=2.5cm,right=2.5cm,marginparwidth=2cm]{geometry}
\usepackage{amssymb}
\usepackage{amsfonts}
\usepackage{amsmath}
\usepackage{graphicx}%
\usepackage{csquotes}
\usepackage{subfigure}
\usepackage[colorlinks=true, allcolors=blue]{hyperref}
\usepackage{tikz}
\providecommand{\U}[1]{\protect\rule{.1in}{.1in}}
\begin{document}

\title{Some Laplacian eigenvalues can be computed by matrix perturbation}
\author{Piet Van Mieghem\thanks{ Faculty of Electrical Engineering, Mathematics and
Computer Science, P.O Box 5031, 2600 GA Delft, The Netherlands; \emph{email}:
P.F.A.VanMieghem@tudelft.nl } \; and Yingyue Ke\thanks{ Faculty of Electrical Engineering, Mathematics and
Computer Science, P.O Box 5031, 2600 GA Delft, The Netherlands; \emph{email}:
y.y.ke@tudelft.nl }}
\date{Delft University of Technology}
\maketitle

\begin{abstract}
Based on matrix perturbation theory and Euler summation, closed-form analytic expansions are studied for a Laplacian eigenvalue of an undirected, possibly weighted graph, which is close to a unique degree in that graph. 
An approximation is presented to provide an analytic estimate of a Laplacian eigenvalue and complements bounds on Laplacian eigenvalues in spectral graph theory. 
We numerically analyze how the structure of graphs influences the convergence of the corresponding Euler series.
Moreover, we obtain the explicit form of the perturbation Taylor series and its Euler series of almost regular graphs in which only one node has a unique degree and all remaining nodes have the same degrees.
We find that the Euler series possesses a superior convergence range than the perturbation Taylor series for almost regular graphs.

\end{abstract}

% \begin{keywords}
% Laplacian eigenvalue, Euler series, Matrix perturbation, Graph theory.
% \end{keywords}
% \begin{AMS}
% 05C50, 15A18. 
% \end{AMS}

% Sample article for the Electronic Journal of Linear Algebra

%%%%%%%%%%%%%%%%%%%%%%%%%%%%%%%%%%%%%%%

%%%%%%%%%%%%%%%%%%%%%%%%%%%%%%%%%%%%%%%%%%%%%%%%%%%%%%%%%%%%%
\section{Introduction} \label{intro-sec}
The eigenvalues of a possibly weighted Laplacian play
an important role in linear dynamics on a graph, as illustrated in
\cite{PVM_graphspectra_second_edition},\cite{PVM_pseudo_inverse_Laplacian} and
in Markov theory \cite{PVM_PAComplexNetsCUP}, where the infinitesimal generator is minus a weighted
Laplacian. There exist many bounds on Laplacian eigenvalues $\mu_1\geq \mu_2\geq\cdots\mu_{N-1}\geq 0$ of an $N\times N$ Laplacian matrix $Q$ of a graph $G$ on $N$ nodes and $L$ links. The
Brouwer-Haemers bound \cite{Brouwer_Haemers_LAA2008} is $\mu_{k}\geq
d_{\left(  k\right)  }-k+2$, where $d_{\left(  k\right)  }$ is the $k$-th
largest degree in the graph, i.e. $d_{\left(  1\right)  }\geq d_{\left(
2\right)  }\geq\ldots\geq d_{\left(  N\right)  }$. An upper bound for the
largest Laplacian eigenvalue \cite{PVM_graphspectra_second_edition} is
$\mu_{1}\leq\min\left(  N,\max_{l\in\mathcal{L}}\left(  d_{l^{+}}+d_{l^{-}%
}\right)  \right)  $, where $d_{l^{+}}$ and $d_{l^{-}}$ are the nodal degrees
at the left- and right-hand node of a link $l$; clearly $d_{l^{+}}+d_{l^{-}%
}\leq2d_{\max}$.

Here, we present analytic expansions of some Laplacian eigenvalues, deduced from matrix perturbation theory around a node with a unique degree. 
From the general matrix perturbation theory in \cite[Sec. 10.10]%
{PVM_graphspectra_second_edition}, Section \ref{sec_theory} derives expansions for a Laplacian eigenvalue. 
After briefly reviewing the theory in Section \ref{sec_brief_review_theory}, Section \ref{sec_coefficients_perturbation_series} derives the explicit coefficients of the perturbation Taylor series (\ref{Taylor_perturbation_series}). 
Section \ref{sec_convergence_perturbation_Taylor_series} investigates the convergence of the perturbation Taylor series, while Section \ref{sec_Euler_series} proposes an Euler $t$-series (\ref{mu_perturbation_Eulersummation_general_t})  with superior convergence range than the perturbation Taylor series and the approximation (\ref{lambda_Euler_upto_order4}), which provides an analytic estimate of a Laplacian eigenvalue close to a unique degree in any graph and complements bounds on Laplacian eigenvalues. 
Unfortunately, even the Euler $t$-series $\xi_{q;K}$ in (\ref{mu_perturbation_Eulersummation_general_t}) does not always converge. 
Since deriving precise mathematical convergence conditions is dauntingly complex, Section \ref{sec_performance_analysis}  numerically investigates the convergence of the Euler series $\xi_{q;K}$. 
Section \ref{sec:almost_regular} provides an explicit form of the perturbation Taylor series (\ref{perturbed_eigenvalue_zeta_almost_reg_graph}) and its Euler series
(\ref{perturbed_eigenvalue_zeta_almost_reg_graph_Euler}) of almost regular graphs in which only one node has a unique degree and all remaining nodes have the same degrees.
The final section \ref{sec_conclusion} concludes our performance evaluation of Euler $t$-series $\xi_{q;K}$.
%%%%%%%%%%%%%%%%%%%%%%%%%%%%%%%%%%%%%%%
\section{Theory}
\label{sec_theory}
\subsection{Brief review of matrix perturbation theory}
\label{sec_brief_review_theory}
Perturbation theory allows us to compute the eigenvalue and corresponding
eigenvector of a $N \times N$ matrix $W\left(  \zeta\right)  =W+\zeta B$ in terms of the
spectrum of the matrix $W$. Perturbation theory \cite[Sec. 10.10]%
{PVM_graphspectra_second_edition} assumes that the perturbation parameter
$\zeta$ is sufficiently small so that we may regard $W\left(  \zeta\right)  $
as the perturbation of the original symmetric matrix $W$ by a matrix $B$,
which is not necessarily symmetric. We limit ourselves to a simple eigenvalue
$\lambda\left(  W\right)  $ of the matrix $W$ with multiplicity one.
Perturbation theory for higher multiplicity eigenvalues of $W$ is involved and
here omitted. We denote by $x\left(  \zeta\right)  $ the $N\times1$
eigenvector of $W\left(  \zeta\right)  $ belonging to the eigenvalue
$\lambda\left(  \zeta\right)  $. As shown in \cite[pp. 60-70]{Wilkinson}, both the eigenvector
$x\left(  \zeta\right)  $ and the eigenvalue $\lambda\left(  \zeta\right)  $ are analytic
functions of $\zeta$ around the point $\zeta_{0}=0$ and can be represented by
a Taylor series%
\begin{align}
x\left(  \zeta\right)   &  =x+\zeta z_{1}+\zeta^{2}z_{2}+\cdots=\sum
_{j=0}^{\infty}z_{j}\zeta^{j}\label{eigenvector_expansion_zeta}\\
\lambda\left(  \zeta\right)   &  =\lambda+\zeta c_{1}+\zeta^{2}c_{2}%
+\cdots=\sum_{j=0}^{\infty}c_{j}\zeta^{j} \label{eigenvalue_expansion_zeta}%
\end{align}
where $x\left(  0\right)  =x$ is the eigenvector of $W$ and $\lambda\left(
0\right)  =\lambda$ is its corresponding simple eigenvalue. We choose
$x=x_{q}$ as the normalized eigenvector of $W$ corresponding to $\lambda
=\lambda_{q}$. The convergence of the Taylor series
(\ref{eigenvector_expansion_zeta}) and (\ref{eigenvalue_expansion_zeta})
assumes that $\zeta$ is \textquotedblleft sufficiently\textquotedblright%
\ small, in particular $\left\vert \zeta\right\vert <R$, where \thinspace$R$
is the radius $R$ of convergence of the above Taylor series. The entire
difficulty in perturbation theory lies in the determination of the radius $R$
of convergence of the Taylor series (\ref{eigenvector_expansion_zeta}) and (\ref{eigenvalue_expansion_zeta}).

Since the vector $z_{j}$ can be written as a linear combination of
the eigenvectors $x_{k}$ of matrix $W$, we have%
\begin{equation*}
z_{j}=\sum_{k=1}^{n}\beta_{jk}x_{k} \label{def_zj_as_lin_comb_xk}%
\end{equation*}
where the coefficients $\beta_{jm}=x_{m}^{T}z_{j}=z_{j}^{T}x_{m}\neq\beta
_{mj}$. The particular case $j=0$, where $z_{0}=x_{q}$, indicates that
$\beta_{0k}=\delta_{kq}$. 
With our eigenvector scaling choice $\beta_{jq}=0$ for $j>0$, the  coefficients $c_{k}$ \cite[Sec. 10.10]
{PVM_graphspectra_second_edition} simplifies considerably to
\begin{equation}
c_{1}=x_{q}^{T}Bx_{q} \label{coeff_c1_in_B}%
\end{equation}
and
\begin{equation}
c_{j}=\sum_{k=1;k\neq q}^{n}\beta_{j-1,k}x_{q}^{T}Bx_{k}\hspace{1cm}\text{for
}j>1 \label{coeff_cj}
\end{equation}
where, 
\begin{equation}
\beta_{1r}=\frac{x_{r}^{T}Bx_{q}}{\lambda_{q}-\lambda_{r}}\hspace
{1cm}\text{for }r\neq q \label{beta_1r_r_not_q}%
\end{equation}
and
\[
\beta_{jr}=\frac{1}{\lambda_{r}-\lambda_{q}}\sum_{l=1;l\neq q}^{n}\left\{
\sum_{k=1}^{j-1}\beta_{kr}\beta_{j-k-1,l}x_{q}^{T}Bx_{l}-\beta_{j-1,l}%
x_{r}^{T}Bx_{l}\right\}  \hspace{1cm}\text{for }r\neq q, j>1.
\]
The scaling choice $\beta_{0,l}=\delta_{lq}$ and $\beta_{jq}=0$ for $j\geq1$
simplifies, for $r\neq q$, to a recursion in $\beta_{jr}$
\begin{equation}
\beta_{jr}=\frac{\beta_{j-1;r}x_{q}^{T}Bx_{q}}{\lambda_{r}-\lambda_{q}}%
+\frac{1}{\lambda_{r}-\lambda_{q}}\sum_{l=1;l\neq q}^{n}\left\{  \sum
_{k=1}^{j-2}\beta_{kr}\beta_{j-k-1,l}x_{q}^{T}Bx_{l}-\beta_{j-1,l}x_{r}%
^{T}Bx_{l}\right\}  \label{beta_jr_recursion}%
\end{equation}
which can be iterated up to any desired integer value of $j$.

We apply the perturbation theory in \cite[Sec. 10.10]%
{PVM_graphspectra_second_edition} to a weighted Laplacian $\widetilde
{Q}=\widetilde{\Delta}-\widetilde{A}$, where $\widetilde{A}$ is a weighted
adjacency matrix and the diagonal matrix is $\widetilde{\Delta}=$ diag$\left(
\widetilde{A}u\right)  $, where $u$ is the all-one vector. In the unweighted
case, we omit the tilde and write the Laplacian $Q=\Delta-A$, where $A$ is the
$N\times N$ zero-one, symmetric adjacency matrix of a simple graph $G$ without
self-loops, i.e. $a_{jj}=0$ for any node $j$ in the graph $G$. The graph $G$
has $N$ nodes and $L$ links.   Thus, in our case, the $N\times N$
matrix $W=$ $\widetilde{\Delta}$ and the $N\times N$ perturbating matrix $B=$
$\widetilde{A}$. Hence, $W\left(  \zeta\right)  =\widetilde{\Delta}%
+\zeta\widetilde{A}$ and $W\left(  -1\right)  =\widetilde{Q}$. We will denote
the eigenvalue $\lambda\left(  \zeta\right)  $ in
(\ref{eigenvalue_expansion_zeta}) of the $N\times N$ matrix $W\left(  \zeta\right)
=\widetilde{\Delta}+\zeta\widetilde{A}$ by $\xi\left(  \zeta\right)  $. 
Similarly, $\zeta(1)$ with perturbation parameter $\zeta = 1$ denotes the eigenvalue   \cite[art. 30]
{PVM_graphspectra_second_edition} of the signless Laplacian $W\left(  -1 \right)  =\overline{Q}$.
The entire challenge is caused by the relatively large perturbation parameter
$| \zeta |=1$, which may lead to diverging series (\ref{eigenvector_expansion_zeta}) and (\ref{eigenvalue_expansion_zeta}).

\subsection{Coefficients and the perturbation Taylor series}
\label{sec_coefficients_perturbation_series}
The diagonal matrix $\widetilde{\Delta}$ has a rather obvious spectral
structure, which greatly simplifies the perturbation computations. The
eigenvalue $\lambda_{j}\left(  W\right)  =\lambda_{j}\left(  \widetilde
{\Delta}\right)  =\widetilde{d}_{j}$ is equal to the $j$-th diagonal element
of $\widetilde{\Delta}$, which corresponds to the nodal strength
$\widetilde{d}_{j}=\left(  \widetilde{A}u\right)  _{j}$ of node $j$. The
corresponding, normalized eigenvector $e_{j}$ is the basic vector with all
zeros, except an entry 1 for the $j$-th component, i.e. $\left(  e_{j}\right)
_{k}=\delta_{jk}$, where the Kronecker $\delta_{jk}=1$ if $k=j$, else
$\delta_{jk}=0$. The eigenvalue $\lambda_{j}\left(  W\right)  $ can possess a
multiplicity $m_{j}$, but the orthogonal eigenvector matrix remains the
identity matrix $I$ for $W=$ $\widetilde{\Delta}$, which avoids
complications of Jordan forms in the general case.

In order to allow direct application of the perturbation theory in \cite[Sec.
10.10]{PVM_graphspectra_second_edition}, we assume that $\lambda
_{q}=\widetilde{d}_{q}$ is a simple eigenvalue, i.e. the weighted degree
should be simple and only node $q$ in the graph $G$ has a degree equal to
$\widetilde{d}_{q}$. This assumption of unique degree is the most confining constraint of the presented method and expansions. Gerschgorin's theorem \cite[art. 245]
{PVM_graphspectra_second_edition} states that there is an eigenvalue
$0\leq\widetilde{\mu}_{q}\leq2\widetilde{d}_{q}$ of the weighted Laplacian
$\widetilde{Q}$ in an interval with length $\widetilde{d}_{q}$ around
$\widetilde{d}_{q}$. We will first apply the general perturbation in
\cite[Sec. 10.10]{PVM_graphspectra_second_edition} and compute the first four
terms in the eigenvalue series in (\ref{eigenvector_expansion_zeta}) analytically.

The first perturbation coefficient $c_{1}=x_{k}^{T}Bx_{q}$ in \cite[Sec.
10.10]{PVM_graphspectra_second_edition} translates to%
\[
\widetilde{c}_{1}\left(  q\right)  =e_{q}^{T}\widetilde{A}e_{q}=\widetilde
{a}_{qq}=0
\]
and, in general, $x_{k}^{T}Bx_{q}$ translates to $e_{k}^{T}\widetilde{A}%
e_{q}=\widetilde{a}_{kq}\geq0$. The second perturbation coefficient%
\begin{equation}
c_{2}=\sum_{k=1;k\neq q}^{N}\frac{\left(  x_{k}^{T}Bx_{q}\right)  ^{2}%
}{\lambda_{q}-\lambda_{k}} \label{coeff_c2_perturbation}%
\end{equation}
becomes%
\[
\widetilde{c}_{2}\left(  q\right)  =\sum_{k=1;k\neq q}^{N}\frac{\left(
\widetilde{a}_{kq}\right)  ^{2}}{\widetilde{d}_{q}-\widetilde{d}_{k}}%
\]
and in the unweighted case%
\[
c_{2}\left(  q\right)  =\sum_{k=1;k\neq q}^{N}\frac{a_{kq}}{d_{q}-d_{k}}%
=\sum_{k\in\mathcal{N}_{q}}\frac{1}{d_{q}-d_{k}}%
\]
where $\mathcal{N}_{q}$ is the set of all direct neighbors of node $q$. Thus,
the perturbation coefficient $c_{2}\left(  q\right)  $ is the sum of the
reciprocal degree difference with $d_{q}$ over all direct neighbors of node
$q$ and clearly illustrates that only node $q$ can have the degree $d_{q}$.
Using $\sum_{k=1;k\neq q}^{N}a_{kq}=d_{q}$, the degree of node $q$, we can
bound $c_{2}$ as%
\[
\left\vert c_{2}\left(  q\right)  \right\vert <\max_{\substack{1\leq k\leq
N\\k\neq q}}\frac{d_{q}}{\left\vert d_{q}-d_{k}\right\vert }=d_{q}=\left(
A^{2}\right)  _{qq}%
\]
The third perturbation coefficient in \cite[Sec.
10.10]{PVM_graphspectra_second_edition}
\begin{equation}
c_{3}=\sum_{r=1;r\neq q}^{N}\frac{x_{q}^{T}Bx_{r}}{\lambda_{q}-\lambda_{r}%
}\sum_{k=1;k\neq q}^{N}\frac{\left(  x_{k}^{T}Bx_{q}\right)  \left(  x_{k}%
^{T}Bx_{r}\right)  }{\lambda_{q}-\lambda_{k}}-\sum_{r=1;r\neq q}^{N}%
\frac{\left(  x_{r}^{T}Bx_{q}\right)  ^{2}\left(  x_{q}^{T}Bx_{q}\right)
}{\left(  \lambda_{q}-\lambda_{r}\right)  ^{2}} \label{coeff_c3_perturbation}%
\end{equation}
reduces, with $\widetilde{a}_{qq}=0$, to%
\[
\widetilde{c}_{3}\left(  q\right)  =\sum_{r=1;r\neq q}^{N}\frac{\widetilde
{a}_{rq}}{\widetilde{d}_{q}-\widetilde{d}_{r}}\sum_{k=1;k\neq q}^{N}%
\frac{\widetilde{a}_{kq}\widetilde{a}_{kr}}{\widetilde{d}_{q}-\widetilde
{d}_{k}}%
\]
and in the unweighted case%
\[
c_{3}\left(  q\right)  =\sum_{r=1;r\neq q}^{N}\sum_{k=1;k\neq q}^{N}%
\frac{a_{rq}}{\left(  d_{q}-d_{r}\right)  }\frac{a_{qk}}{\left(  d_{q}%
-d_{k}\right)  }a_{kr}=\sum_{r\in\mathcal{N}_{q}}\sum_{k\in\mathcal{N}_{q}%
\cap\mathcal{N}_{r}}\frac{1}{d_{q}-d_{r}}\frac{1}{d_{q}-d_{k}}%
\]
which is the sum of the product of reciprocal degree differences over all
mutually connected neighbors of node $q$. We bound $c_{3}$ as%
\begin{align*}
\left\vert c_{3}\left(  q\right)  \right\vert  &  =\left\vert \sum_{r=1;r\neq
q}^{N}\sum_{k=1;k\neq q}^{N}\frac{a_{rq}}{d_{q}-d_{r}}\frac{a_{qk}a_{kr}%
}{d_{q}-d_{k}}\right\vert \leq\sum_{r=1;r\neq q}^{N}\sum_{k=1;k\neq q}%
^{N}\frac{a_{rq}}{\left\vert d_{q}-d_{r}\right\vert }\frac{a_{qk}a_{kr}%
}{\left\vert d_{q}-d_{k}\right\vert }\\
&  <\max_{\substack{1\leq k\leq N\\k\neq q}}\frac{1}{\left\vert d_{q}%
-d_{k}\right\vert ^{2}}\sum_{r=1}^{N}a_{rq}\sum_{k=1}^{N}a_{qk}a_{kr}%
\end{align*}
With $\sum_{r=1}^{N}a_{rq}\sum_{k=1}^{N}a_{qk}a_{kr}=\sum_{r=1}^{N}%
a_{rq}\left(  A^{2}\right)  _{qr}=\left(  A^{3}\right)  _{qq}$, the number of
closed walks of length 3 from node $j$ back to itself
\cite{PVM_graphspectra_second_edition}, we find%
\[
\left\vert c_{3}\left(  q\right)  \right\vert <\max_{\substack{1\leq k\leq
N\\k\neq q}}\frac{\left(  A^{3}\right)  _{qq}}{\left\vert d_{q}-d_{k}%
\right\vert ^{2}}\leq\left(  A^{3}\right)  _{qq}%
\]
The fourth perturbation coefficient, not explicitly given in \cite[Sec.
10.10]{PVM_graphspectra_second_edition}, but derived from the recursion (\ref{coeff_ch(q)}) below,%
\begin{align}
c_{4}  &  =\left(  x_{q}^{T}Bx_{q}\right)  ^{2}\sum_{r=1;r\neq q}^{N}%
\frac{\left(  x_{r}^{T}Bx_{q}\right)  ^{2}}{\left(  \lambda_{q}-\lambda
_{r}\right)  ^{3}}-2\left(  x_{q}^{T}Bx_{q}\right)  \sum_{r=1;r\neq q}%
^{N}\frac{\left(  x_{q}^{T}Bx_{r}\right)  }{\left(  \lambda_{q}-\lambda
_{r}\right)  ^{2}}\sum_{k=1;k\neq q}^{N}\frac{\left(  x_{k}^{T}Bx_{q}\right)
\left(  x_{k}^{T}Bx_{r}\right)  }{\lambda_{q}-\lambda_{k}}\nonumber\\
&  \hspace{0.5cm}-\sum_{r=1;r\neq q}^{N}\frac{\left(  x_{r}^{T}Bx_{q}\right)
^{2}}{\left(  \lambda_{q}-\lambda_{r}\right)  ^{2}}\sum_{k=1;k\neq q}^{N}%
\frac{\left(  x_{k}^{T}Bx_{q}\right)  ^{2}}{\lambda_{q}-\lambda_{k}%
}\nonumber\\
&  \hspace{0.5cm}+\sum_{r=1;r\neq q}^{N}\frac{x_{q}^{T}Bx_{r}}{\lambda
_{q}-\lambda_{r}}\sum_{l=1;l\neq q}^{N}\frac{x_{r}^{T}Bx_{l}}{\lambda
_{q}-\lambda_{l}}\sum_{k=1;k\neq q}^{N}\frac{\left(  x_{k}^{T}Bx_{q}\right)
\left(  x_{k}^{T}Bx_{l}\right)  }{\lambda_{q}-\lambda_{k}}
\label{coeff_c4_perturbation}%
\end{align}
reduces, with $\widetilde{a}_{qq}=0$, to%
\[
\widetilde{c}_{4}\left(  q\right)  =\sum_{r=1;r\neq q}^{N}\frac{\widetilde
{a}_{rq}}{\widetilde{d}_{q}-\widetilde{d}_{r}}\sum_{l=1;l\neq q}^{N}%
\frac{\widetilde{a}_{rl}}{\widetilde{d}_{q}-\widetilde{d}_{l}}\sum_{k=1;k\neq
q}^{N}\frac{\widetilde{a}_{kq}\widetilde{a}_{kl}}{\widetilde{d}_{q}%
-\widetilde{d}_{k}}-\sum_{r=1;r\neq q}^{N}\frac{\left(  \widetilde{a}%
_{rq}\right)  ^{2}}{\left(  \widetilde{d}_{q}-\widetilde{d}_{r}\right)  ^{2}%
}\sum_{k=1;k\neq q}^{N}\frac{\left(  \widetilde{a}_{kq}\right)  ^{2}%
}{\widetilde{d}_{q}-\widetilde{d}_{k}}%
\]
and in the unweighted case
\[
c_{4}\left(  q\right)  =\sum_{r=1;r\neq q}^{N}\frac{a_{rq}}{d_{q}-d_{r}}%
\sum_{l=1;l\neq q}^{N}\frac{a_{rl}}{d_{q}-d_{l}}\sum_{k=1;k\neq q}^{N}%
\frac{a_{kq}a_{kl}}{d_{q}-d_{k}}-\sum_{r=1;r\neq q}^{N}\frac{a_{rq}}{\left(
d_{q}-d_{r}\right)  ^{2}}\sum_{k=1;k\neq q}^{N}\frac{a_{kq}}{d_{q}-d_{k}}%
\]
Conservatively bounding yields
\begin{align*}
\left\vert c_{4}\left(  q\right)  \right\vert  &  <\sum_{r=1;r\neq q}^{N}
\sum_{l=1;l\neq q}^{N}\sum_{k=1;k\neq q}^{N}\frac{a_{rq}}{\left\vert
d_{q}-d_{r}\right\vert }\frac{a_{rl}}{\left\vert d_{q}-d_{l}\right\vert }
\frac{a_{kq}a_{kl}}{\left\vert d_{q}-d_{k}\right\vert }
+ \sum_{r=1;r\neq q}^{N}\sum_{k=1;k\neq q}^{N}\frac{a_{rq}}{\left(
d_{q}-d_{r}\right)  ^{2}}\frac{a_{kq}}{\left\vert d_{q}-d_{k} \right\vert}\\
&  <\max_{\substack{1\leq k\leq N\\k\neq q}}\frac{1}{\left\vert d_{q}
-d_{k}\right\vert ^{3}}\sum_{r=1}^{N}\sum_{l=1}^{N}\sum_{k=1}^{N}a_{qr}
a_{rl}a_{lk}a_{kq}
+\max_{\substack{1\leq k\leq N\\k\neq q}}\frac{1}{\left\vert d_{q}%
-d_{k}\right\vert ^{3}}\sum_{r=1}^{N}\sum_{k=1}^{N}a_{rq}a_{kq}\\
& =\max_{\substack{1\leq k\leq N\\k\neq q}}\frac{\left(
A^{4}\right)_{qq}+\left( A^{2}\right)_{qq}}{\left\vert d_{q}-d_{k}\right\vert ^{3}}
\leq\left(A^{4}\right)_{qq}+\left( A^{2}\right)_{qq}%
\end{align*}
If $d_{q}=d_{\max}$, then $c_{2}\left(  q\right)  $ and $c_{3}\left(
q\right)  $ are positive, but the sign of $c_{4}\left(  q\right)  $ may be
negative. In summary, for the unweighted case and up to order $\zeta^{4}$ in the perturbation parameter $\zeta$, we
find that the eigenvalue expansion $\xi_{q}\left(  \zeta\right)  $ in
(\ref{eigenvalue_expansion_zeta}) of the matrix $\Delta+\zeta A$ around the
unique degree $d_{q}$ of node $q$ is
\begin{align}
\xi_{q}\left(  \zeta\right)   
&  =d_{q}+\zeta^{2}c_{2}\left(  q\right)
+\zeta^{3}c_{3}\left(  q\right)  +\zeta^{4}c_{4}\left(  q\right)  +O\left(
\zeta^{5}\right) \nonumber\\
&  =d_{q}+\zeta^{2}\sum_{k=1;k\neq q}^{N}\frac{a_{kq}}{d_{q}-d_{k}}+\zeta
^{3}\sum_{r=1;r\neq q}^{N}\frac{a_{rq}}{d_{q}-d_{r}}\sum_{k=1;k\neq q}
^{N}\frac{a_{qk}a_{kr}}{d_{q}-d_{k}}\nonumber\\
& +\zeta^{4}\left(  \sum_{r=1;r\neq q}^{N}\frac{a_{rq}}
{d_{q}-d_{r}}\sum_{l=1;l\neq q}^{N}\frac{a_{rl}}{d_{q}-d_{l}}\sum_{k=1;k\neq
q}^{N}\frac{a_{kq}a_{kl}}{d_{q}-d_{k}}
-\sum_{r=1;r\neq q}^{N}\frac{a_{rq}
}{\left(  d_{q}-d_{r}\right)  ^{2}}\sum_{k=1;k\neq q}^{N}\frac{a_{kq}}
{d_{q}-d_{k}}\right)  +O\left(  \zeta^{5}\right)
\label{lambda_zeta_upto_order4}
\end{align}

We may continue, as in the numerical evaluation sections below, to evaluate
higher orders in the perturbation parameter $\zeta$. For the matrix $W=$
$\widetilde{\Delta}$ and the perturbation matrix $B=$ $\widetilde{A}$, 
the perturbation Talyor series of a unique degree in (\ref{eigenvalue_expansion_zeta}) becomes
\begin{equation}
\label{Taylor_perturbation_series}
    \xi_q(\zeta) = d_q + \sum_{j=2}^{\infty}c_j(q)\zeta^{j}
\end{equation}
where the recursion of the coefficients in \cite[eq. (A.201-203) on p. 397]{PVM_graphspectra_second_edition} becomes
\begin{equation}
c_{j}\left(  q\right)  =\sum_{k=1;k\neq q}^{N}\beta_{j-1,k}\widetilde{a}%
_{qk}\hspace{0.5cm}\text{for }j>1\label{coeff_ch(q)}%
\end{equation}
in the eigenvalue expansion (\ref{eigenvalue_expansion_zeta}) and for
$r\neq q$,%
\[
\left\{
\begin{array}
[c]{lc}%
\beta_{1r}=\frac{\widetilde{a}_{rq}}{\widetilde{d}_{q}-\widetilde{d}_{r}} & \\
\beta_{jr}=\frac{1}{\widetilde{d}_{r}-\widetilde{d}_{q}}\sum_{l=1;l\neq q}%
^{N}\left\{  \sum_{k=1}^{j-2}\beta_{kr}\beta_{j-k-1,l}\widetilde{a}_{ql}%
-\beta_{j-1,l}\widetilde{a}_{rl}\right\}   & \text{for }j>1
\end{array}
\right.
\]
which can be computed up to any desired value of $j$. 
Both weighted and
unweighted cases have the same structure. 
The main difference lies in the
powers $\left(  \widetilde{a}_{rl}\right)  ^{m}$, which disappear in the zero-one
unweighted case, because $\left(  a_{rl}\right)  ^{m}=a_{rl}$. In addition,
all coefficients $c_{j}\left(  q\right)  \in\mathbb{Q}$ are rational numbers
in the unweighted case and computer computations with integers can be performed exactly up to any desired accuracy.
Therefore, in the sequel, we continue with the unweighted case, thus, omitting
the tildes, which simplifies the notation.

\subsection{Convergence of the perturbation Taylor series}
\label{sec_convergence_perturbation_Taylor_series}
Assuming that $\left\vert c_{j}\left(  q\right)  \right\vert \leq$
$\kappa^{j-1}\left(  A^{j}\right)  _{qq}$, where
\begin{equation}
    \kappa=\max_{\substack{1\leq
k\leq N\\k\neq q}}\frac{1}{\left\vert d_{q}-d_{k}\right\vert }=\frac{1}
{\min_{1\leq k\neq q\leq N}\left\vert d_{q}-d_{k}\right\vert }
\label{degree_difference_kappa}
\end{equation}
holds for any
integer $j$, then the eigenvalue expansion in (\ref{eigenvalue_expansion_zeta}%
) of the matrix $\Delta+\zeta A$,%
\begin{equation}
\xi_{q}\left(  \zeta\right)  =d_{q}+\zeta^{2}c_{2}\left(  q\right)  +\zeta
^{3}c_{3}\left(  q\right)  +\zeta^{4}c_{4}\left(  q\right)  +\cdots=d_{q}%
+\sum_{j=2}^{\infty}c_{j}\left(  q\right)  \zeta^{j}%
\end{equation}
is bounded as%
\[
\left\vert \xi_{q}\left(  \zeta\right)  \right\vert <d_{q}+\frac{1}{\kappa
}\sum_{j=2}^{\infty}\left(  A^{j}\right)  _{qq}\left(  \kappa\zeta\right)
^{j}=d_{q}+\frac{1}{\kappa}\left(  \sum_{j=2}^{\infty}A^{j}\left(  \kappa
\zeta\right)  ^{j}\right)  _{qq}%
\]
Introducing the eigenvalue decomposition $A=X\Lambda X^{T}=\sum_{k=1}%
^{N}\lambda_{k}\left(  A\right)  x_{k}x_{k}^{T}$, where $x_{k}$ is the
normalized eigenvector of the adjacency matrix $A$ belonging to eigenvalue
$\lambda_{k}\left(  A\right)  $ and assuming the ordering $\lambda_{1}\left(
A\right)  \geq\lambda_{2}\left(  A\right)  \geq\cdots\geq\lambda_{N}\left(
A\right)  $, then%
\[
\sum_{j=2}^{\infty}A^{j}\left(  \kappa\zeta\right)  ^{j}=\sum_{k=1}^{N}%
x_{k}x_{k}^{T}\sum_{j=2}^{\infty}\left(  \lambda_{k}\left(  A\right)
\kappa\zeta\right)  ^{j}%
\]
The geometric $j$-series converges, provided $\left\vert \lambda_{k}\left(
A\right)  \kappa\zeta\right\vert <1$ for any integer $k$, i.e. $\left\vert
\zeta\right\vert <\frac{1}{\kappa\lambda_{1}\left(  A\right)  }$. In the
unweighted case, $\kappa=\max_{\substack{1\leq k\leq N\\k\neq q}}\frac
{1}{\left\vert d_{q}-d_{k}\right\vert }\leq1$ and equality is reached if there
exists at least one node $k$ whose degree differs from the degree $d_{q}$ of
node $q$ by one; no degree $d_{k}$ of node $k\neq q$ can be equal to $d_{q}$
due to the multiplicity one condition of the perturbation expansion. In other
words, the larger the minimum difference $\left\vert d_{q}-d_{k}\right\vert $,
the smaller $\kappa$ and, consequently, the larger the convergence radius for
the perturbation series (\ref{eigenvalue_expansion_zeta}). The spectral radius
$\lambda_{1}\left(  A\right)  $ is bounded (see e.g.
\cite{PVM_graphspectra_second_edition}) by the average and maximum degree, $d_{av}=\frac{2L}{N}\leq\lambda
_{1}\left(  A\right)  \leq d_{\max}$, and the minimum difference $\max_{\substack{1\leq k\leq N}}\left\vert
d_{q}-d_{k}\right\vert <d_{av}$, which implies that $\kappa\lambda_{1}\left(
A\right)  >1$. In summary, if the bounds $\left\vert c_{j}\left(  q\right)
\right\vert \leq$ $\kappa^{j-1}\left(  A^{j}\right)  _{qq}$ holds for any
integer $j$, then the perturbation parameter $\zeta$ must be smaller than 1 in absolute value, i.e., $\left\vert \zeta\right\vert <1$ in order that the eigenvalue
expansion (\ref{eigenvalue_expansion_zeta}) of the matrix $\Delta+\zeta A$ as well as the matrix $\Delta-\zeta A$ converges. The bounds $\left\vert c_{j}\left(  q\right)  \right\vert \leq$
$\kappa^{j-1}\left(  A^{j}\right)  _{qq}$ are conservative, because, if
$d_{q}\neq d_{\max}$, then there will be negative terms and cancellations in
each sum of the coefficients $c_{j}\left(  q\right)  $ for $j\geq2$.
Nevertheless, numerical computations for $\zeta=-1$ indeed reveal that the
expansion $\xi_{q}\left(  -1\right)  =d_{q}+\sum_{j=2}^{\infty}\left(
-1\right)  ^{j}c_{j}\left(  q\right)  $, corresponding to an eigenvalue of the
Laplacian $Q=\Delta-A$ around degree $d_{q}$ (assumed to be unique or simple),
converges in some cases that we have considered.
For $N=10$ and $N=20$, we found numerically that the first 4 coefficients
$c_{k}\left(  q\right)  $ decrease and just from $c_{5}\left(  q\right)  $ on
increase. Limiting the series $\xi_{q;K}\left(  -1\right)  =d_{q}+\sum
_{j=2}^{K}\left(  -1\right)  ^{j}c_{j}\left(  q\right)  $ up to $K=4$ terms
seems \textquotedblleft reasonably\textquotedblright\ accurate\ and limiting a
divergent series (as in Stirling's approximation
\cite{PVM_Binet_convergentseries}) to the point, where the terms start to 
increase, may be adequate. 
\iffalse
Numerical evaluation seems that the accuracy of
$\xi_{q;K}\left(  -1\right)  $ improves with increasing size $N$ of the graph,
but for not too high link density $p=\frac{L}{\binom{N}{2}}$.
\fi

\subsection{Euler summation and the Euler series}
\label{sec_Euler_series}
In the sequel, we will only consider the perturbation parameter $\zeta = -1$ in the eigenvalue expansion $\xi_{q}\left(  \zeta\right)  $ in (\ref{eigenvalue_expansion_zeta}), but emphasize that the above convergence analysis holds for $\zeta=1$ as well, hence both for the Laplacian $Q=\Delta-A$ and the signless Laplacian $\overline{Q}=\Delta+A$.
In order to simplify the notation, we further write $\xi_{q} $ for $\xi_{q}\left( -1\right)  $.

Instead of summing $\xi_{q}  =d_{q}+\lim_{K\rightarrow\infty
}\sum_{j=2}^{K}c_{j}\left(  -1\right)  ^{j}$, Euler summation \cite{Hardy_div} 
\begin{equation}
 \xi_{q;K}  =d_{q}+\sum_{m=2}^{K}\left(  \sum_{k=2}^{m}%
\binom{m-1}{k-1}\left(  -1\right)  ^{k}c_{k}\right)  \frac{1}{2^{m}%
}\label{mu_perturbation_Eulersummation_K}%
\end{equation}
yields considerably better results, because Euler summation
\cite{Hardy_div} is known to extend the convergence range. If the Euler series (\ref{mu_perturbation_Eulersummation_K}) converges, then %
\begin{equation}
\mu_{k}=\lim_{K\rightarrow\infty}\xi_{q;K}  =d_{q}+\sum
_{m=2}^{\infty}\left(  \sum_{k=2}^{m}\binom{m-1}{k-1}\left(  -1\right)
^{k}c_{k}\right)  \frac{1}{2^{m}}\label{mu_perturbation_Eulersummation}%
\end{equation}
Thus, in case of convergence, the sum (\ref{mu_perturbation_Eulersummation_K}) for sufficiently large $K$ tends
to the $k$-th Laplacian eigenvalue $\mu_{k}$ close to
$d_{q}$ for some integer $k$, where the eigenvalues of the Laplacian matrix $Q$ are $\mu_{1}\geq
\mu_{2}\geq\cdots\geq\mu_{N}=0$. 
The Euler series (\ref{mu_perturbation_Eulersummation_K}) demonstrates that the worst-case complexity $C=O(K^2N^2)$ for  large $N$. Indeed, the most consuming part lies in the computation of $\beta_{jr}$ in the recursion (\ref{coeff_ch(q)}). The computational complexity of the recursion (\ref{coeff_ch(q)}) for each $\beta_{jr}$  for $1 \leq r \leq N$ and $1 \leq j \leq K$ is of the order $O(KN)$.

For $K=4$ terms in (\ref{mu_perturbation_Eulersummation_K}) and if the node
$q$ has a sufficiently large degree and unique $d_{q}$, then%
\begin{align}
\xi_{q;4}   &  \approx d_{q}+\frac{11}{16}c_{2}-\frac{5}%
{16}c_{3}+\frac{1}{16}c_{4}\nonumber\\
&  =d_{q}+\frac{11}{16}\sum_{k=1;k\neq q}^{N}\frac{a_{kq}}{d_{q}-d_{k}}%
-\frac{5}{16}\sum_{r=1;r\neq q}^{N}\frac{a_{rq}}{d_{q}-d_{r}}\sum_{k=1;k\neq
q}^{N}\frac{a_{qk}a_{kr}}{d_{q}-d_{k}}\nonumber\\
&  \hspace{0.5cm}+\frac{1}{16}\left(  \sum_{r=1;r\neq q}^{N}\frac{a_{rq}%
}{d_{q}-d_{r}}\sum_{l=1;l\neq q}^{N}\frac{a_{rl}}{d_{q}-d_{l}}\sum_{k=1;k\neq
q}^{N}\frac{a_{kq}a_{kl}}{d_{q}-d_{k}}-\sum_{r=1;r\neq q}^{N}\frac{a_{rq}%
}{\left(  d_{q}-d_{r}\right)  ^{2}}\sum_{k=1;k\neq q}^{N}\frac{a_{kq}}%
{d_{q}-d_{k}}\right)  \label{lambda_Euler_upto_order4}%
\end{align}
is a reasonable estimate for a Laplacian eigenvalue $\mu_{k}\approx\xi
_{q;4}  $. As illustrated in Section \ref{sec_performance_analysis} below, we found numerically that the
largest eigenvalue $\mu_{1}$ is retrieved from $d_{q}=d_{\max}$, even with 4
coefficients, quite accurately! In addition, Euler summation
(\ref{mu_perturbation_Eulersummation}) also seems to converge for other large
degrees or small $\kappa=\max_{\substack{1\leq k\leq N\\k\neq q}}\frac
{1}{\left\vert d_{q}-d_{k}\right\vert }$. On the other hand, expansion around
$d_{q}=d_{\min}$ is considerably less accurate and Euler summation
(\ref{mu_perturbation_Eulersummation}) does not seem to converge anymore.

\subsection{Extension of the Euler series and analogy with Lagrange series}
\label{sec_analogy_Lagrange}
The more general Euler transform of a Taylor series $f(z)=f_{0}+\sum_{m=1}^{\infty}f_{k}z^{k}$
disposes of a tune-able parameter $t$,%
\begin{equation}
f(z)=f_{0}+\sum_{m=1}^{\infty}\left[  \sum_{k=1}^{m}{\binom{m-1}{k-1}}%
\,f_{k}\,t^{m-k}\right]  \;\left(  \frac{z}{1+tz}\right)  ^{m}%
\label{Eulertransform}%
\end{equation}
Thus, (\ref{mu_perturbation_Eulersummation}) generalizes to an Euler $t$-series,%
\begin{equation}
\xi_{q;K}  =d_{q}+\lim_{K\rightarrow\infty}\sum_{m=2}%
^{K}\left(  \sum_{k=2}^{m}\binom{m-1}{k-1}t^{m-k}c_{k}\right)  \left(
\frac{1}{t-1}\right)  ^{m}\label{mu_perturbation_Eulersummation_general_t}%
\end{equation}
and (\ref{mu_perturbation_Eulersummation_general_t}) reduces to (\ref{mu_perturbation_Eulersummation}) for $t=-1$.  The Euler
summation (\ref{mu_perturbation_Eulersummation_general_t}) of the matrix
perturbation series in (\ref{Taylor_perturbation_series}) bears
resemblance to a Lagrange series, where expansion around different points may
converge to a same zero. Lagrange's series for the inverse $f^{-1}(z)$ of a
function $f\left(  z\right)  $ is \cite[II, pp. 88]{Markushevich}
\begin{equation}
f^{-1}(z)=z_{0}+\sum_{m=1}^{\infty}\frac{1}{m!}\left.  \left[  \frac{d^{m-1}%
}{dw^{m-1}}\left(  \frac{w-z_{0}}{f(w)-f(z_{0})}\right)  ^{m}\right]
\right\vert _{w=z_{0}}(z-f(z_{0}))^{m}\label{Lagrange}%
\end{equation}
A zero $y$ of $f\left(  z\right)  $, obeying $f\left(  y\right)  =0$ and
$y=f^{-1}(0)$, has the Lagrange series%
\[
y=f^{-1}(0)=z_{0}+\sum_{m=1}^{\infty}\frac{1}{m!}\left.  \left[  \frac
{d^{m-1}}{dw^{m-1}}\left(  \frac{w-z_{0}}{f(w)-f(z_{0})}\right)  ^{m}\right]
\right\vert _{w=z_{0}}(-f(z_{0}))^{m}%
\]
provided that $z_{0}$ is sufficiently close to the zero $y$, else the Lagrange
series may diverge or converge towards a different, more nearby zero of
$f\left(  z\right)  $. Since the characteristic polynomial $c_{Q}\left(
z\right)  =\det\left(  Q-zI\right)  =\sum_{k=0}^{N}\gamma_{k}z^{k}=\prod
_{k=1}^{N}\left(  \mu_{k}-z\right)  $ corresponds to the function $f\left(
z\right)  $, it is tempting to infer by comparing the Euler summation
(\ref{mu_perturbation_Eulersummation_general_t}) and the Lagrange series
(\ref{Lagrange}) of $\mu_{k}=c_{Q}^{-1}\left(  0\right)  $ that $z_{0}=d_{q}$,
$c_{Q}\left(  d_{q}\right)  =\frac{-1}{t-1}$ and $\sum_{k=2}^{m}\binom
{m-1}{k-1}t^{m-k}c_{k}=\frac{1}{m!}\left.  \left[  \frac{d^{m-1}}{dw^{m-1}%
}\left(  \frac{w-d_{q}}{c_{Q}\left(  w\right)  -c_{Q}\left(  d_{q}\right)
}\right)  ^{m}\right]  \right\vert _{w=d_{q}}$. However, we are unable to
prove this speculation, although it is correct for the special case in Sec. \ref{sec:almost_regular}.
%%%%%%%%%%%%%%%%%%%%%%%%%%%%%%%%%%%%%%%%%%%%%%%%%%%%%%%%%%%%%%

\section{Performance evaluation of the Euler series $\xi_{q;K}$}
\label{sec_performance_analysis}
We assess the convergence of the Euler series $\xi_{q;K}$ in (\ref{mu_perturbation_Eulersummation_general_t}), whose special case for $t=-1$ is in (\ref{mu_perturbation_Eulersummation}).
We can provide some insights by numerical computation of various cases in the following subsections.

\subsection{Impact of the tuning parameter $t$ on the convergence of the Euler series $\xi_{q;K} $}

To evaluate how good the Euler $t$-series $\xi_{q ; K}$ in (\ref{mu_perturbation_Eulersummation_general_t}) approximates the Laplacian eigenvalue $\mu_q$ for different tuning parameters $t$, we define the difference 
\begin{equation*}
    10^{\alpha} : = \left| \xi_{q ; K}(-1) - \mu_{m} \right|
\end{equation*}
The Euler $t$-series $\xi_{q ; K}$ converges to an eigenvalue $\mu_m$ if the accuracy $\alpha$ becomes increasingly negative as the number of terms $K$ increases.

\textbf{Example 1.} A tree on $N=5$ nodes, with the adjacency matrix
\begin{equation}\label{eq:adj_tree}
A=\left[
\begin{array}
[c]{ccccc}%
0 & 0 & 1 & 1 & 1\\
0 & 0 & 0 & 0 & 1\\
1 & 0 & 0 & 0 & 0\\
1 & 0 & 0 & 0 & 0\\
1 & 1 & 0 & 0 & 0
\end{array}
\right]
\end{equation}
has a degree vector $d=\left(  3,1,1,1,2 \right)  $ and Laplacian eigenvalue
vector $\mu=\left(  4.17009,2.31111,1.,0.518806,0 \right)  $. The infinite
series (\ref{mu_perturbation_Eulersummation}) does not seem to converge, but
the Euler sum (\ref{lambda_Euler_upto_order4}) up to $K=4$ terms equals
$\xi_{5;4}\left(  -1\right)  =2.125$ for $d_{5}=2$ and $\xi_{5;5}\left(
-1\right)  =2.375$ for $K=5$ terms and $\xi_{1;4}\left(  -1\right)  =4.21875$
for $d_{1}=3$. Also, all odd coefficients $c_{2m+1}(q)=0$ in
(\ref{coeff_ch(q)}) for odd $m\geq0$ are zero, possibly agreeing with the fact
\cite[p. 217]{PVM_graphspectra_second_edition} that the characteristic
polynomial of the adjacency matrix of a tree is even, $c_{A_{\text{tree}}%
}\left(  z\right)  =\det\left(  A_{\text{tree}}-zI\right)  =c_{A_{\text{tree}%
}}\left(  -z\right)  $.
\begin{figure}[htbp]
    \centering
    \subfigure[$d_1=3$, $\mu_1=4.17009$]{ \label{fig:e11} \includegraphics[width=0.47\textwidth]{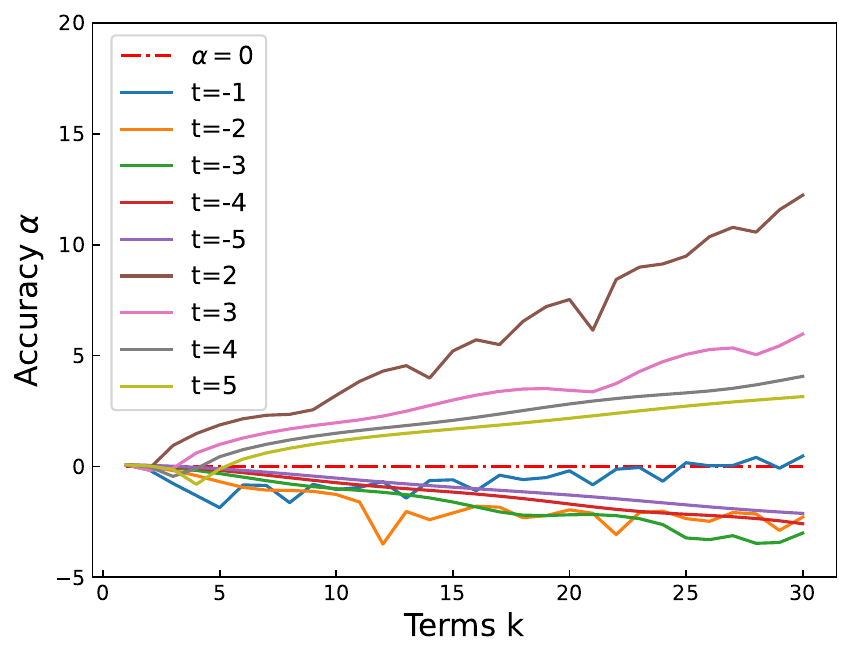} }
    \subfigure[$d_5=2$, $\mu_2=2.31111$]{ \label{fig:e12} \includegraphics[width=0.47\textwidth]{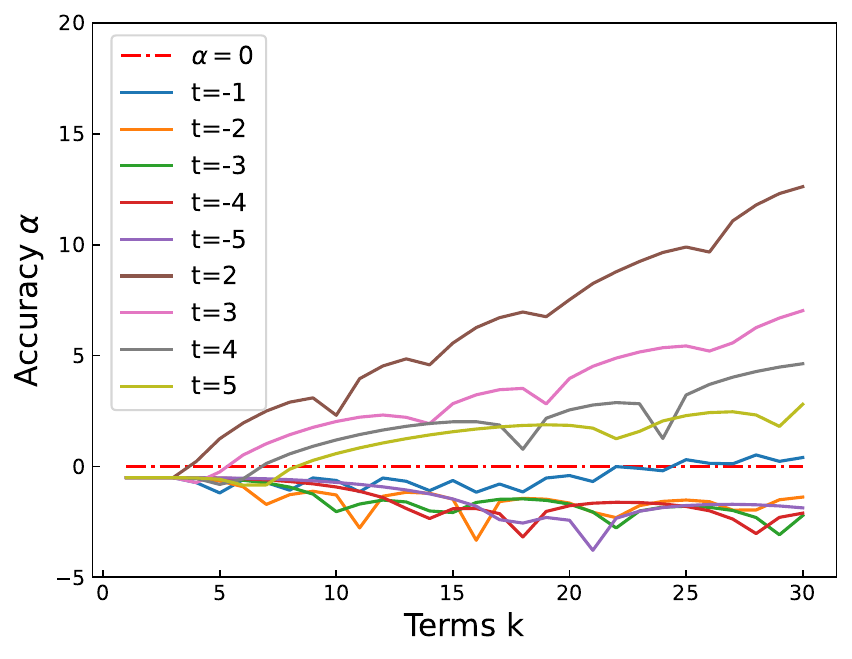} }
    \caption{The accuracy $\alpha$ of the Euler $t$-series $\xi_{q; K}$ versus the number of terms $K$ and tuning parameter $t$ for two nodes with a unique degree in Example 1.}
    \label{fig:e1}
\end{figure}

Figure \ref{fig:e1} shows the influence of the tuning parameter $t$ on the convergence of the Euler $t$-series $\xi_{q; K}$ for Example 1.
The dash-dot line $\alpha=0$ represents the eigenvalue $\mu_m$. 
The tree with adjacency matrix (\ref{eq:adj_tree}) has two nodes with a unique degree $d_1=3$ and $d_5=2$. 
As the number of terms $K$ increases, the Euler $t$-series $\xi_{1; K}$ with a tuning parameter $t=-1$ does not converge to the Laplacian eigenvalue $\mu_1$, but the more general Euler $t$-series with a tuning parameter $t<-1$ converges to the  Laplacian eigenvalue $\mu_1$.
For another unique degree $d_5=2$, the Euler series $\xi_{5; K}$ exhibits a similar performance.
Both Euler $t$-series $\xi_{1; K}$ and $\xi_{5; K}$ converge faster as the negative tuning parameter $t$ decreases.
Moreover, Fig. \ref{fig:e11} and \ref{fig:e12} indicate that the Euler $t$-series $\xi_{q; K}$ with a positive parameter $t$ always diverges but diverges slower for a larger positive parameter $t$.
In summary, Example 1 illustrates that the convergence of the Euler $t$-series $\xi_{q; K}$ in (\ref{mu_perturbation_Eulersummation_general_t}) can be improved by negative tuning parameters $t$.

\textbf{Example 2.} An instance of an Erd\H{o}s-R\'{e}nyi graph $G_{p}\left(
N\right)  $ on $N=20$ nodes and link density $p=0.3$ has the adjacency matrix
\[
A=\left[
\begin{array}
[c]{cccccccccccccccccccc}%
0 & 0 & 1 & 0 & 1 & 0 & 1 & 0 & 0 & 0 & 0 & 0 & 0 & 1 & 0 & 0 & 0 & 0 & 0 &
0\\
0 & 0 & 0 & 0 & 1 & 0 & 0 & 0 & 0 & 0 & 1 & 0 & 0 & 0 & 1 & 1 & 0 & 0 & 0 &
0\\
1 & 0 & 0 & 0 & 1 & 0 & 0 & 0 & 1 & 1 & 0 & 0 & 1 & 1 & 1 & 0 & 1 & 0 & 0 &
0\\
0 & 0 & 0 & 0 & 0 & 0 & 0 & 0 & 0 & 0 & 0 & 1 & 1 & 0 & 0 & 1 & 0 & 0 & 0 &
0\\
1 & 1 & 1 & 0 & 0 & 1 & 0 & 0 & 1 & 1 & 0 & 0 & 0 & 0 & 0 & 0 & 0 & 1 & 0 &
0\\
0 & 0 & 0 & 0 & 1 & 0 & 1 & 0 & 0 & 1 & 0 & 0 & 1 & 0 & 0 & 0 & 0 & 0 & 0 &
0\\
1 & 0 & 0 & 0 & 0 & 1 & 0 & 0 & 1 & 0 & 1 & 1 & 1 & 1 & 1 & 0 & 1 & 1 & 1 &
1\\
0 & 0 & 0 & 0 & 0 & 0 & 0 & 0 & 1 & 0 & 1 & 0 & 0 & 0 & 0 & 0 & 0 & 0 & 1 &
1\\
0 & 0 & 1 & 0 & 1 & 0 & 1 & 1 & 0 & 0 & 1 & 0 & 1 & 0 & 0 & 0 & 0 & 0 & 1 &
0\\
0 & 0 & 1 & 0 & 1 & 1 & 0 & 0 & 0 & 0 & 1 & 0 & 0 & 1 & 0 & 0 & 0 & 0 & 0 &
1\\
0 & 1 & 0 & 0 & 0 & 0 & 1 & 1 & 1 & 1 & 0 & 0 & 1 & 0 & 1 & 0 & 0 & 0 & 0 &
0\\
0 & 0 & 0 & 1 & 0 & 0 & 1 & 0 & 0 & 0 & 0 & 0 & 0 & 1 & 0 & 1 & 1 & 0 & 0 &
1\\
0 & 0 & 1 & 1 & 0 & 1 & 1 & 0 & 1 & 0 & 1 & 0 & 0 & 1 & 0 & 1 & 1 & 1 & 0 &
0\\
1 & 0 & 1 & 0 & 0 & 0 & 1 & 0 & 0 & 1 & 0 & 1 & 1 & 0 & 0 & 0 & 0 & 0 & 0 &
0\\
0 & 1 & 1 & 0 & 0 & 0 & 1 & 0 & 0 & 0 & 1 & 0 & 0 & 0 & 0 & 1 & 0 & 1 & 0 &
0\\
0 & 1 & 0 & 1 & 0 & 0 & 0 & 0 & 0 & 0 & 0 & 1 & 1 & 0 & 1 & 0 & 0 & 0 & 1 &
0\\
0 & 0 & 1 & 0 & 0 & 0 & 1 & 0 & 0 & 0 & 0 & 1 & 1 & 0 & 0 & 0 & 0 & 1 & 0 &
1\\
0 & 0 & 0 & 0 & 1 & 0 & 1 & 0 & 0 & 0 & 0 & 0 & 1 & 0 & 1 & 0 & 1 & 0 & 0 &
0\\
0 & 0 & 0 & 0 & 0 & 0 & 1 & 1 & 1 & 0 & 0 & 0 & 0 & 0 & 0 & 1 & 0 & 0 & 0 &
0\\
0 & 0 & 0 & 0 & 0 & 0 & 1 & 1 & 0 & 1 & 0 & 1 & 0 & 0 & 0 & 0 & 1 & 0 & 0 & 0
\end{array}
\right]
\]
The corresponding degree vector is
\[
d=\left(  4,4,8,3,7,4,12,4,7,6,7,6,10,6,6,6,6,5,4,5\right)
\]
ranked in decreasing order as $\left(
12,10,8,7,7,7,6,6,6,6,6,6,5,5,4,4,4,4,4,3\right)  $ to see uniqueness. The
Laplacian eigenvalue vector is
\begin{align*}
\mu &  =\left(
13.3514,11.6199,9.80641,9.32872,7.6586,7.46193,7.11613,6.92149,\right.  \\
&  6.3782,6.07484,5.80058,5.29648,4.59557,4.05486,3.58036,3.50647\\
&  \left.  2.83079,2.39082,2.22645,0\right)
\end{align*}
The largest eigenvalue of the adjacency matrix $A$ is $\lambda_{1}\left(
A\right)  =6.67615$. The Euler summation $\xi_{q;K} $ in
(\ref{mu_perturbation_Eulersummation}) converges for node $q=13$ and $d_{q}=10$, for which $\kappa=\frac{1}{2}$.
Indeed, $\xi_{13;K}$ as a function of the number $K$ of terms converges to
$\mu_{2}=11.6199127895910$ (with 15 digits accurate) as%
\[%
\begin{array}
[c]{cccc}%
\xi_{13;2}=10.48154762 & \xi_{13;7}=11.61206362 & \xi_{13;12}=11.61699285 &
\xi_{13;17}=11.62009681\\
\xi_{13;3}=11.00138889 & \xi_{13;8}=11.62002740 & \xi_{13;13}=11.61921713 &
\xi_{13;18}=11.61968541\\
\xi_{13;4}=11.33195709 & \xi_{13;9}=11.61508019 & \xi_{13;14}=11.62029217 &
\xi_{13;19}=11.61958213\\
\xi_{13;5}=11.49508126 & \xi_{13;10}=11.61181587 & \xi_{13;15}=11.62070805 &
\xi_{13;20}=11.61970380\\
\xi_{13;6}=11.57496760 & \xi_{13;11}=11.61364728 & \xi_{13;16}=11.62057580 &
\xi_{13;30}=11.61991367
\end{array}
\]
Furthermore, for $K$ up to 100, the table gives the approximation up to 15
digits and last column specifies the difference $\mu_{2}-\xi_{13;K}$, as an
indication of the accuracy,%

\[%
\begin{array}
[c]{ccc}%
\xi_{13;10} & 11.6118158710925 & 0.00809692\\
\xi_{13;20} & 11.6197037971111 & 0.000208992\\
\xi_{13;30} & 11.6199136700045 & -\text{8.8041349$\;10^{-7}$}\\
\xi_{13;40} & 11.6199135474613 & -\text{7.5787030$\;10^{-7}$}\\
\xi_{13;50} & 11.6199128258523 & -\text{3.6261285$\;10^{-8}$}\\
\xi_{13;60} & 11.6199127874211 & \text{2.1699282$\;10^{-9}$}\\
\xi_{13;70} & 11.6199127892558 & \text{3.3520741$\;10^{-10}$}\\
\xi_{13;80} & 11.6199127895867 & \text{4.3485215$\;10^{-12}$}\\
\xi_{13;90} & 11.6199127895931 & -\text{2.0765611$\;10^{-12}$}\\
\xi_{13;100} & 11.6199127895912 & -\text{1.5099033$\;10^{-13}$}%
\end{array}
\]

The Euler summation $\xi_{q;K}$ converges faster for node $q=7$ with the
maximum degree $d_{q}=12$. Indeed, $\xi_{7;K}$ as function of $K$ converges to
$\mu_{1}=13.3513926733482839961128497270$ (accurate up to 30 digits) as%
\[%
\begin{array}
[c]{cccc}%
\xi_{7;2}=12.55684524 & \xi_{7;7}=13.32451888 & \xi_{7;12}=13.35071509 &
\xi_{7;17}=13.35134956\\
\xi_{7;3}=12.92105159 & \xi_{7;8}=13.33893469 & \xi_{7;13}=13.35094032 &
\xi_{7;18}=13.35137642\\
\xi_{7;4}=13.10777862 & \xi_{7;9}=13.34549561 & \xi_{7;14}=13.35114508 &
\xi_{7;19}=13.35138894\\
\xi_{7;5}=13.22029144 & \xi_{7;10}=13.34889571 & \xi_{7;15}=13.35126516 &
\xi_{7;20}=13.35139125\\
\xi_{7;6}=13.28981543 & \xi_{7;11}=13.35029701 & \xi_{7;16}=13.35131598 &
\xi_{7;30}=13.35139267
\end{array}
\]
where all presented digits of $\xi_{7;30}$ are correct. Furthermore, the table
up to $K=100$ with the difference $\mu_{1}-\xi_{7;K}$ in the last column gives
is
\[%
\begin{array}
[c]{ccc}%
\xi_{7;10} & 13.3488957090710433527956303093 & 0.00249696\\
\xi_{7;20} & 13.3513912536915562912298783841 & \text{1.41966\;$10^{-6}$}\\
\xi_{7;30} & 13.3513926692587838827462033968 & \text{4.08949\;$10^{-9}$}\\
\xi_{7;40} & 13.3513926733102103638002761686 & \text{3.80691\;$10^{-11}$}\\
\xi_{7;50} & 13.3513926733476992672999617420 & \text{5.79092\;$10^{-13}$}\\
\xi_{7;60} & 13.3513926733482839312804648760 & \text{6.48323\;$10^{-17}$}\\
\xi_{7;70} & 13.3513926733482839615285740731 & \text{3.45842\;$10^{-17}$}\\
\xi_{7;80} & 13.3513926733482839965575848138 & \text{-4.4473\;$10^{-19}$}\\
\xi_{7;90} & 13.3513926733482839961083203713 & \text{4.52935\;$10^{-21}$}\\
\xi_{7;100} & 13.3513926733482839961129243912 & \text{-7.4664234\;$10^{-23}$}%
\end{array}
\]

For node $q=3$ with degree $d_{3}=8$, the Euler summation
(\ref{mu_perturbation_Eulersummation}) seems to diverge. Indeed, $\xi_{3;K}$
initially tends to converge to $\mu_{3}=9.80641$, but diverges for larger $K$,%
\[%
\begin{array}
[c]{cccc}%
\xi_{3;2}=8.937500000 & \xi_{3;7}=9.961090970 & \xi_{3;12}=7.834417220 &
\xi_{3;17}=20.15673862\\
\xi_{3;3}=9.593750000 & \xi_{3;8}=9.152494535 & \xi_{3;13}=11.10756619 &
\xi_{3;18}=19.02124175\\
\xi_{3;4}=9.541536458 & \xi_{3;9}=9.649234801 & \xi_{3;14}=13.44103998 &
\xi_{3;19}=-15.77306388\\
\xi_{3;5}=9.632552083 & \xi_{3;10}=10.87231670 & \xi_{3;15}=5.881312597 &
\xi_{3;20}=-0.7480476187\\
\xi_{3;6}=10.14137146 & \xi_{3;11}=9.574716849 & \xi_{3;16}=3.636664041 &
\xi_{3;30}=-1883.697136
\end{array}
\]
All mentioned values of $\xi_{q;4}$, corresponding to the explicitly form in
(\ref{lambda_Euler_upto_order4}), indicate that
(\ref{lambda_Euler_upto_order4}) is a reasonably accurate estimate for a
Laplacian eigenvalue.

\begin{figure}[htbp]
    \centering
    \subfigure[$d_7=12$, $\mu_1=13.3514$]{\label{fig:e21} \includegraphics[width=0.46\textwidth]{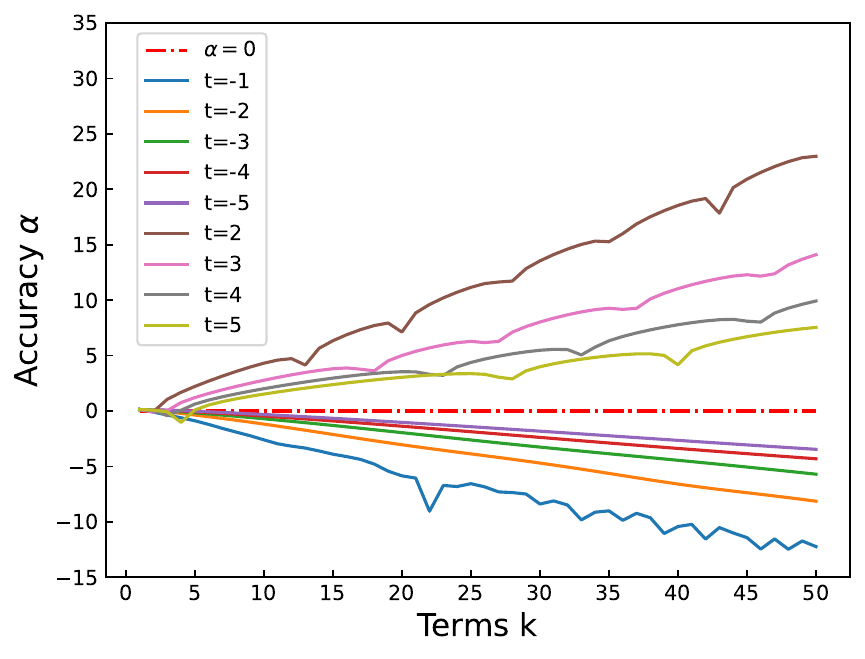}} 
    \subfigure[$d_{13}=10$, $\mu_2=11.6199$]{\label{fig:e22} \includegraphics[width=0.46\textwidth]{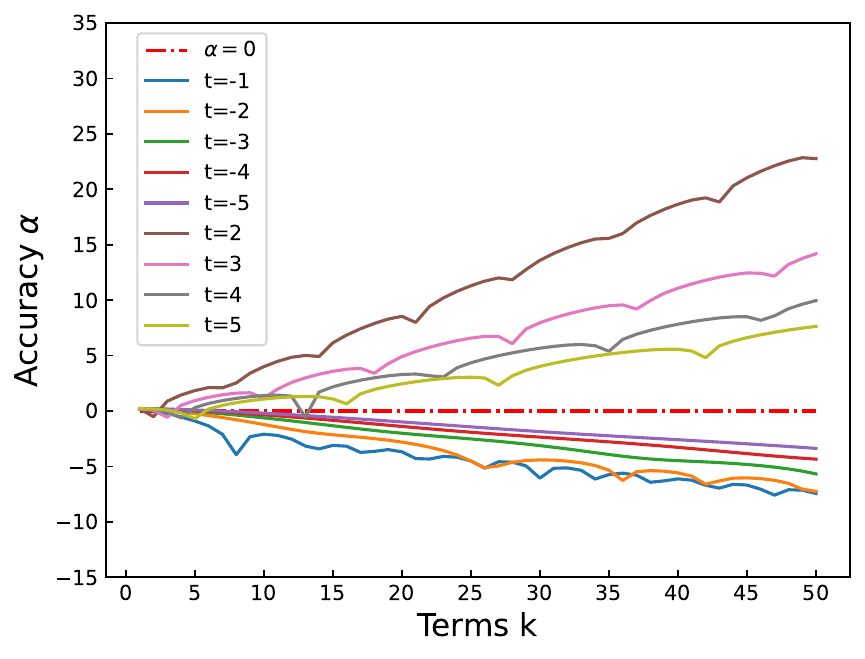}}
    \subfigure[$d_3=8$, $\mu_3=9.80641$]{\label{fig:e23} \includegraphics[width=0.46\textwidth]{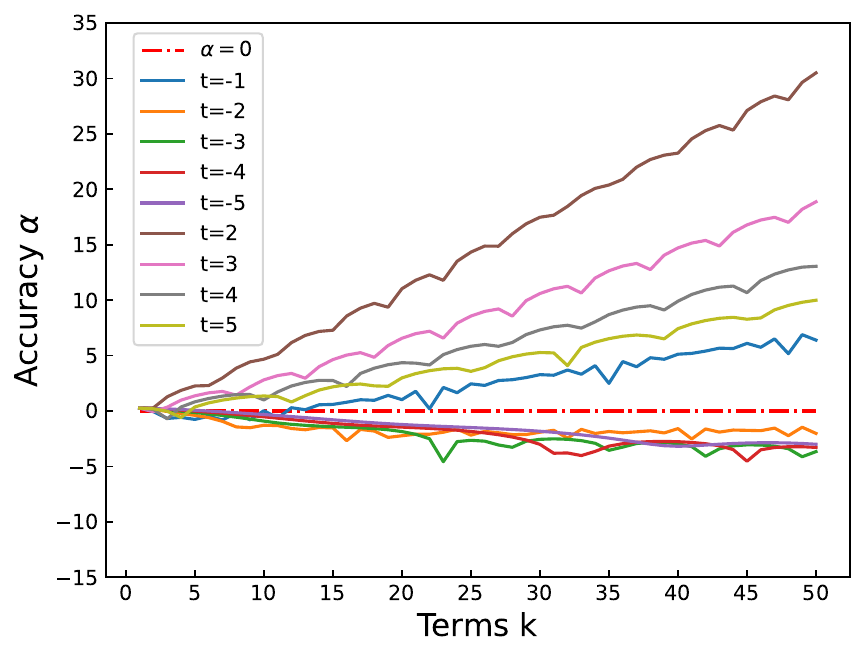}} 
    \subfigure[$d_4=3$, $\mu_{20}=0$]{\label{fig:e24} \includegraphics[width=0.46\textwidth]{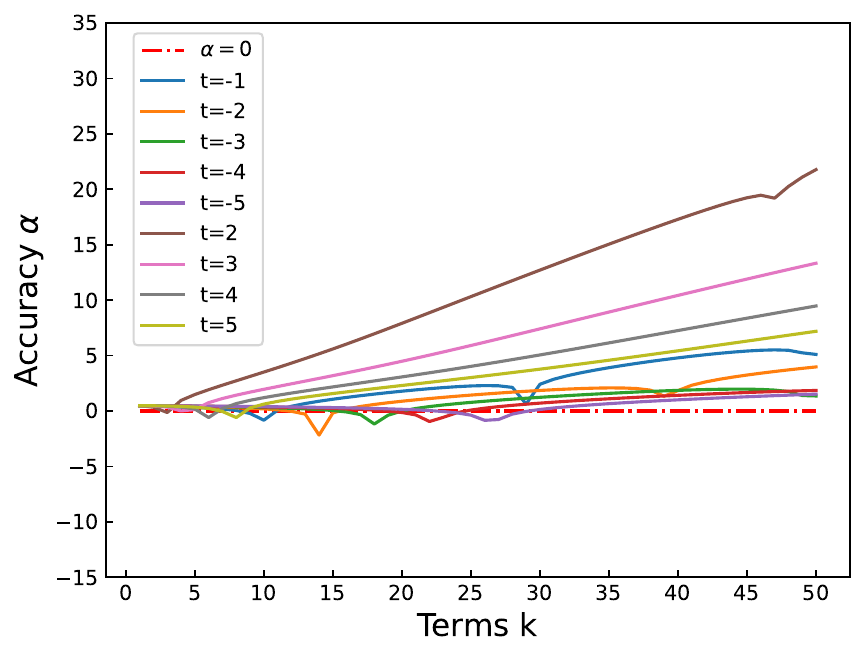}}
    \caption{The accuracy $\alpha$ of the Euler $t$-series $\xi_{q; K}$  versus the number of terms $K$ and tuning parameter $t$ for four nodes with a unique degree in Example 2.}
    \label{fig:e2}
\end{figure}

The convergence of the Euler $t$-series $\xi_{q; K}$ for the four nodes $q=3,4,7$ and $12$ with a unique degree is shown in Fig. \ref{fig:e2} for Example 2. 
We observe that the Euler $t$-series $\xi_{7; K}$ and  $\xi_{13; K}$ in (\ref{mu_perturbation_Eulersummation_general_t}) converge to the Laplacian eigenvalue $\mu_1$ and $\mu_2$ respectively when the tuning parameter $t$ is negative. 
In addition, Figure \ref{fig:e23} shows that the Euler $t$-series $\xi_{3; K}$ diverges with $t=-1$ but converges with $t<-1$, which agrees with our observation in Example 1 that the convergence of the Euler $t$-series $\xi_{q; K}$ can be improved by the tuning parameter $t$.
For the node $q=4$ with the smallest and unique degree $d_4=3$, Figure \ref{fig:e24} depicts that the Euler series $\xi_{4; K}$ diverges, no matter whether the tuning parameter $t$ is.
The Euler $t$-series $\xi_{q; K}$ with a positive tuning parameter $t$  does not converge for the four nodes in Example 2, similar to what we observe in Example 1.
The example implies that for a node with a high and unique degree, the Euler $t$-series $\xi_{q; K}$ in (\ref{mu_perturbation_Eulersummation_general_t}) probably converges to a Laplacian eigenvalue with a suitable parameter $t$, but the convergence hardly happens for a node with a small and unique degree.

The two examples provide us with an understanding of the impact of the tuning parameter $t$ on the convergence of the Euler $t$-series $\xi_{q; K}$. 
First, the convergence of the Euler series $\xi_{q; K}$ is related to the parameter $t$.
With a negative tuning parameter $t$, the Euler $t$-series $\xi_{q; K}$ likely converges to a Laplacian eigenvalue for nodes with large and unique degrees.
If the Euler series $\xi_{q; K}$ converges, then the speed of convergence changes when parameter $t$ changes, for example, node $q=7$ with $d_7=12$ and node $q=13$ with $d_{13}=10$ in Example 2.
However, the Euler $t$-series $\xi_{q; K}$  hardly converges for any parameter $t$ for some nodes with a small and unique degree, for instance, node $q=4$ with a degree $d_4=3$ in Example 2. 
Moreover, the convergence of the Euler $t$-series $\xi_{q; K}$ can hardly happen with a positive parameter $t$.

\subsection{Impact of the node degree on the convergence of the Euler series $\xi_{q;K} $ }
We have shown that the convergence of the Euler $t$-series $\xi_{q; K}$ of a node with a large and unique degree performs differently from a node with a small and unique degree, which motivates us to investigate the role of the node degree itself in the convergence of the Euler series $\xi_{q;K}$.
We confine ourselves to antiregular graphs, where 
each node, except for two nodes, has a different degree.  
Thus, antiregular graphs are reasonable examples to examine the relationship between the convergence of the Euler series  $\xi_{q;K}$ and the node degree. Moreover, antiregular graphs \cite[Sec.3.11]{brouwer2011spectra} possess integer Laplacian eigenvalues, which are composed of 0 and $N$ and all unique node degrees for an antiregular graph with $N$ nodes.

\textbf{Example 3.} An antiregular  graph $G$ of $N=10$ nodes has the adjacency matrix
\begin{equation}
    A=\left[\begin{array}{llllllllll}
0 & 1 & 0 & 1 & 0 & 1 & 0 & 1 & 0 & 1 \\
1 & 0 & 0 & 1 & 0 & 1 & 0 & 1 & 0 & 1 \\
0 & 0 & 0 & 1 & 0 & 1 & 0 & 1 & 0 & 1 \\
1 & 1 & 1 & 0 & 0 & 1 & 0 & 1 & 0 & 1 \\
0 & 0 & 0 & 0 & 0 & 1 & 0 & 1 & 0 & 1 \\
1 & 1 & 1 & 1 & 1 & 0 & 0 & 1 & 0 & 1 \\
0 & 0 & 0 & 0 & 0 & 0 & 0 & 1 & 0 & 1 \\
1 & 1 & 1 & 1 & 1 & 1 & 1 & 0 & 0 & 1 \\
0 & 0 & 0 & 0 & 0 & 0 & 0 & 0 & 0 & 1 \\
1 & 1 & 1 & 1 & 1 & 1 & 1 & 1 & 1 & 0
\end{array}\right]
\end{equation}
The degree vector 
\begin{equation}\label{eq:d}
    d = (5,5,4, 6, 3, 7, 2, 8, 1, 9)
\end{equation}
ranks in decreasing order as
\begin{equation*}
    d = (9,8,7,6,5,5,4,3,2,1).
\end{equation*}
The Laplacian eigenvalue vector is
\begin{equation*}
    \mu=(10,9,8,7,6,4,3,2,1,0)
\end{equation*}

\begin{figure}[htbp]
    \centering
    \subfigure[$d_9=1$, $\mu_{10}=0$]{\label{fig:ea1} \includegraphics[width=0.4\textwidth]{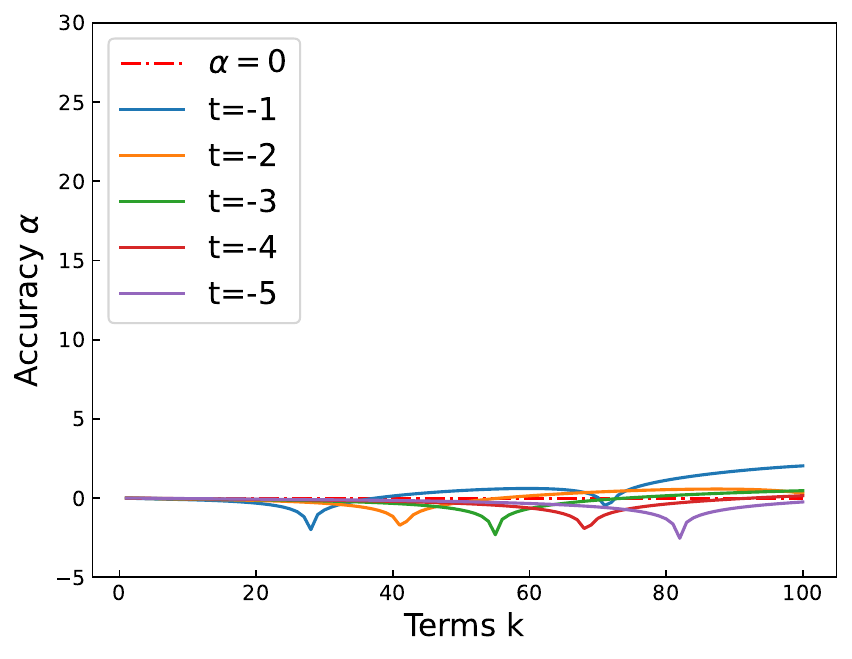}} 
    \subfigure[$d_7=2$, $\mu_9=1$]{\label{fig:ea2} \includegraphics[width=0.4\textwidth]{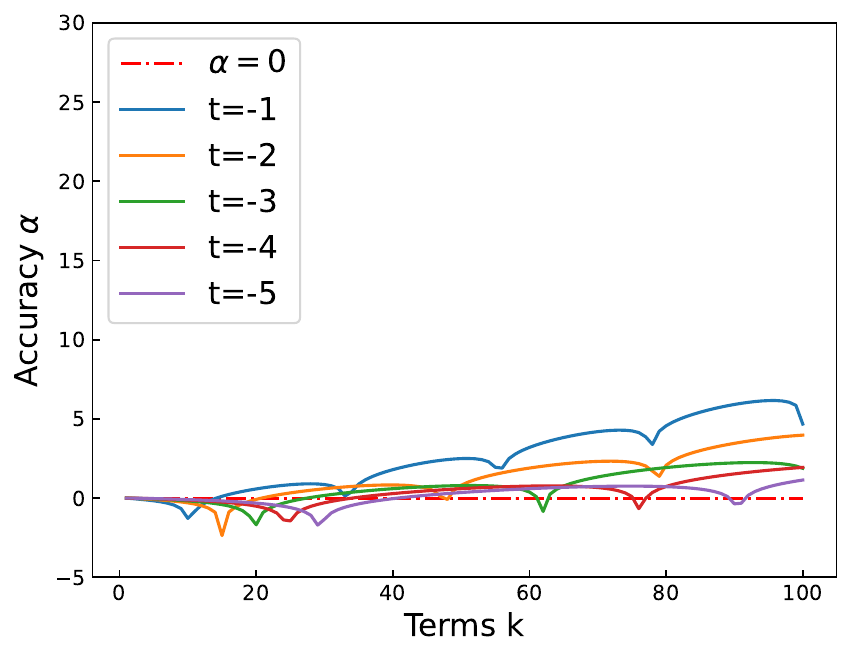}} 
    \subfigure[$d_5=3$, $\mu_8=2$]{\label{fig:ea3} \includegraphics[width=0.4\textwidth]{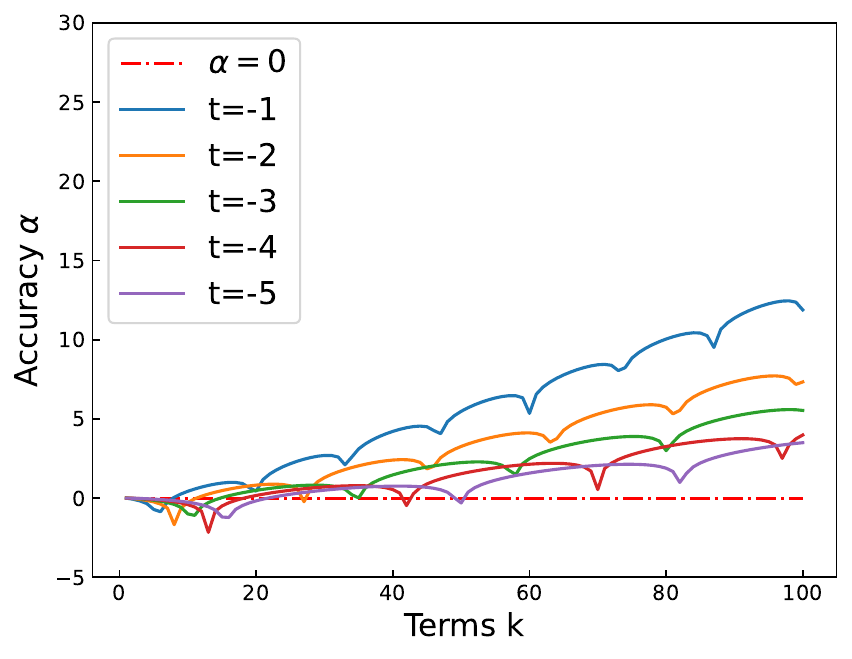}} 
    \subfigure[$d_3=4$, $\mu_7=3$]{\label{fig:ea4} \includegraphics[width=0.4\textwidth]{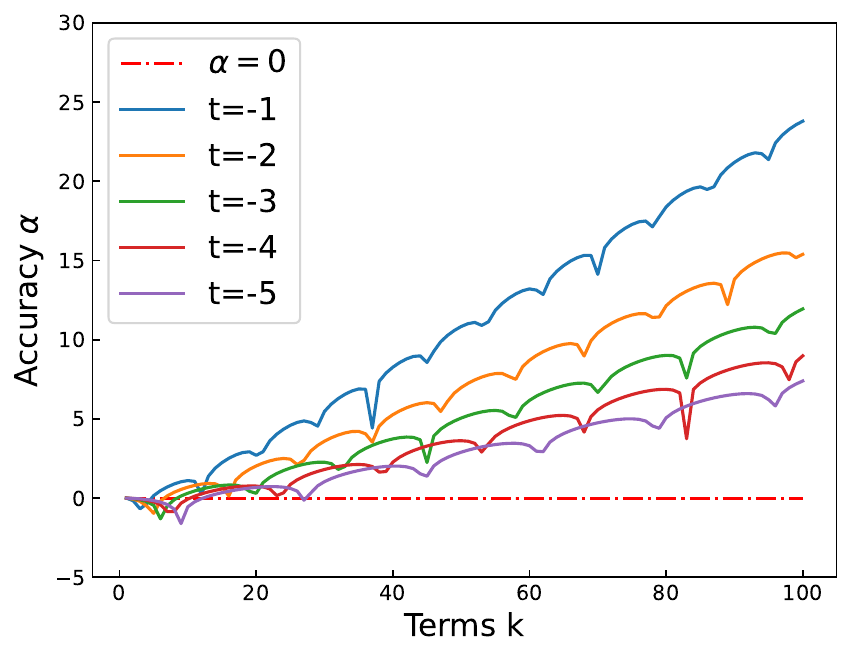}} 
    \subfigure[$d_4=6$, $\mu_4=7$]{\label{fig:ea6} \includegraphics[width=0.4\textwidth]{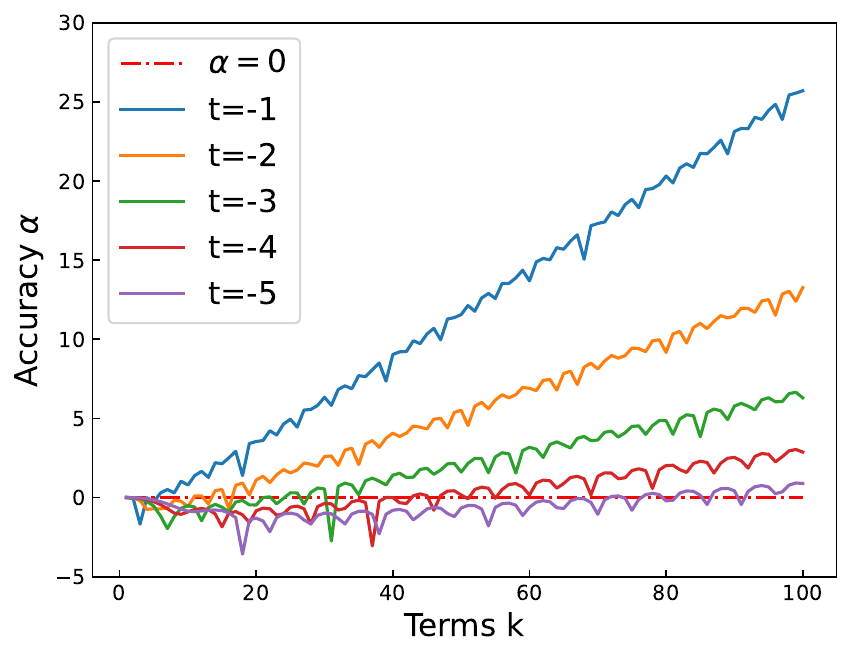}} 
    \subfigure[$d_8=7$, $\mu_3=8$]{\label{fig:ea7} \includegraphics[width=0.4\textwidth]{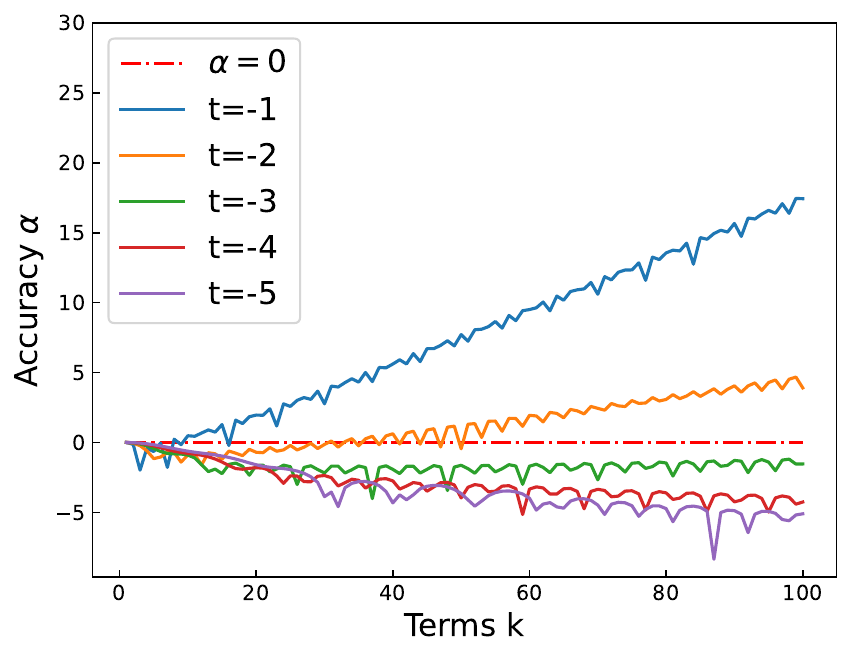}} 
    \subfigure[$d_9=8$, $\mu_2=9$]{\label{fig:ea8} \includegraphics[width=0.4\textwidth]{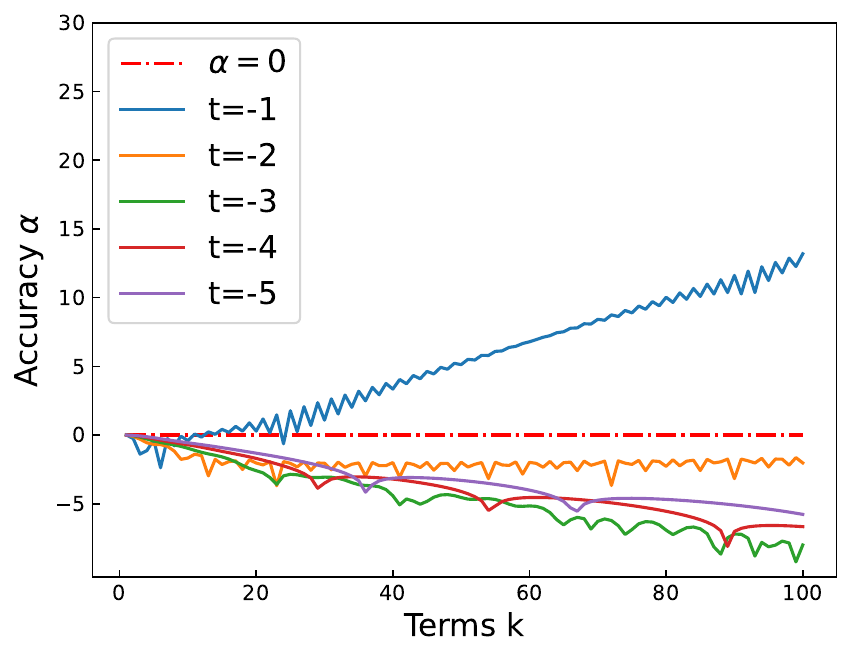}} 
    \subfigure[$d_{10}=9$, $\mu_1=10$]{\label{fig:ea9} \includegraphics[width=0.4\textwidth]{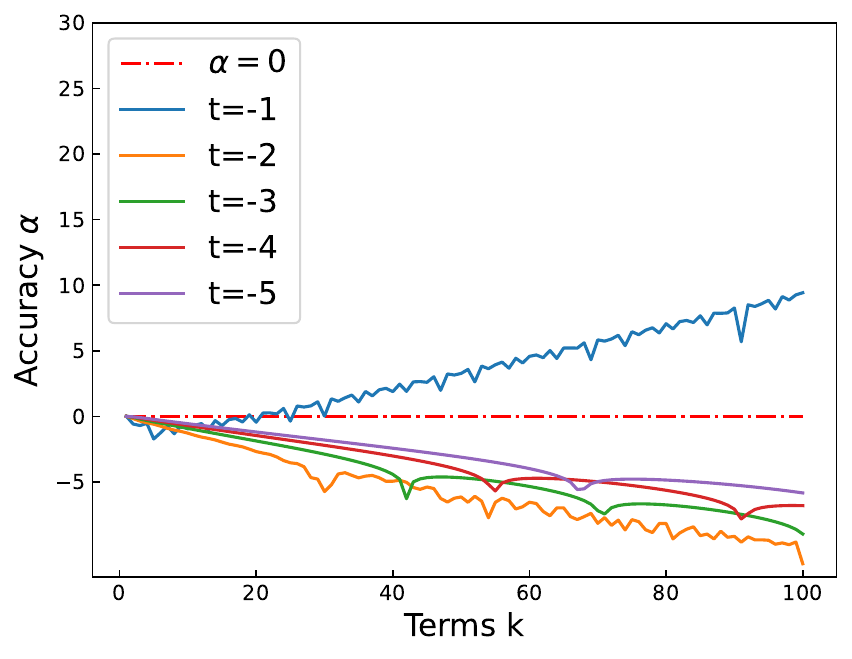}} 
    \caption{The accuracy $\alpha$ versus the number of terms $K$ and tuning parameter $t$ for each node with a unique degree in Example 3.}
    \label{fig:anti}
\end{figure}

The antiregular graph in Example 3 has 10 nodes, eight of which have a unique degree, except for two nodes $q=1$ and $q=2$, which have the same degrees $d_1=d_2=5$.
Figure \ref{fig:anti} demonstrates the convergence of the Euler series $\xi_{q;K}$ for each node with a unique degree in Example 3.
The Euler $t$-series $\xi_{q;K}$ diverges in Fig. \ref{fig:ea1}-\ref{fig:ea6} for degrees less than 7 but converges to a specific Laplacian eigenvalue in Fig. \ref{fig:ea7}-\ref{fig:ea9} with a tuning parameter $t<-1$ for node degrees no less than 7. 
The observation suggests that the increasing node degree transforms the Euler $t$-series $\xi_{q;K}$ from a diverging series into a converging series.
Another finding is that the performance of the accuracy $\alpha$ versus the number $K$ of terms fluctuates less and becomes almost linear as the node degree increases.
Overall, convergence is more likely to happen for nodes with large degrees.

\newpage
\subsection{Impact of the degree difference $\kappa$ on the convergence  of $\xi_{q;K} $ } \label{sec_impact_degree_difference}
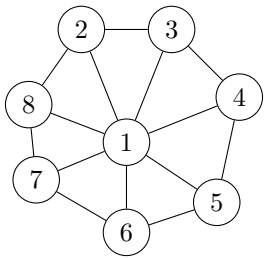
\begin{figure}[htbp]
    \centering
        \begin{tikzpicture}
            \node[draw,circle, draw, minimum size=0.6cm] (a) at (0,0) {1};
            \node[draw,circle, draw, minimum size=0.6cm] (b) at (-0.6,1.5) {2};
            \node[draw,circle, draw, minimum size=0.6cm] (c) at (0.6,1.5) {3};
            \node[draw,circle, draw, minimum size=0.6cm] (d) at (1.5,0.6) {4};
            \node[draw,circle, draw, minimum size=0.6cm] (e) at (1.2,-0.8) {5};
            \node[draw,circle, draw, minimum size=0.6cm] (f) at (0,-1.2) {6};
            \node[draw,circle, draw, minimum size=0.6cm] (g) at(-1.2,-0.5) {7};
            \node[draw,circle, draw, minimum size=0.6cm] (h) at (-1.3,0.5) {8};
            \foreach \x in {b,c,d,e,f,g,h} \draw (a) -- (\x);
            \draw (b) -- (c); \draw (c) -- (d);  \draw (d) -- (e);  \draw (e) -- (f);   \draw (g) -- (f);  \draw (h) -- (g);   \draw (b) -- (h); 
        \end{tikzpicture}
    \caption{An example of a graph $G_{(N,k)}$ with $N=8$ nodes and $k=1$.}
    \label{fig:wheel}
\end{figure}

In the last section, we have examined the influence of the node degree itself on the convergence of the Euler $t$-series $\xi_{q;K} $. 
Section \ref{sec_convergence_perturbation_Taylor_series} has shown that the degree difference between a unique degree $d_q$ of a node $q$ and degrees of other nodes, quantified by the degree difference $\kappa=\max_{\substack{1\leq
k\leq N\\k\neq q}}\frac{1}{\left\vert d_{q}-d_{k}\right\vert }\leq 1$ in (\ref{degree_difference_kappa}), influences the convergence of the perturbation series $\xi_{q;K} $.
To explore the impact of the degree difference $\kappa$ on the convergence of the Euler series $\xi_{q;K} $, we perform the computation on the graph $G_{(N,k)}$ that consists of a $2k$ nearest neighbor regular graph with $N-1$ nodes and a core node. 
A $2k$ nearest neighbor regular graph is a graph with each node joining its $k$ nearest neighbors. 
A core node \cite[p. 92]{PVM_graphspectra_second_edition} is a node that links to all other nodes in a graph except itself. 
Thus, only one node $q=N$ has a unique degree $d_{N}=N-1$ and other $N-1$ nodes have the same degree  $d_{i} = 2k+1$ in a graph $G_{(N,k)}$. 
Then,the degree difference $\kappa$ in (\ref{degree_difference_kappa}) in a graph $G_{(N,k)}$ is $\kappa = \frac{1}{N-2k-2}$ that increases as the parameter $k$ increases.
A graph $G_{(N,k)}$ with $N=8$ nodes and $k=1$ is shown in Fig. \ref{fig:wheel}. 

\begin{figure}[htpp]
    \centering
    \subfigure[$t=-1$]{\includegraphics[width=0.46\textwidth]{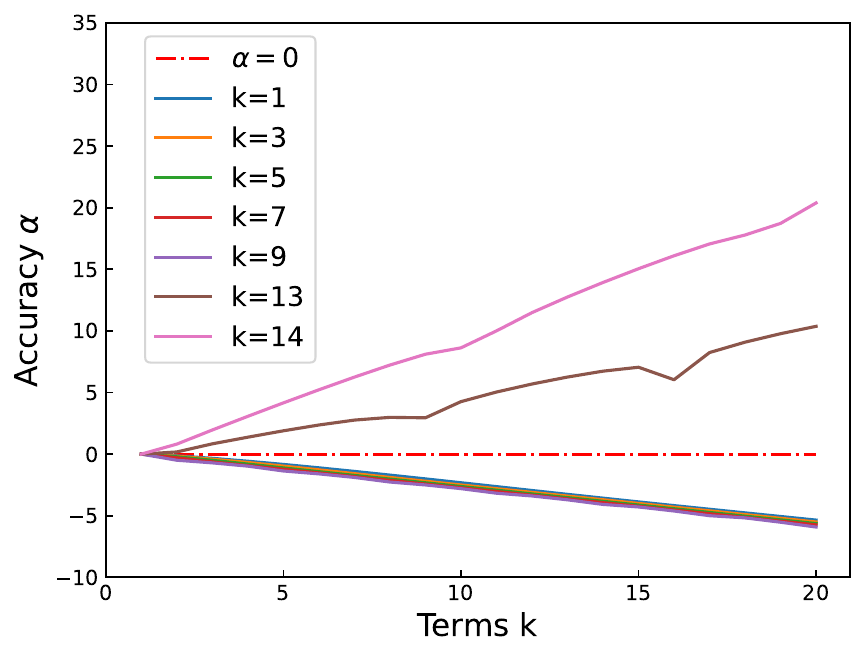} } 
    \subfigure[$t=-2$]{\includegraphics[width=0.46\textwidth]{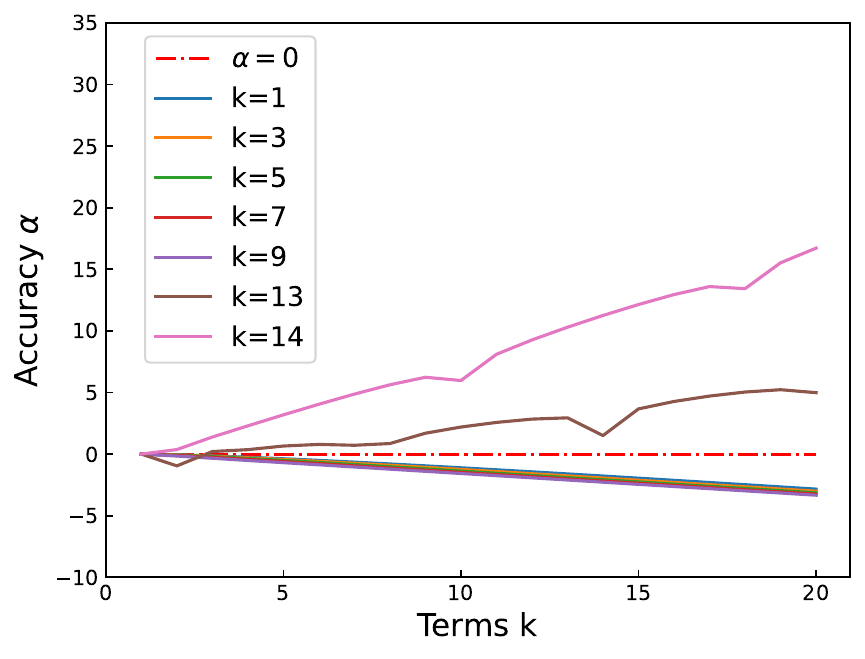} }
    \subfigure[$t=-3$]{\includegraphics[width=0.46\textwidth]{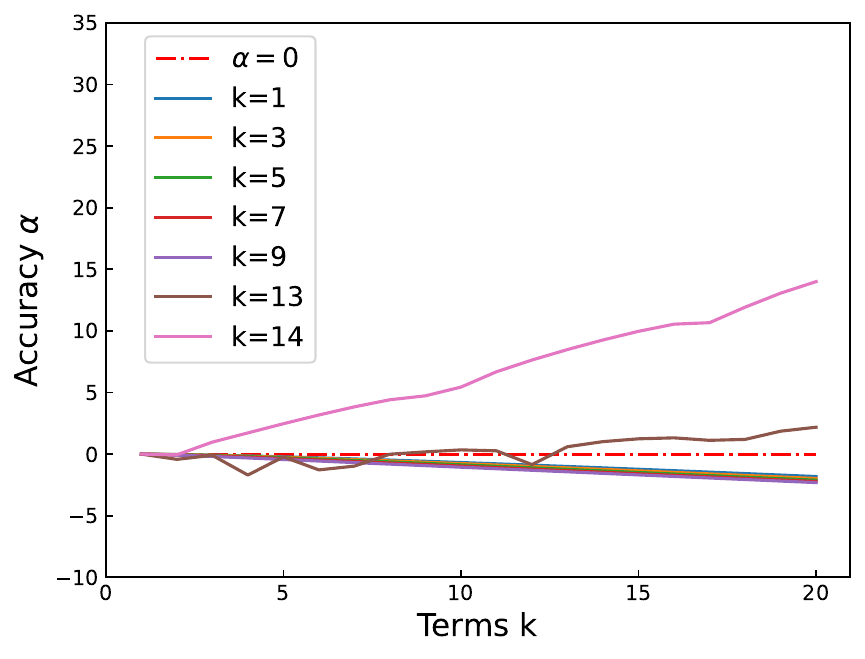} } 
    \subfigure[$t=-4$]{\includegraphics[width=0.46\textwidth]{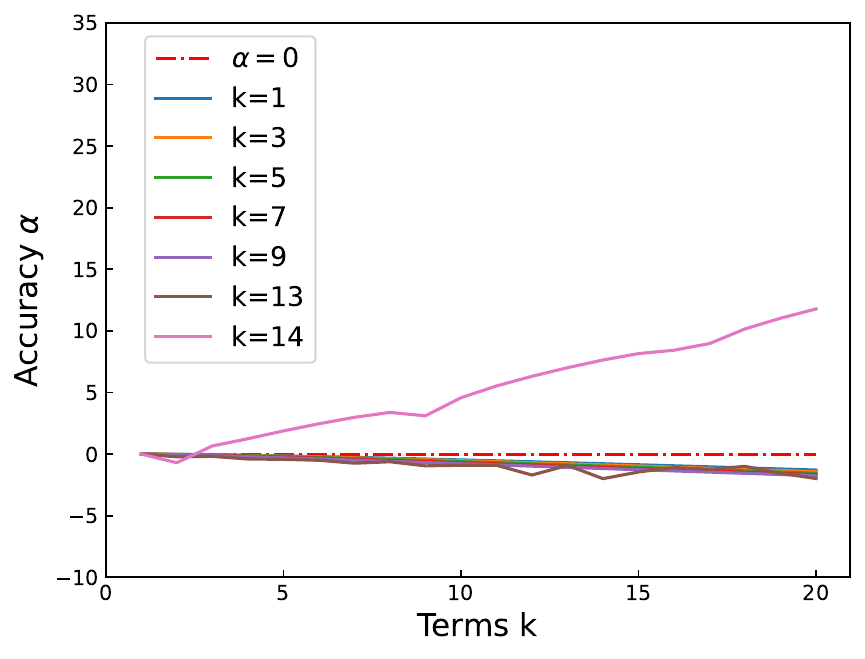} }
    \subfigure[$t=-5$]{\includegraphics[width=0.46\textwidth]{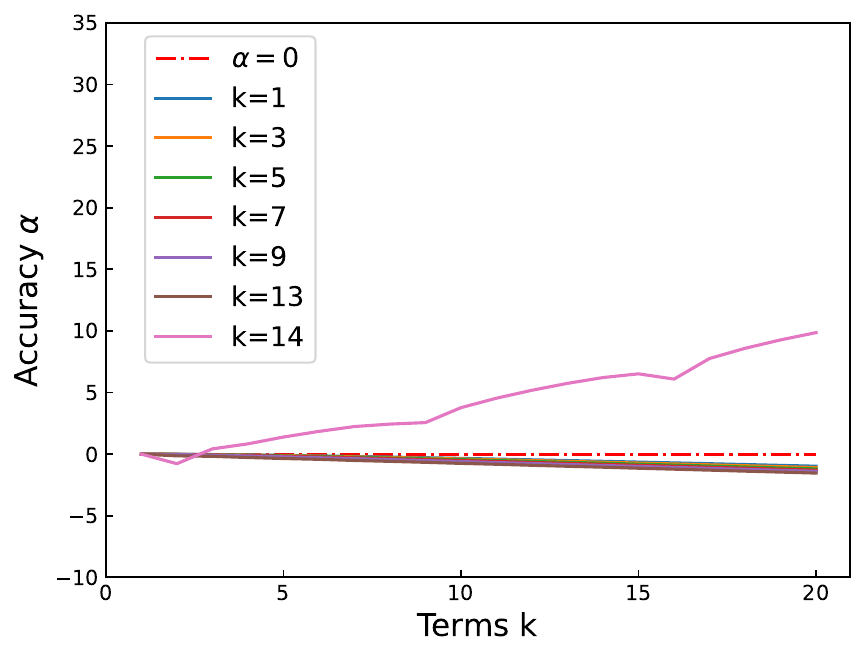} } 
    \subfigure[$t=-6$]{\includegraphics[width=0.46\textwidth]{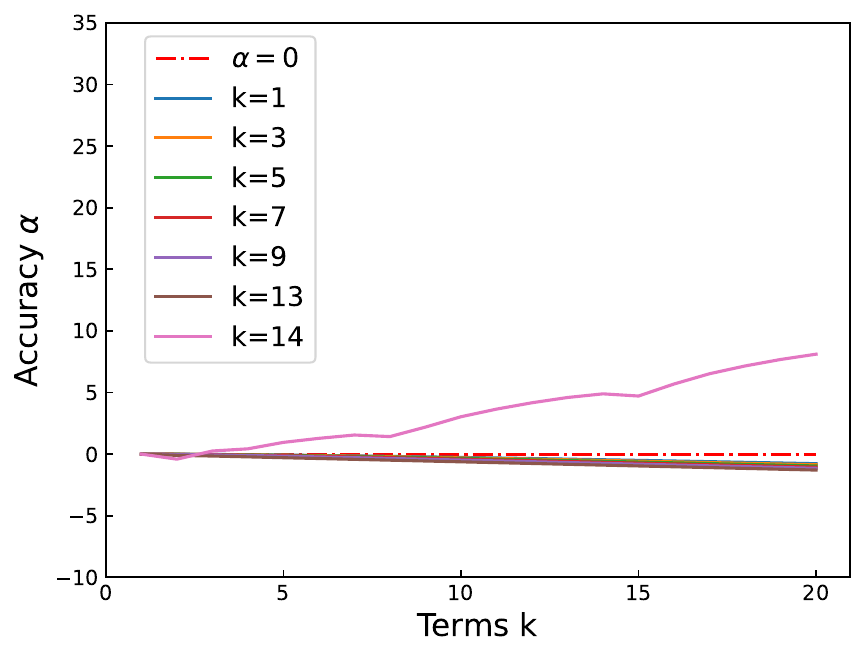} }
    \caption{The accuracy $\alpha$ of the Euler $t$-series $\xi_{31;K}$ of the node $q=31$ with the largest degree $d_{31}=30$ versus the number of terms $K$ and tuning parameter $t$ for each parameter $k$ in graphs $G_{(31, k)}$.}
    \label{fig:re}
\end{figure}

Figure \ref{fig:re} illustrates the convergence of the Euler $t$-series $\xi_{31; K}$ on graphs $G_{(N,k)}$ with $N=31$ and $1 \leq k \leq 14$. 
We notice that the Euler $t$-series $\xi_{31; K}$ transforms from convergence to divergence as the parameter $k$ increases while the tuning parameter $t$ is fixed.
Also, the decrease of the negative tuning parameter $t$ transforms the Euler series $\xi_{31; K}$ from diverging into converging when the value of $k$ is given, for example, $k=13$.
Another noteworthy finding is that if the Euler $t$-series  $\xi_{31; K}$ converges for some $k$, then the value of $k$ does not obviously impact the convergence speed.
The Euler series $\xi_{31; K}$ converges slower when the tuning parameter $t$ decreases from $t=-1$ to $-6$, which confirms our previous observation that the parameter $t$ influences the speed of convergence of the Euler series $\xi_{q; K}$.
The simulation conveys that the convergence is more likely to happen or be faster for a smaller $\kappa $.

\begin{figure}[htbp]
    \centering
    \subfigure[]{\label{fig:ak}    
    \includegraphics[width=0.45\textwidth]{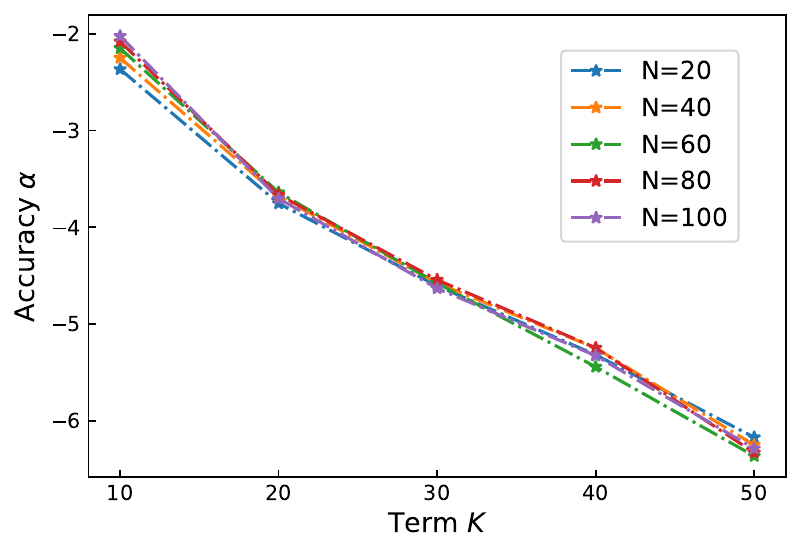}}\hspace{1cm}
    \subfigure[]{\label{fig:fe}
    \includegraphics[width=0.45\textwidth]{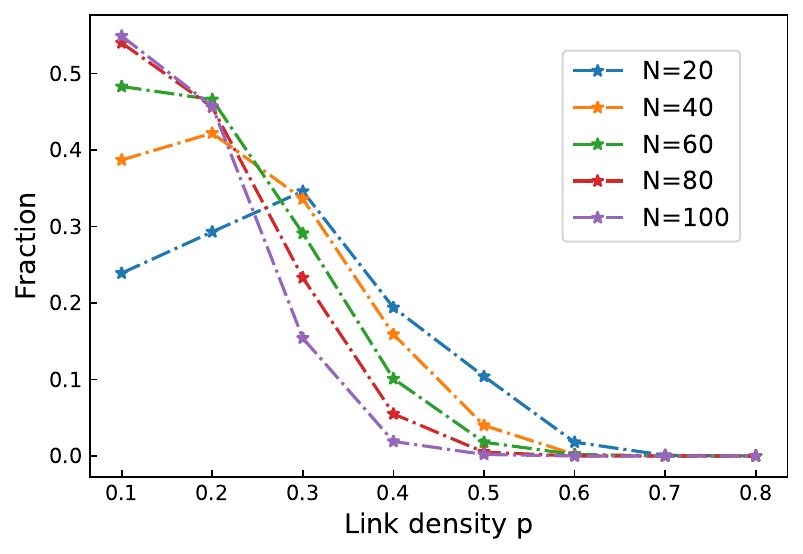}}
    \caption{(a) The accuracy $\alpha$ of the Euler $t$-series $\xi_{q;K}$ of the node with the largest and unique degree versus the number of terms $K$ on Erd\H{o}s-R\'{e}nyi graphs $G_{0.2}(N)$ with the parameter $t=-1$. (b) The fraction of convergent Euler $t$-series $\xi_{q; K}$ of the node with the largest and unique degree on Erd\H{o}s-R\'{e}nyi graphs $G_{p}(N)$ with the parameter $t=-1$.} 
    \label{fig:fra}
\end{figure}

Fig. \ref{fig:ak} depicts the accuracy $\alpha$ of convergent Euler $t$-series $\xi_{q; K}$ with $t=-1$ of the node with the largest and unique degree versus the number of terms $K$ in $10^3$ Erd\H{o}s-R\'{e}nyi graphs $G_{p}(N)$ of $N = \left\{ 20,40,60,80,100\right\}$ nodes and link density $p = 0.2$. 
The accuracy $\alpha$ decreases almost linearly versus the number of terms $K$.
The graph size $N$ hardly impacts the speed of convergence.
Fig. \ref{fig:fe} presents the fraction of convergent Euler $t$-series $\xi_{q; K}$ with $t=-1$ of the node with the largest and unique degree in $10^3$ Erd\H{o}s-R\'{e}nyi graphs $G_{p}(N)$ on $N = \left\{ 20,40,60,80,100 \right\}$ nodes and link density $p = \left\{ 0.1,0.2,0.3,0.4,0.5,0.6,0.7,0.8 \right\}$. 
We set the accuracy  $\alpha=-4$ and the number of terms $K=30$.
The Euler series  $\xi_{q; K}$ is more likely to converge to the largest Laplacian eigenvalue if the node is in graphs with larger sizes or lower link densities.
As the link density $p$ increases, convergence will likely happen in smaller graphs.
The fraction of convergent Euler series becomes very low when the link density $p \geq 0.7$ regardless of the graph size.
A reasonable explanation is that the increasing degree difference $\kappa$ in (\ref{degree_difference_kappa}) caused by an increasing link density $p$ lowers the possibility of convergence, consistent with what we observed in the last section.
Another finding is that there exists an optimal link density $p$ for the convergence of small graphs because Fig. \ref{fig:fe} shows that the fraction of convergent Euler series first increases and then decreases with the link density $p$ for graph sizes $N \leq 60$.

\section{Almost regular graph with one high degree node}
\label{sec:almost_regular}
In Section \ref {sec_coefficients_perturbation_series}, we emphasize that a unique degree is necessary for the expansion discussed in this work. 
We consider a special case of graphs that contain only one node with a unique degree--almost a regular graph.
An almost regular graph $G$ with $d_{j}=r$ for $2\leq j\leq N$ and
only node 1 has degree $d_{1}=d_{\max}>r$. 
Since the perturbation expansion
only holds for a unique degree $d_{q}$, it is clear that node $q$ here equals
node $1$.

\subsection{Deduction of the general form for the coefficient $c_{j}$}

The recursion in (\ref{beta_1r_r_not_q}) and (\ref{beta_jr_recursion})  becomes, for $l\neq1$,%
\[
\beta_{1l}=\frac{a_{l1}}{d_{1}-d_{l}}=\frac{a_{l1}}{d_{\max}-r}%
\]
and, for $j>1$,%
\begin{align*}
\beta_{jl}  &  =\frac{1}{r-d_{\max}}\sum_{m=1;m\neq1}^{N}\left\{  \sum
_{k=1}^{j-2}\beta_{kl}\beta_{j-k-1,m}a_{1m}-\beta_{j-1,m}a_{lm}\right\} \\
&  =\frac{1}{r-d_{\max}}\left(  \sum_{m=1}^{N}\left\{  \sum_{k=1}^{j-2}%
\beta_{kl}\beta_{j-k-1,m}a_{1m}-\beta_{j-1,m}a_{lm}\right\}  -\beta
_{j-1,1}a_{l1}\right) \\
&  =\frac{1}{d_{\max}-r}\sum_{m=2}^{N}\beta_{j-1,m}a_{lm}-\frac{1}{d_{\max}%
-r}\sum_{m=1}^{N}a_{1m}\sum_{k=1}^{j-2}\beta_{kl}\beta_{j-k-1,m}%
\end{align*}
However, it is not true that $\beta_{j1}=0$ (as follows from the computations
below). The coefficients (\ref{coeff_cj}) for $j>1$ are, for $q=1$,%
\[
c_{j}\left(  q\right)  =\sum_{k=1;k\neq q}^{N}\beta_{j-1,k}a_{qk}=\sum
_{l=1}^{N}\beta_{j-1,l}a_{1l}%
\]
which can be computed up to any desired value of $j$.

We compute now for several values of $j$ the coefficients $\beta_{jl}$ with
$x=d_{\max}-r$, by iterating the recursion%
\begin{equation}
\beta_{jl}=\frac{1}{x}\left(  \sum_{m=1}^{N}\beta_{j-1,m}a_{lm}-\beta
_{j-1,1}a_{l1}-\sum_{m=1}^{N}a_{1m}\sum_{k=1}^{j-2}\beta_{kl}\beta
_{j-k-1,m}\right)  \label{beta_jl_recursion_almost_regular_graph}%
\end{equation}
together with the corresponding coefficient $c_{j+1}\left(  q\right)  $ with
$q=1$. We rewrite the recursion in terms of $c_{j}\left(  q\right)  $ as%
\begin{align*}
\beta_{jl}  &  =\frac{1}{x}\left(  \sum_{m=2}^{N}\beta_{j-1,m}a_{lm}%
-\sum_{k=1}^{j-2}\beta_{kl}\sum_{m=1}^{N}a_{1m}\beta_{j-k-1,m}\right) \\
&  =\frac{1}{x}\left(  \sum_{m=2}^{N}\beta_{j-1,m}a_{lm}-\sum_{k=1}^{j-2}%
\beta_{kl}c_{j-k}\left(  q\right)  \right)
\end{align*}
Also,%
\begin{align*}
c_{j+1}\left(  q\right)   &  =\sum_{l=1}^{N}\beta_{jl}a_{1l}=\frac{1}{x}%
\sum_{l=1}^{N}\left(  \sum_{m=2}^{N}\beta_{j-1,m}a_{lm}a_{1l}-\sum_{k=1}%
^{j-2}\beta_{kl}c_{j-k}\left(  q\right)  a_{1l}\right) \\
&  =\frac{1}{x}\left(  \sum_{m=2}^{N}\beta_{j-1,m}\sum_{l=1}^{N}a_{lm}%
a_{1l}-\sum_{k=1}^{j-2}c_{j-k}\left(  q\right)  \sum_{l=1}^{N}\beta_{kl}%
a_{1l}\right) \\
&  =\frac{1}{x}\left(  \sum_{m=2}^{N}\beta_{j-1,m}\left(  A^{2}\right)
_{m1}-\sum_{k=1}^{j-2}c_{j-k}\left(  q\right)  c_{k+1}\left(  q\right)
\right)
\end{align*}
We arrive at the pseudo-recursion for the coefficient $c_{j}\left(  q\right)
$
\begin{equation}
c_{j}\left(  q\right)  =\frac{1}{x}\left(  \sum_{m=1}^{N}\beta_{j-2,m}\left(
A^{2}\right)  _{m1}-\beta_{j-2,1}\left(  A^{2}\right)  _{11}-\sum_{k=1}%
^{j-3}c_{j-1-k}\left(  q\right)  c_{k+1}\left(  q\right)  \right)
\label{pseudo_recursion_c_j}%
\end{equation}

Both the coefficients $\beta_{jl}$ and $c_{k}\left(  q\right)  $ are computed and summarized 
in Appendix \ref{sec_computation_coeff_almost_regular_graphs} up to $j=9$.
Thus, we find the sequence in terms of the number
$\left(  A^{m}\right)  _{11}$ of closed walks with length $m$ at node $q=1$,%

We guess the general solution and propose, inspired by our characteristic
coefficients (chc) \cite{PVM_charcoef}, first defined in \cite{PVM_ASYM}, the form%
\[
c_{m}=\frac{1}{x^{m-1}}\sum_{k=1}^{\frac{m}{2}}g_{k}\left(  m\right)
\mathcal{A}\left[  k,m\right]
\]
where%
\begin{equation}
\mathcal{A}[k,m]=\sum_{\sum_{i=1}^{k}j_{i}=m;j_{i}>0}\prod_{i=1}^{k}\left(
A^{j_{i}}\right)  _{11} \label{chc_adjacency_closed_walks_node_q}.%
\end{equation}
Thus,
\begin{align*}
\mathcal{A}[1,m]  &  =\sum_{\sum_{i=1}^{1}j_{1}=m;j_{i}>0}\prod_{i=1}%
^{1}\left(  A^{j_{i}}\right)  _{11}=\left(  A^{m}\right)  _{11}\\
\mathcal{A}[2,m]  &  =\sum_{\sum_{i=1}^{2}j_{i}=m;j_{i}>0}\prod_{i=1}%
^{2}\left(  A^{j_{i}}\right)  _{11}=\sum_{j_{1}+j_{2}=m;j_{i}>0}\left(
A^{j_{1}}\right)  _{11}\left(  A^{j_{2}}\right)  _{11}=\sum_{j=1}^{m-1}\left(
A^{j}\right)  _{11}\left(  A^{m-j}\right)  _{11}\\
\mathcal{A}[3,m]  &  =\sum_{\sum_{i=1}^{3}j_{i}=m;j_{i}>0}\prod_{i=1}%
^{3}\left(  A^{j_{i}}\right)  _{11}=\sum_{j_{1}+j_{2}+j_{3}=m;j_{i}>0}\left(
A^{j_{i}}\right)  _{11}\left(  A^{j_{2}}\right)  _{11}\left(  A^{j_{3}%
}\right)  _{11}\\
&  =\sum_{j=1}^{m-2}\sum_{i=1}^{m-j-1}\left(  A^{i}\right)  _{11}\left(
A^{j}\right)  _{11}\left(  A^{m-i-j}\right)  _{11}%
\end{align*}
and all other combinations follow the chc recursion.
We have placed the theory in Appendix. \ref{sec_ recursion_characteristic_coefficient_A}.

\subsection{Closed form for the coefficient $c_{m}$}

Comparison with the computations of $c_{k \leq 9}$ in Appendix. \ref{sec_computation_coeff_almost_regular_graphs} indicates
that $g_{k}\left(  m\right)  =\frac{\left(  -1\right)  ^{k-1}}{k}%
\binom{m+(k-2)}{k-1}$ and $g_{1}\left(  m\right)  =1$. Hence, we arrive at%
\begin{equation}
c_{m}=\frac{1}{x^{m-1}}\left(  \sum_{k=1}^{\frac{m}{2}}\frac{\left(
-1\right)  ^{k-1}}{k}\binom{m+(k-2)}{k-1}\mathcal{A}\left[  k,m\right]
\right)  \label{estimated_analytic_form}%
\end{equation}
The equation( \ref{estimated_analytic_form}) has been numerically verified for $m < 10$.

Since $\left(  A\right)  _{11}=a_{11}=0$, we can also extend\footnote{Indeed,
assume for the summation conditions $\sum_{i=1}^{k}j_{i}=m$ and $j_{i}>0$ that
all $j_{i}\geq2$, then $\sum_{i=1}^{k}j_{i}=m$ implies that $\sum_{i=1}%
^{k}j_{i}\geq2k$ or $m\geq2k$ and, thus, $k\leq\frac{m}{2}$. Hence, if
$k>\frac{m}{2}$, then there must be at least one $j_{i}$-index that equals
$j_{i}=1$. But $\left(  A^{j_{i}}\right)  _{11}=\left(  A^{1}\right)  _{11}%
=0$, which means that $\mathcal{A}[k,m]=0$ for $k>\frac{m}{2}$.} the
upperbound from $\frac{m}{2}$ to $m$, thus%
\begin{equation}
c_{m}=\frac{1}{x^{m-1}}\left(  \sum_{k=1}^{m}\frac{\left(  -1\right)  ^{k-1}%
}{k}\binom{m+k-2}{k-1}\mathcal{A}\left[  k,m\right]  \right)
\label{c_k_almost_regular_analytic}%
\end{equation}
Indeed, for $m=1$, it holds that $c_{1}=0$. Series computations with upper
bound $m$ are considerably easier than with $\frac{m}{2}$. The coefficient
$c_{m}$ in (\ref{estimated_analytic_form}) has been confirmed by numerical
computations. Since $\binom{m+k-2}{k-1}\mathcal{A}\left[  k,m\right]  $ is
non-negative, the $k$-sum in (\ref{c_k_almost_regular_analytic}) is
alternating in $k$.

We introduce the contour integral (\ref{s_general}) into
(\ref{c_k_almost_regular_analytic})%
\begin{align*}
c_{m}  &  =\frac{1}{x^{m-1}}\frac{1}{2\pi i}\int_{C(z_{0})}\frac{dz}%
{(z-z_{0})^{m+1}}\sum_{k=1}^{m}\frac{\left(  -1\right)  ^{k-1}}{k}%
\binom{m+k-2}{k-1}[f(z)-f(z_{0})]^{k}\\
&  =\,\frac{1}{x^{m-1}}\frac{1}{2\pi i}\int_{C(z_{0})}\frac{dz}{(z-z_{0}%
)^{m+1}}\sum_{k=1}^{\infty}\frac{\left(  -1\right)  ^{k-1}}{k}\binom
{m+k-2}{k-1}[f(z)-f(z_{0})]^{k}%
\end{align*}
because $\int_{C(z_{0})}\frac{[f(z)-f(z_{0})]^{k}}{(z-z_{0})^{m+1}}dz=0$ for
$k>m$, since the integrand is analytic around $z_{0}$. Since
\[
\int_{0}^{z}\left(  1+u\right)  ^{-s}du=\sum_{k=0}^{\infty}\binom{-s}{k}%
\frac{z^{k+1}}{k+1}=\sum_{k=1}^{\infty}\binom{-s}{k-1}\frac{z^{k}}{k}%
\]
indicates that%
\[
\frac{\left(  1+z\right)  ^{1-s}-1}{1-s}=\sum_{k=1}^{\infty}\binom{-s}%
{k-1}\frac{z^{k}}{k}%
\]
With the binomial property, valid for any complex $s$,%
\begin{equation}
{\binom{-s}{j}}=(-1)^{j}\;\frac{\Gamma(s+j)}{j!\,\Gamma(s)}=(-1)^{j}%
\;\binom{s-1+j}{j} \label{binomial_min_z_in_positive_z}%
\end{equation}
We find that%
\[
\frac{\left(  1+z\right)  ^{1-s}-1}{1-s}=\sum_{k=1}^{\infty}\frac{\left(
-1\right)  ^{k-1}}{k}\binom{s+k-2}{k-1}z^{k}%
\]
Hence,%
\[
c_{m}=\frac{1}{x^{m-1}}\frac{1}{2\pi i}\int_{C(z_{0})}\frac{dz}{(z-z_{0}%
)^{m+1}}\frac{\left(  1+f(z)-f(z_{0})\right)  ^{1-m}-1}{1-m}%
\]
and with $f\left(  z_{0}\right)  =1$ for $z_{0}=0$, in $f\left(  z\right)
=\left(  A^{z}\right)  _{11}$, we arrive, for $m>1$, at%
\begin{align}
c_{m}  &  =\frac{1}{\left(  m-1\right)  x^{m-1}}\frac{1}{2\pi i}\int
_{C(0)}\frac{\left(  1-\left(  f(z)\right)  ^{1-m}\right)  dz}{z^{m+1}%
}\nonumber\\
&  =\frac{-1}{\left(  m-1\right)  x^{m-1}}\frac{1}{2\pi i}\int_{C(0)}\frac
{dz}{\left(  f(z)\right)  ^{m-1}z^{m+1}} \label{c_contour_integral}%
\end{align}
which also equals%
\[
c_{m}=\frac{-1}{\left(  m-1\right)  x^{m-1}}\frac{1}{m!}\left.  \frac{d^{m}%
}{dz^{m}}\left(  \frac{1}{\left(  f(z)\right)  ^{m-1}}\right)  \right\vert
_{z=0}%
\]
and which is related to the Lagrange series coefficient in (\ref{Lagrange}).

\subsection{Perturbed eigenvalue $\xi_{q}\left(  \zeta\right)  $ of the
matrix $\Delta+\zeta A$}

The perturbed eigenvalue of the matrix $\Delta+\zeta A$ in (\ref{Taylor_perturbation_series})
becomes, for an almost regular graph with degree $\left(  A^{2}\right)
_{jj}=r$ for $2\leq j\leq N$ and one node $q$ with degree $d_{q}=d_{\max
}=\left(  A^{2}\right)  _{11}$,%
\begin{equation}
\xi_{q}\left(  \zeta\right)  =\left(  A^{2}\right)  _{11}+x\sum_{m=2}^{\infty
}\left(  \sum_{k=1}^{m}\frac{\left(  -1\right)  ^{k-1}}{k}\binom{m+k-2}%
{k-1}\mathcal{A}\left[  k,m\right]  \right)  \left(  \frac{\zeta}{x}\right)
^{m} \label{perturbed_eigenvalue_zeta_almost_reg_graph}%
\end{equation}
where $x=d_{\max}-r\geq1$. If $\zeta=-1$ and the series in
(\ref{perturbed_eigenvalue_zeta_almost_reg_graph}) converges, then $\xi
_{q}\left(  -1\right)  $ represents a Laplacian eigenvalue $\mu_{q}$ of an
almost regular graph that is closest to the degree $d_{q}=d_{\max}=\left(
A^{2}\right)  _{11}$. On the other hand, if $\zeta=1$, then $\xi_{q}\left(
1\right)  $ equals the eigenvalue of the signless Laplacian $\overline{Q}%
=\Delta+A$ of an almost regular graph. However, as argued above, it is
unlikely that the series in (\ref{perturbed_eigenvalue_zeta_almost_reg_graph})
converges for perturbation parameter $\left\vert \zeta\right\vert =1$.
However, numerical computations seem to indicate that the Euler transformation
of (\ref{perturbed_eigenvalue_zeta_almost_reg_graph}) converges.

Introducing the contour integral (\ref{c_contour_integral}) in
(\ref{Taylor_perturbation_series}) yields%
\begin{align*}
\xi_{q}\left(  \zeta\right)   &  =d_{q}-\frac{1}{2\pi i}\int_{C(0)}\sum
_{j=2}^{\infty}\left(  \frac{1}{\left(  j-1\right)  x^{j-1}}\frac{\zeta^{j}%
}{\left(  f(z)\right)  ^{j-1}z^{j+1}}\right)  dz\\
&  =d_{q}-\frac{\zeta}{2\pi i}\int_{C(0)}\sum_{j=1}^{\infty}\left(  \frac
{1}{j}\left(  \frac{\zeta}{xzf\left(  z\right)  }\right)  ^{j}\right)
\frac{dz}{z^{2}}%
\end{align*}
Provided $\left\vert \frac{\zeta}{xzf\left(  z\right)  }\right\vert <1$, then
we can use the Taylor series $\log\left(  1-y\right)  =-\sum_{j=1}^{\infty
}\frac{y^{j}}{j}$ and we obtain%
\[
\xi_{q}\left(  \zeta\right)  =d_{q}+\frac{\zeta}{2\pi i}\int_{C(0)}\log\left(
1-\frac{\zeta}{xzf\left(  z\right)  }\right)  \frac{dz}{z^{2}}%
\]
Finally, since $f\left(  z\right)  =\left(  A^{z}\right)  _{11}$, we arrive at
the complex integral for the eigenvalue $\xi_{q}\left(  \zeta\right)  $ of the
matrix $\Delta+\zeta A$ close to the degree $d_{q}$ in an almost regular graph
with adjacency matrix $A$,%
\begin{equation}
\xi_{q}\left(  \zeta\right)  =d_{q}+\frac{\zeta}{2\pi i}\int_{C(0)}\log\left(
1-\frac{\zeta}{xz\left(  A^{z}\right)  _{11}}\right)  \frac{dz}{z^{2}%
}\label{eigenvalue_as_contour_integral}%
\end{equation}
Evaluting the contour in (\ref{eigenvalue_as_contour_integral}) around a
circle at the origin with radius $r$, small enough to avoid enclosing the
smallest (in absolute value) zero of $f\left(  z\right)  $%
\begin{equation}
\xi_{q}\left(  \zeta\right)  =d_{q}+\frac{\zeta}{2\pi r}\int_{0}^{2\pi
}e^{-i\theta}\log\left(  1-\frac{\zeta}{xr}\frac{e^{-i\theta}}{\left(
A^{re^{i\theta}}\right)  _{11}}\right)  d\theta
\label{eigenvalue_as_contour_integral_circle}%
\end{equation}
If the function $f\left(  z\right)  =\left(  A^{z}\right)  _{11}$ is known
analytically, then the integral (\ref{eigenvalue_as_contour_integral_circle})
can be computed numerically.

Choosing $q$ equal to the highest degree node $d_{q}=d_{\max}=\left(
A^{2}\right)  _{11}$, then the Euler $t$-transform in (\ref{mu_perturbation_Eulersummation_general_t}) of the perturbed eigenvalue
for an almost regular graph is
\begin{align}
\xi_{q}\left(  \zeta\right)   &  =d_{q}+\sum_{m=2}^{\infty}\left(  \sum
_{k=2}^{m}\binom{m-1}{k-1}\frac{1}{x^{k-1}}\left(  \sum_{j=1}^{k}\frac{\left(
-1\right)  ^{j-1}}{j}\binom{k+j-2}{j-1}\mathcal{A}\left[  j,k\right]  \right)
t^{m-k}\right)  \left(  \frac{\zeta}{1+t\zeta}\right)  ^{m}\nonumber\\
&  =d_{q}+x\sum_{m=1}^{\infty}\left(  \sum_{k=1}^{m}\sum_{j=1}^{k}
\frac{\left(  -1\right)  ^{j-1}}{j}\binom{m-1}{k-1}\binom{k+j-2}
{j-1}\mathcal{A}\left[  j,k\right]  \left(  tx\right)  ^{-k}\right)  \left(
\frac{t\zeta}{1+t\zeta}\right)  ^{m}
\label{perturbed_eigenvalue_zeta_almost_reg_graph_Euler}
\end{align}

\subsection{Numerical assessment of the perturbation series (\ref{perturbed_eigenvalue_zeta_almost_reg_graph}) and its corresponding Euler series
(\ref{perturbed_eigenvalue_zeta_almost_reg_graph_Euler})}
\label{sec:numerical_almost_regular}

\begin{figure}[htbp]
    \centering
    \subfigure[$k=1$ ]{ \label{fig:c1} \includegraphics[width=0.47\textwidth]{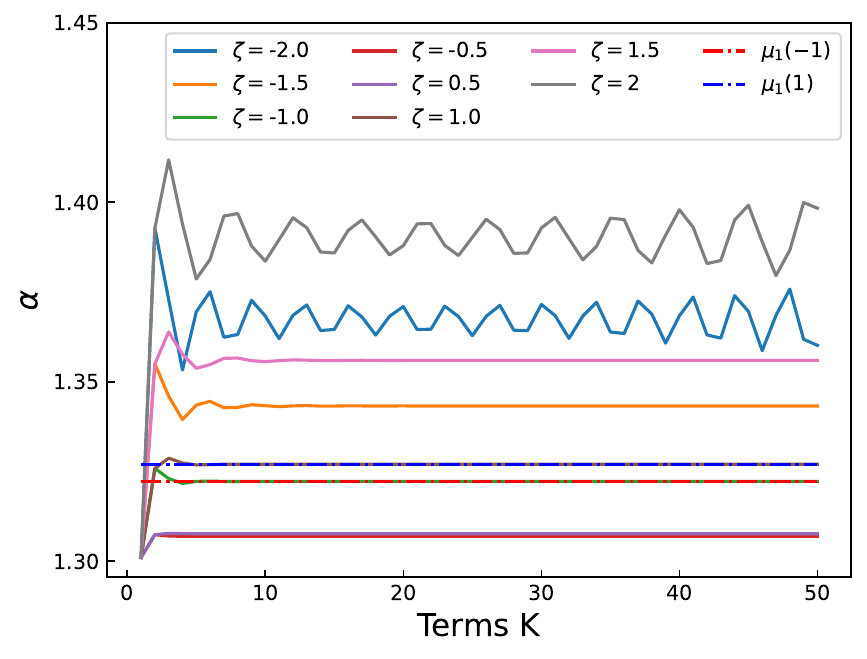} }
    \subfigure[$k=9$]{ \label{fig:c2} \includegraphics[width=0.47\textwidth]{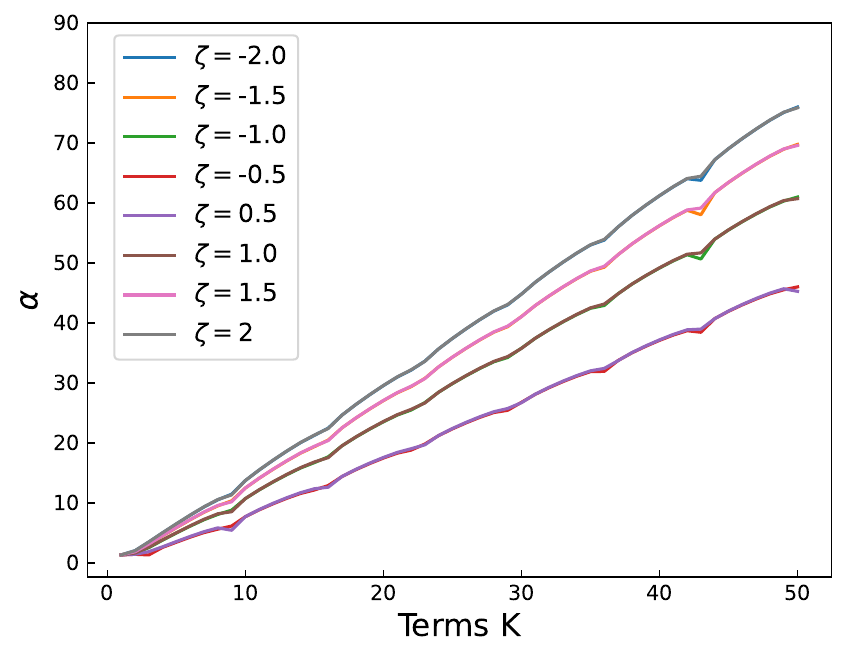} }
    \caption{The value of $\alpha $  of the series
    (\ref{perturbed_eigenvalue_zeta_almost_reg_graph}) versus parameter $\zeta$ and the number of terms $K$ of graph $G_{(21,k)}$.}
    \label{fig:cc}
\end{figure}

We compare the convergence between the series
(\ref{perturbed_eigenvalue_zeta_almost_reg_graph}) and its Euler transform
(\ref{perturbed_eigenvalue_zeta_almost_reg_graph_Euler}) on the graph $G_{(N,k)}$ defined in Section \ref{sec_impact_degree_difference}.
Thus, only one node has a unique degree $d_{max}=N-1$ and the other $N-1$ nodes have the same degree  $r = 2k+1$ in a graph $G_{(N,k)}$. 
The degree difference is $x=d_{\max}-r=N-2k$, which decreases as the parameter $k$ increases.
We define $10^{\alpha} := \left| \xi_{q}\left(  \zeta\right) \right|$ and present the value of $\alpha$ of the perturbation series
(\ref{perturbed_eigenvalue_zeta_almost_reg_graph}) in graphs $G_{(21,k)}$ with $k=1$ and $ 9$ under different $\zeta$ in Fig. \ref{fig:cc}.
The dash-dot line represents the largest eigenvalue $\mu_{1}(\zeta)$ of the matrix $Q=\Delta+\zeta A$.
Figure. \ref{fig:c1} shows that the perturbation series
(\ref{perturbed_eigenvalue_zeta_almost_reg_graph}) coverges when $k=1$ and $\left| \zeta \right| \leq 1$.
Moreover, the perturbation series
(\ref{perturbed_eigenvalue_zeta_almost_reg_graph}) converges to the lagrest eigenvalue $\mu_1(-1)$ of Laplacian $Q=\Delta+\zeta A$ when $\zeta=-1$.
Similarly, the perturbation series
(\ref{perturbed_eigenvalue_zeta_almost_reg_graph}) converges to the lagrest eigenvalue $\mu_1(1)$ of signless Laplacian $Q=\Delta+\zeta A$ when $\zeta=1$.
However, if $k$ is large and the degree difference $x$ is relatively small, which means all degrees are close to each other, then the perturbation series
(\ref{perturbed_eigenvalue_zeta_almost_reg_graph}) hardly converges for $ \left| \zeta \right|$ around 1 as shown in Figure \ref{fig:c2}.
\begin{figure}[htbp]
    \centering
    \subfigure[$k=1$]{ \label{fig:ec1} \includegraphics[width=0.47\textwidth]{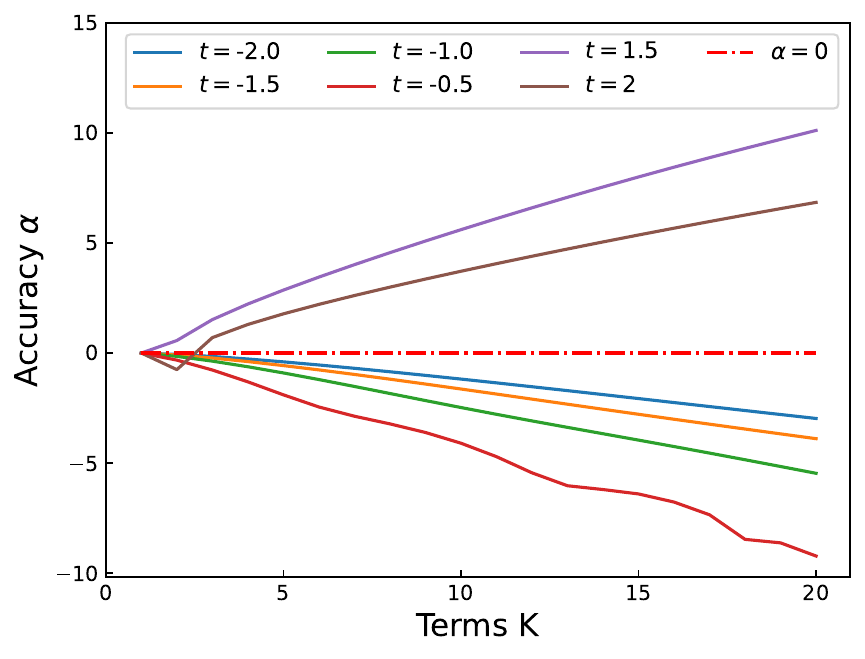} }
    \subfigure[$k=9$]{ \label{fig:ec2} \includegraphics[width=0.47\textwidth]{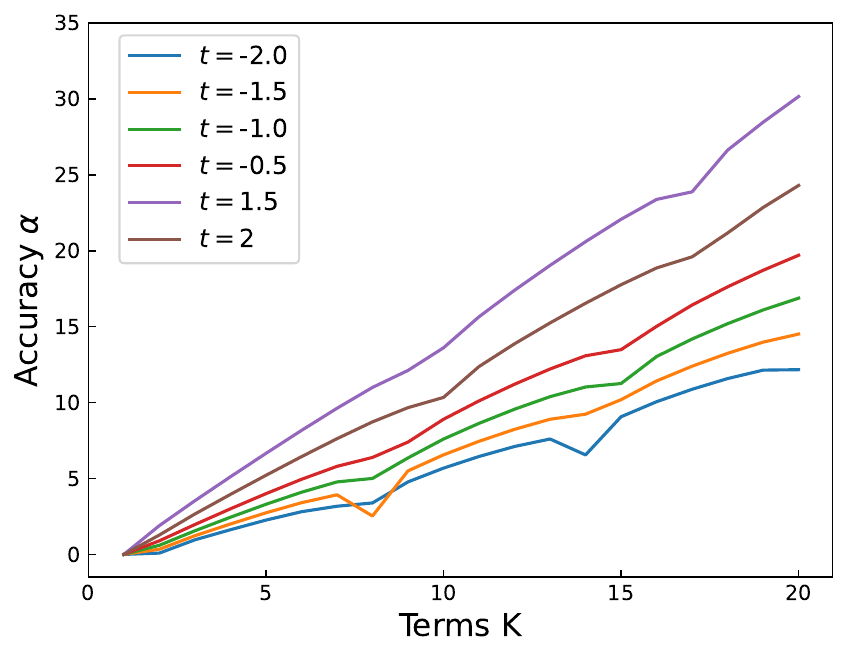} }
    \caption{The accuracy $\alpha$ of the Euler series   
    (\ref{perturbed_eigenvalue_zeta_almost_reg_graph_Euler}) versus parameter $t$ and the number of terms $K$ of graph $G_{(21,k)}$ with $\zeta=-1$.}
    \label{fig:ec}
\end{figure}

Next, we examine the convergence of the Euler series in
(\ref{perturbed_eigenvalue_zeta_almost_reg_graph_Euler}).
Similarly, we define the difference between the Euler series and the largest  eigenvalue $\mu_{1}(\zeta)$ as
\begin{equation*}
\label{eq:accuracy_delta}
    10^{\alpha} := \left| \xi_{q}\left(  \zeta\right)  - \mu_{1}(\zeta)\right|
\end{equation*}
where the parameter $\alpha$ measures the accuracy of the series $\xi_{q}\left(  \zeta\right) $. 
If the accuracy $\alpha$ becomes increasingly negative, then the series $\xi_{q}\left(  \zeta\right) $ converges to the eigenvalue $\mu_{1}(\zeta)$ .
Figure. \ref{fig:ec} illustrates that the Euler series (\ref{perturbed_eigenvalue_zeta_almost_reg_graph_Euler}) can also converge with a suitable parameter $t$ when the perturbation series
(\ref{perturbed_eigenvalue_zeta_almost_reg_graph}) converges with $k=1$ and $\zeta=-1$.
The Euler series also hardly converges when the perturbation series diverges with $k=9$ and $\zeta=-1$ no matter what the parameter $t$ is.

\begin{figure}[htbp]
    \centering
    \subfigure[$k=1$]{ \label{fig:ec11} \includegraphics[width=0.47\textwidth]{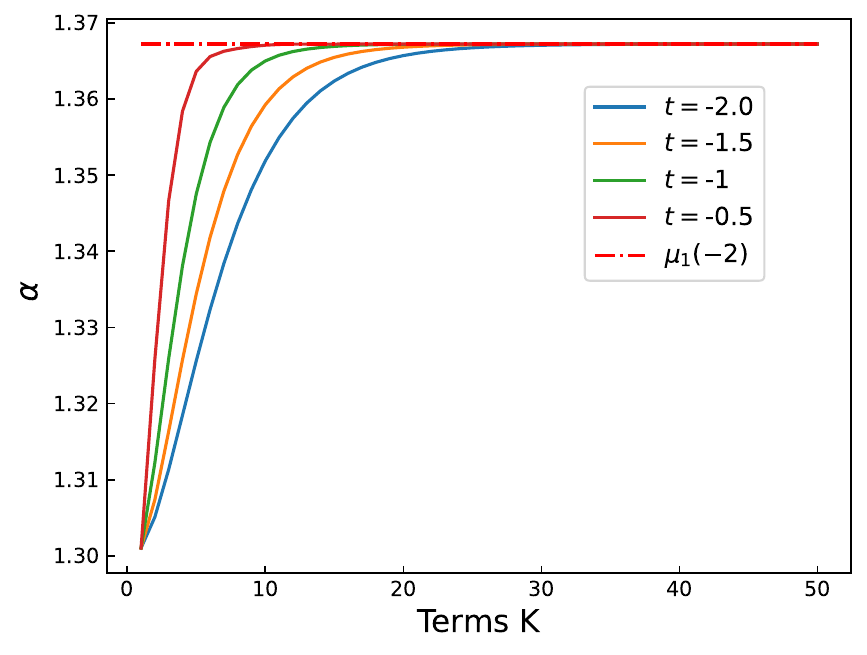} }
    \subfigure[$k=9$]{ \label{fig:ec12} \includegraphics[width=0.47\textwidth]{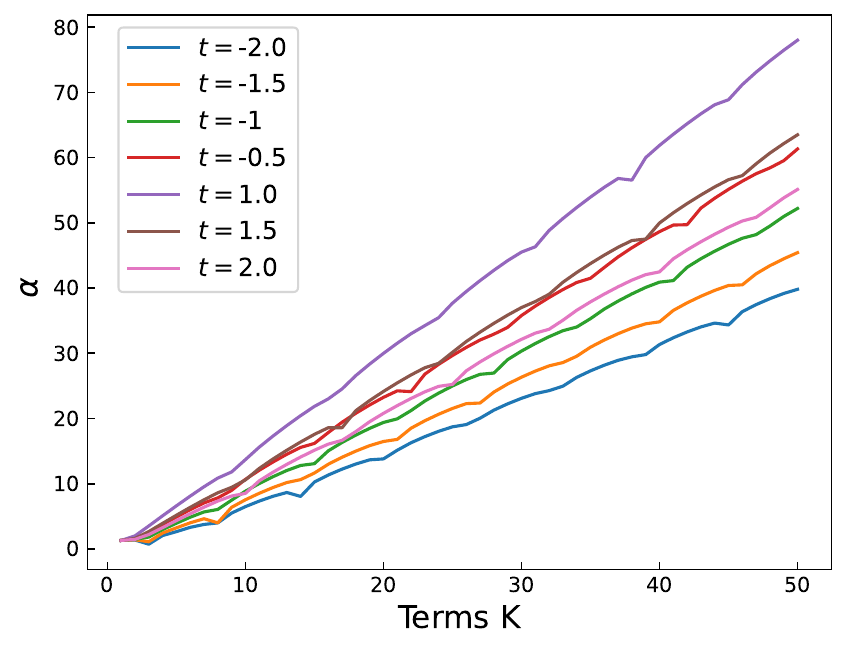} }
    \caption{The value of $\alpha $  of the Euler series    
    (\ref{perturbed_eigenvalue_zeta_almost_reg_graph_Euler}) versus parameter $t$ and the number of terms $K$ of graph $G_{(21,k)}$ with $\zeta=-2$.}
\end{figure}

We notice that the perturbation series (\ref{perturbed_eigenvalue_zeta_almost_reg_graph}) with $k=1$ and $\zeta=-2$ does not converge but fluctuates within a certain range as the number of terms $K$ increases in Figure. \ref{fig:c1}. 
Thus, we look into the convergence of its 
corresponding Euler series
(\ref{perturbed_eigenvalue_zeta_almost_reg_graph_Euler}) and
surprisingly find that the Euler series converges to the largest eigenvalue $\mu_1(-2)$ of Laplacian $ Q=\Delta+\zeta A$ in Fig. \ref{fig:ec11}, which implies that the Euler series (\ref{perturbed_eigenvalue_zeta_almost_reg_graph_Euler}) can provide a superior convergence range than the perturbation series (\ref{perturbed_eigenvalue_zeta_almost_reg_graph}) for some $\zeta$.
However, if the perturbation series (\ref{perturbed_eigenvalue_zeta_almost_reg_graph}) diverges with $k=9$ and $\zeta=-2$, then the Euler series
(\ref{perturbed_eigenvalue_zeta_almost_reg_graph_Euler}) also diverges, as shown in Figure. \ref{fig:ec12}, no matter what the parameter $t$ is.

\section{Conclusion}
\label{sec_conclusion}
This paper, divided into three parts, presents in the first part in Section \ref{sec_theory} an Euler $t$-series $\xi_{q;K}$ in (\ref{mu_perturbation_Eulersummation_general_t}) as a function of a unique node degree $d_q$ that approaches a Laplacian eigenvalue $\mu_k$ close to $d_q$ under certain convergence conditions. 
The second part in Section \ref{sec_performance_analysis} tries to unravel under which conditions the Euler series $\xi_{q;K}$ converges.
In general, we find that the Euler series $\xi_{q;K}$ converges for a high degree node $q$ and converges faster for larger degree differences of $d_q$ with the other degrees in the graph.
The third part in Section \ref{sec:almost_regular} discusses almost regular graphs, where only one node has a unique node degree while the remaining nodes all have the same degrees.
We provide the explicit perturbation Taylor series (\ref{perturbed_eigenvalue_zeta_almost_reg_graph})
and its Euler series $\xi_Q(\zeta)$ 
(\ref{perturbed_eigenvalue_zeta_almost_reg_graph_Euler}).
The numerical results show that the Euler series
(\ref{perturbed_eigenvalue_zeta_almost_reg_graph_Euler}) can converge even if the perturbation Taylor series (\ref{perturbed_eigenvalue_zeta_almost_reg_graph})
diverges for some $\zeta$.
However, we do not succeed in finding precise convergence conditions for the Euler series $\xi_{q;K}$ in (\ref{mu_perturbation_Eulersummation_general_t}) that hold for any given graph and any unique degree $d_q$. 
Nevertheless, the analytic closed-form approximation (\ref{lambda_Euler_upto_order4}) of a Laplacian eigenvalue close to a unique degree $d_q$ is deemed useful as an estimate and complements the many bounds for Laplacian eigenvalues in spectral graph theory.

\bigskip
{\bf Acknowledgment.} P. Van Mieghem has been funded by the
European Research Council (ERC) under the European Union's Horizon 2020
research and innovation programme (grant agreement No 101019718).

%%%%%%%%%%%%%%%%%%%%%%%%%%%%%%%%%%%%%%%%%%%%%%%%%%%%%%%%%%%%%

%%%%%%%%%%%%%%%%%%%%%%%%%%%%%%%%%%%%%%%%%%%%%%%%%%%%%%%%%%%%%

\appendix
\section*{Appendix}
\section{Computations of coefficients in almost regular graphs}
\label{sec_computation_coeff_almost_regular_graphs}
First,
\[
\beta_{1l}=\frac{a_{l1}}{x}%
\]
and
\[
c_{2}\left(  q\right)  =\sum_{l=1}^{N}\beta_{1,l}a_{1l}=\frac{1}{x}\sum
_{l=1}^{N}a_{l1}a_{1l}=\frac{d_{\max}}{x}%
\]

For $j=2$, the recursion (\ref{beta_jl_recursion_almost_regular_graph}) yields%
\[
\beta_{2l}=\frac{1}{x}\sum_{m=2}^{N}\beta_{1,m}a_{lm}=\frac{1}{x^{2}}%
\sum_{m=2}^{N}a_{m1}a_{lm}=\frac{1}{x^{2}}\sum_{m=1}^{N}a_{m1}a_{lm}%
\]
and%
\[
\beta_{2l}=\frac{\left(  A^{2}\right)  _{1l}}{x^{2}}%
\]
The corresponding coefficient $c_{j+1}\left(  q\right)  $ is%
\[
c_{3}\left(  q\right)  =\sum_{l=1}^{N}\beta_{2,l}a_{1l}=\frac{1}{x^{2}}%
\sum_{l=1}^{N}\left(  A^{2}\right)  _{1l}a_{1l}=\frac{\left(  A^{3}\right)
_{11}}{x^{2}}%
\]
The pseudo-recursion (\ref{pseudo_recursion_c_j}) becomes%
\[
c_{3}\left(  q\right)  =\frac{1}{x}\sum_{m=1}^{N}\beta_{1,m}\left(
A^{2}\right)  _{m1}=\frac{\left(  A^{3}\right)  _{11}}{x^{2}}%
\]

For $j=3$,%
\begin{align*}
    \beta_{3l}  
    &  =\frac{1}{x}\left(  \sum_{m=2}^{N}\beta_{2,m}a_{lm}-\sum
    _{m=1}^{N}a_{1m}\beta_{1l}\beta_{1,m}\right) \\
    &  =\frac{1}{x^{3}}\left(  \sum_{m=2}^{N}\left(  A^{2}\right)  _{1m}%
    a_{lm}-a_{l1}\sum_{m=1}^{N}a_{1m}\right) \\
    &  =\frac{1}{x^{3}}\left(  \sum_{m=1}^{N}\left(  A^{2}\right)  _{1m}%
    a_{lm}-\left(  A^{2}\right)  _{11}a_{l1}-a_{l1}d_{\max}\right)
\end{align*}
and%
\[
\beta_{3l}=\frac{1}{x^{3}}\left(  \left(  A^{3}\right)  _{1l}-2d_{\max}%
a_{l1}\right)
\]
The coefficient $c_{4}\left(  q\right)  $ is%
\[
c_{4}\left(  q\right)  =\sum_{l=1}^{N}\beta_{3,l}a_{1l}=\frac{1}{x^{3}}\left(
\sum_{l=1}^{N}\left(  A^{3}\right)  _{1l}a_{1l}-2d_{\max}\sum_{l=1}^{N}%
a_{l1}a_{1l}\right)  =\frac{1}{x^{3}}\left(  \left(  A^{4}\right)
_{11}-2d_{\max}^{2}\right)
\]
The pseudo-recursion (\ref{pseudo_recursion_c_j}) verifies%
\begin{align*}
    c_{4}\left(  q\right)   
    &  =\frac{1}{x}\left(  \sum_{m=1}^{N}\beta
    _{2,m}\left(  A^{2}\right)  _{m1}-\beta_{2,1}\left(  A^{2}\right)  _{11}%
    -c_{2}\left(  q\right)  c_{2}\left(  q\right)  \right) \\
    &  =\frac{1}{x^{3}}\left(  \sum_{m=1}^{N}\left(  A^{2}\right)  _{1m}\left(
    A^{2}\right)  _{m1}-\left(  A^{2}\right)  _{11}\left(  A^{2}\right)
    _{11}-d_{\max}^{2}\right)
\end{align*}

For $j=4$,%
\begin{align*}
    \beta_{4l}  
    &  =\frac{1}{x}\left(  \sum_{m=2}^{N}\beta_{3,m}a_{lm}-\sum
    _{m=1}^{N}a_{1m}\sum_{k=1}^{2}\beta_{kl}\beta_{3-k,m}\right) \\
    &  =\frac{1}{x}\left(  \sum_{m=2}^{N}\beta_{3,m}a_{lm}-\sum_{m=1}^{N}%
    a_{1m}\beta_{1l}\beta_{2,m}-\sum_{m=1}^{N}a_{1m}\beta_{2l}\beta_{1m}\right) \\
    &  =\frac{1}{x^{4}}\left(  \sum_{m=2}^{N}\left(  \left(  A^{3}\right)
    _{1m}-2a_{m1}d_{\max}\right)  a_{lm}-a_{l1}\sum_{m=1}^{N}a_{1m}\left(
    A^{2}\right)  _{1m}-\left(  A^{2}\right)  _{1l}d_{\max}\right) \\
    &  =\frac{1}{x^{4}}\left(  \sum_{m=2}^{N}a_{lm}\left(  A^{3}\right)
    _{1m}-2d_{\max}\left(  A^{2}\right)  _{1l}-a_{l1}\left(  A^{3}\right)
    _{11}-\left(  A^{2}\right)  _{1l}d_{\max}\right)
\end{align*}
and%
\[
\beta_{4l}=\frac{1}{x^{4}}\left(  \left(  A^{4}\right)  _{l1}-3d_{\max}\left(
A^{2}\right)  _{1l}-2\left(  A^{3}\right)  _{11}a_{l1}\right)
\]
The coefficient $c_{5}\left(  q\right)  $ is%
\begin{align*}
    c_{5}\left(  q\right)   
    &  =\sum_{l=1}^{N}a_{1l}\beta_{4,l}=\frac{1}{x^{4}%
    }\sum_{l=1}^{N}a_{1l}\left(  \left(  A^{4}\right)  _{l1}-3d_{\max}\left(
    A^{2}\right)  _{1l}-2\left(  A^{3}\right)  _{11}a_{l1}\right) \\
    &  =\frac{1}{x^{4}}\left(  \sum_{l=1}^{N}a_{1l}\left(  A^{4}\right)
    _{l1}-3d_{\max}\sum_{l=1}^{N}a_{1l}\left(  A^{2}\right)  _{1l}-2\left(
    A^{3}\right)  _{11}\sum_{l=1}^{N}a_{1l}a_{l1}\right) \\
    &  =\frac{1}{x^{4}}\left(  \left(  A^{5}\right)  _{11}-3d_{\max}\left(
    A^{3}\right)  _{11}-2\left(  A^{3}\right)  _{11}d_{\max}\right) \\
    &  =\frac{1}{x^{4}}\left(  \left(  A^{5}\right)  _{11}-5d_{\max}\left(
    A^{3}\right)  _{11}\right)
\end{align*}
The pseudo-recursion (\ref{pseudo_recursion_c_j}) is%
\begin{align*}
    c_{5}\left(  q\right)   
    &  =\frac{1}{x}\left(  \sum_{m=1}^{N}\beta
    _{3,m}\left(  A^{2}\right)  _{m1}-\beta_{3,1}\left(  A^{2}\right)  _{11}%
    -\sum_{k=1}^{2}c_{4-k}\left(  q\right)  c_{k+1}\left(  q\right)  \right) \\
    &  =\frac{1}{x^{4}}\left(  \sum_{m=1}^{N}\left(  \left(  \left(  A^{3}\right)
    _{1m}-2d_{\max}a_{1m}\right)  \right)  \left(  A^{2}\right)  _{m1}-3\left(
    A^{3}\right)  _{11}\left(  A^{2}\right)  _{11}\right) \\
    &  =\frac{1}{x^{4}}\left(  \left(  A^{5}\right)  _{11}-2d_{\max}\left(
    A^{3}\right)  _{11}-3\left(  A^{3}\right)  _{11}\left(  A^{2}\right)
    _{11}\right)
\end{align*}
but does not appear to be more efficient. Only, that we can compute one order
in $j$ higher.

For $j=5$, the recursion (\ref{beta_jl_recursion_almost_regular_graph})
becomes
\begin{align*}
\beta_{5l}  
    &  =\frac{1}{x}\left(  \sum_{m=2}^{N}\beta_{4,m}a_{lm}-\sum
    _{m=1}^{N}a_{1m}\sum_{k=1}^{3}\beta_{kl}\beta_{4-k,m}\right) \\
    &  =\frac{1}{x}\left(  \sum_{m=2}^{N}\beta_{4,m}a_{lm}-\beta_{1l}\sum
    _{m=1}^{N}a_{1m}\beta_{3,m}-\beta_{2l}\sum_{m=1}^{N}a_{1m}\beta_{2,m}%
    -\beta_{3l}\sum_{m=1}^{N}a_{1m}\beta_{1,m}\right) \\
    &  = \frac{1}{x} \left(\sum_{m=2}^{N}\left(  \frac{1}{x^{4}}\left(  \left(  A^{4}\right)
        _{m1}-3d_{\max}\left(  A^{2}\right)  _{1m}-2\left(  A^{3}\right)  _{11}%
        a_{m1}\right)  \right)  a_{lm} \right.\\
    &  \hspace{30pt} 
        -\frac{a_{l1}}{x}\sum_{m=1}^{N}a_{1m}\frac{1}{x^{3}}\left(  \left(
        A^{3}\right)  _{1m}-2d_{\max}a_{m1}\right)  -\frac{\left(  A^{2}\right)
        _{1l}}{x^{2}}\sum_{m=1}^{N}a_{1m}\frac{\left(  A^{2}\right)  _{1m}}{x^{2}}\\
    & \hspace{30pt}
        \left. -\frac{1}{x^{3}}\left(  \left(  A^{3}\right)  _{1l}-2d_{\max}a_{l1}\right)
        \sum_{m=1}^{N}a_{1m}\frac{a_{m1}}{x}\right)\\
    &  =\frac{1}{x^{5}}\left(
        \sum_{m=1}^{N}\left(  A^{4}\right)  _{m1}a_{lm}-3d_{\max}\sum_{m=1}^{N}\left(
        A^{2}\right)  _{1m}a_{lm}-2\left(  A^{3}\right)  _{11}\sum_{m=1}^{N}%
        a_{m1}a_{lm} \right.\\
    &  \hspace{30pt} 
        - \left(  \left(  A^{4}\right)  _{11}-3d_{\max}\left(  A^{2}\right)_{11}-2\left(  A^{3}\right)  _{11}a_{11} \right)  a_{l1}\\
    & \hspace{30pt} 
        -a_{l1}\sum_{m=1}^{N}a_{1m}\left(  A^{3}\right)  _{1m}+2a_{l1}d_{\max}%
        \sum_{m=1}^{N}a_{1m}-\left(  A^{2}\right)  _{1l}\sum_{m=1}^{N}a_{1m}\left(
        A^{2}\right)  _{1m}\\
    &  \hspace{30pt}  
        \left. -\left(  \left(  A^{3}\right)_{1l}-2d_{\max}a_{l1}\right)  \sum_{m=1}%
        ^{N}a_{1m} \right) \\
    &  =\frac{1}{x^{5}}\left(
        \left(  A^{5}\right)  _{1l}-3\left(  A^{2}\right)  _{11}\left(  A^{3}\right)
        _{1l}-2\left(  A^{3}\right)  _{11}\left(  A^{2}\right)  _{1l}-\left(
        A^{4}\right)  _{11}a_{l1}+3\left(  A^{2}\right)  _{11}^{2}a_{l1} \right.\\
    & \hspace{30pt} 
        \left.   -a_{l1}\left(  A^{4}\right) _{11}+2a_{l1}\left(  A^{2}\right)_{11}
        ^{2}-\left(  A^{2}\right)_{1l}\left(  A^{3}\right)_{11}-\left(
        A^{3}\right)  _{1l}d_{\max}+2d_{\max}^{2}a_{l1} \right) \\
    &  = \frac{1}{x^{5}}\left( 
        \left(  A^{5}\right)  _{1l}-4\left(  A^{2}\right)
        _{11}\left(  A^{3}\right)  _{1l}-3\left(  A^{3}\right)  _{11}\left(
        A^{2}\right)  _{1l}-2\left(  A^{4}\right)  _{11}a_{l1}+7\left(  A^{2}\right)
        _{11}^{2}a_{l1}  \right) 
\end{align*}

The coefficient $c_{6}\left(  q\right)  $ is
\begin{align*}
    c_{6}\left(  q\right)   
    &  =\sum_{l=1}^{N}\beta_{5,l}a_{1l}\\
    &  =\frac{1}{x^{5}}\sum_{l=1}^{N}\left(  \left(  A^{5}\right)  _{1l}-4\left(
    A^{2}\right)  _{11}\left(  A^{3}\right)  _{1l}-3\left(  A^{3}\right)
    _{11}\left(  A^{2}\right)  _{1l}-2\left(  A^{4}\right)  _{11}a_{l1}+7\left(
    A^{2}\right)  _{11}^{2}a_{l1}\right)  a_{1l}\\
    &  =\frac{1}{x^{5}}\left(  \left(  A^{6}\right)  _{11}-4\left(  A^{2}\right)
    _{11}\left(  A^{4}\right)  _{11}-3\left(  A^{3}\right)  _{11}^{2}-2\left(
    A^{4}\right)  _{11}\left(  A^{2}\right)  _{11}+7\left(  A^{2}\right)
    _{11}^{3}\right) \\
    &  =\frac{1}{x^{5}}\left(  \left(  A^{6}\right)  _{11}-6\left(  A^{2}\right)
    _{11}\left(  A^{4}\right)  _{11}-3\left(  A^{3}\right)  _{11}^{2}+7\left(
    A^{2}\right)  _{11}^{3}\right)
\end{align*}

For $j=6$, the recursion (\ref{beta_jl_recursion_almost_regular_graph})
becomes
\begin{align*}
    \beta_{6l}  
    &  =\frac{1}{x}\left(  \sum_{m=2}^{N} \beta_{5m} a_{lm}-\sum
    _{m=1}^{N}a_{1m}\sum_{k=1}^{4} \beta_{kl}\beta_{5-k,m}\right) \\
    &  =\frac{1}{x}\left(  
    \sum_{m=2}^{N} \beta_{5m} a_{lm} -\beta_{1l}\sum
    _{m=1}^{N}a_{1m}\beta_{4m} -\beta_{2l}\sum_{m=1}^{N}a_{1m}\beta_{3m} 
    \left. -\beta_{3l}\sum_{m=1}^{N}a_{1m}\beta_{2m} -\beta_{4l}\sum_{m=1}^{N}a_{1m}%
    \beta_{1m}\right) \right.\\
    &  =\frac{1}{x}\left(
        \frac{1}{x^{5}}\sum_{m=2}^{N} a_{lm}\left[ (A^{5})_{1m} -4(A^{2})_{11}
        (A^{3})_{1m}-3(A^{3})_{11}(A^{2})_{1m}  -2(A^{4})_{11}a_{m1}+7(A^{2})^{2} _{11}a_{m1}  \right]  \right.\\
    &  \hspace{30pt}
        -\frac{a_{l1}}{x^{5}}\sum_{m=1}^{N}a_{1m}\left[  (A^{4})_{m1} -3(A^{2}
        )_{11}(A^{2})_{1m}-2(A^{3})_{11}a_{m1}\right] \\
    &  \hspace{30pt}
        -\frac{(A^{2})_{1l}}{x^{5}} \sum_{m=1}^{N}a_{1m}\left[  (A^{3})_{1m}
        -2(A^{2})_{11}a_{m1}\right] 
        -\frac{ (A^{3})_{1l} -2(A^{2})_{11}a_{l1} }{x^{5}}\sum_{m=1}^{N}a_{1m}
        (A^{2})_{1m}\\
    &  \hspace{30pt}
        \left. -\frac{ (A^{4})_{l1} -3(A^{2})_{11}(A^{2})_{1l}-2(A^{3})_{11}a_{l1}}{x^{5}
        }\sum_{m=1}^{N}a_{1m}a_{m1}\right) \\
    &  =\frac{1}{x^{6}} \left(
        (A^{6})_{l1} -4(A^{4})_{l1}(A^{2})_{11}-3(A^{3})_{l1}(A^{3})_{11} 
        -2(A^{2})_{l1}(A^{4})_{11}+7(A^{2})_{l1}(A^{2})^{2}_{11} \right.\\
    &  \hspace{30pt}
        + \left[  -(A^{5})_{11}+7(A^{2})_{11}(A^{3})_{11} \right]  a_{l1}
        -\left[  (A^{5})_{11}-5(A^{2})_{11}(A^{3})_{11} \right]  a_{l1}\\
    &  \hspace{30pt}
        -(A^{2})_{1l} (A^{4})_{11} + 2(A^{2})_{1l}(A^{2})^{2}_{11}
        - (A^{3})_{1l}(A^{3})_{11} + 2(A^{2})_{11}(A^{3})_{11}a_{l1}\\
    & \hspace{30pt}
        \left.  - (A^{4})_{l1}(A^{2})_{1l} + 3(A^{2})_{1l}(A^{2})^{2}_{11} + 2(A^{2}
        )_{11}(A^{3})_{11}a_{l1}
    \right) \\
    &  =\frac{1}{x^{6}}\left(
        (A^{6})_{l1} -5(A^{4})_{l1}(A^{2})_{11}-4(A^{3})_{l1}(A^{3})_{11}
         -3(A^{2})_{l1}(A^{4})_{11}+12(A^{2})_{l1}(A^{2})^{2}_{11} \right.\\
    &  \hspace{30pt} 
        + \left.  \left[  -2(A^{5})_{11}+16(A^{2})_{11}(A^{3})_{11} \right]  a_{l1}\right)
\end{align*}

The coefficient $c_{7}\left(  q\right)  $ is
\begin{align*}
    c_{7}\left(  q \right)   
    &  =\sum_{l=1}^{N}\beta_{6l}a_{1l}\\
    &  =\frac{1}{x^{6}}\sum_{l=1}^{N}\left(
        (A^{6})_{l1} -5(A^{4})_{l1}(A^{2})_{11}-4(A^{3})_{l1}(A^{3})_{11}
        -3(A^{2})_{l1}(A^{4})_{11}+12(A^{2})_{l1}(A^{2})^{2}_{11} \right.\\
    &  \hspace{50pt} + \left.\left[  -2(A^{5})_{11}+16(A^{2})_{11}(A^{3})_{11} \right]  a_{l1}%
    \right)  a_{1l}\\
    &  =\frac{1}{x^{6}}\left(  (A^{7})_{11}-7(A^{5})_{11}(A^{2})_{11}%
    -7(A^{4})_{11}(A^{3})_{11} +28(A^{2})^{2}_{11}(A^{3})_{11}\right)
\end{align*}

\newpage
For $j=7$, the recursion (\ref{beta_jl_recursion_almost_regular_graph})
becomes
\begin{align*}
    \beta_{7l}  
    % ==1
    &  =\frac{1}{x}\left(  \sum_{m=2}^{N} \beta_{6m} a_{lm}-\sum
    _{m=1}^{N}a_{1m}\sum_{k=1}^{5} \beta_{kl}\beta_{6-k,m} \right) \\
    % ==2
    & =\frac{1}{x} \left(
         \sum_{m=2}^{N} \beta_{6m} a_{lm} -\beta_{1l}\sum
        _{m=1}^{N}a_{1m}\beta_{5m} -\beta_{2l}\sum_{m=1}^{N}a_{1m}\beta_{4m} \right.\\
    &  \hspace{30pt} 
        \left. -\beta_{3l}\sum_{m=1}^{N}a_{1m}\beta_{3m} -\beta_{4l}\sum_{m=1}^{N}a_{1m}%
        \beta_{2m} -\beta_{5l}\sum_{m=1}^{N}a_{1m}\beta_{1m} \right) \\
    % ==3
    &  =\frac{1}{x}\left(
        \frac{1}{x^{6}}\sum_{m=2}^{N} a_{lm} \left( (A^{6})_{m1} -5(A^{4})_{m1}(A^{2})_{11}-4(A^{3})_{m1}(A^{3})_{11} \right. \right.\\
    &    \hspace{90pt}     
        \left. -3(A^{2})_{m1}(A^{4})_{11}+12(A^{2})_{m1}(A^{2})^{2}_{11}  
            -2(A^{5})_{11}a_{m1}+16(A^{2})_{11}(A^{3})_{11} a_{m1} \right) \\
    &   \hspace{30pt} 
        -\frac{a_{l1}}{x^{6}}\sum_{m=1}^{N}a_{1m}\left[  
        (A^{5})_{1m} -4(A^{3})_{1m}(A^{2})_{11}-3(A^{2})_{1m}(A^{3})_{11} 
        -2(A^{4})_{11}a_{m1}+7(A^{2})^{2}_{11}a_{m1}     \right] \\
    &    \hspace{30pt} 
        -\frac{(A^{2})_{1l}}{x^{6}}\sum_{m=1}^{N}a_{1m}\left[  (A^{4})_{m1}
        -3(A^{2})_{1m}(A^{2})_{11}-2(A^{3})_{11}a_{m1}\right] \\
    &   \hspace{30pt} 
        -\frac{ (A^{3})_{1l} -2(A^{2})_{11}a_{l1} }{x^{6}}\sum_{m=1}^{N}a_{1m}\left[
        (A^{3})_{1m} -2(A^{2})_{11}a_{m1}\right] \\
    &    \hspace{30pt} 
        -\frac{ (A^{4})_{l1} -3(A^{2})_{1l}(A^{2})_{11}-2(A^{3})_{11}a_{l1}}{x^{6}%
        }\sum_{m=1}^{N}a_{1m}(A^{2})_{1m}\\
    &   \hspace{30pt} 
        \left. -\frac{(A^{5})_{1l} -4(A^{3})_{1l}(A^{2})_{11}-3(A^{2})_{1l}(A^{3}%
        )_{11}-2(A^{4})_{11}a_{l1}+7(A^{2})^{2}_{11}a_{l1}}{x^{6}}\sum_{m=1}^{N}%
        a_{1m}a_{m1}    \right)   \\  
    &  =\frac{1}{x^{7}}
        \left((A^{7})_{l1} -5(A^{5})_{l1}(A^{2})_{11}-4(A^{4})_{l1}(A^{3})_{11}
            +(A^{3})_{l1}\left[  -3(A^{4})_{11}+12(A^{2})^{2}_{11}\right] \right.\\
    &   \hspace{30pt} 
        +(A^{2})_{l1} \left[  -2(A^{5})_{11}+16(A^{2})_{11}(A^{3})_{11}\right] \\
    &   \hspace{30pt} 
        + \left[  -(A^{6})_{11}+8(A^{2})_{11}(A^{4})_{11}+4(A^{3})_{11}^{2}
        -12(A^{2})^{3}_{11} \right]  a_{l1} \\
    &   \hspace{30pt} 
        +\left[  -(A^{6})_{11} +6(A^{2})_{11}(A^{4})_{11}+3(A^{3})_{11}^{2}%
        -7(A^{2})^{3}_{11}\right]  a_{l1}\\
    &   \hspace{30pt} 
        -(A^{2})_{1l}(A^{5})_{11}+5(A^{2})_{1l}(A^{2})_{11}(A^{3})_{11}
        -(A^{3})_{1l}(A^{4})_{11} +2(A^{3})_{1l}(A^{2})_{11}^{2} +2(A^{2})_{11}
        (A^{4})_{11} a_{l1}\\
    &    \hspace{30pt} 
        -4(A^{2})_{11}^{3} a_{l1}
        -(A^{4})_{l1}(A^{3})_{11} +3(A^{2})_{1l}(A^{2})_{11}(A^{3})_{11} +2(A^{3}%
        )^{2}_{11}a_{l1}\\
    & \hspace{30pt} \left.
        -(A^{5})_{1l}(A^{2})_{11} +4(A^{3})_{1l}(A^{2})^{2}_{11}
        +3(A^{2})_{1l}(A^{2})_{11}(A^{3})_{11} + 2(A^{2})_{11}(A^{4})_{11}a_{l1}-7(A^{2})^{3} 
        _{11}a_{l1}  \right) \\
    % ==5
    &  =\frac{1}{x^{7}}\left(
        (A^{7})_{l1} -6(A^{5})_{l1}(A^{2})_{11}-5(A^{4})_{l1}(A^{3})_{11}
        -4(A^{3})_{l1}(A^{4})_{11}-3(A^{2})_{l1}(A^{5})_{11} +18(A^{3})_{1l}(A^{2})^{2}_{11}\right.\\
    &   \hspace{30pt}\left. 
        +27(A^{2})_{l1}(A^{2})_{11}(A^{3})_{11}
        + \left[  -2(A^{6})_{11}+18(A^{2})_{11}(A^{4})_{11}+9(A^{3})_{11}^{2}
        -30(A^{2})_{11}^{3}\right]  a_{l1}\right)
\end{align*}

The coefficient $c_{8}\left(  q\right)  $ is
\begin{align*}
    c_{8}\left(  q \right)   
    &  =\sum_{l=1}^{N}\beta_{7l}a_{1l}\\
    &  =\frac{1}{x^{7}}\sum_{l=1}^{N}\left(
        (A^{7})_{l1} -6(A^{5})_{l1}(A^{2})_{11}-5(A^{4})_{l1}(A^{3})_{11}
        -4(A^{3})_{l1}(A^{4})_{11}\right.\\
    & \hspace{50pt}   
        -3(A^{2})_{l1}(A^{5})_{11} +18(A^{3})_{1l}(A^{2})^{2}_{11}+27(A^{2})_{l1}(A^{2})_{11}(A^{3})_{11} \\
    & \hspace{50pt} 
        + \left.  \left[  -2(A^{6})_{11}+18(A^{2})_{11}(A^{4})_{11}+9(A^{3})_{11}^{2}
        -30(A^{2})_{11}^{3}\right]  a_{l1}  \right)  a_{1l}\\
    &  =\frac{1}{x^{7}}\left(  
        (A^{8})_{11}-8(A^{6})_{11}(A^{2})_{11}-8(A^{5})_{11}(A^{3})_{11}
        -4(A^{4})_{11}^{2} \right.\\
    & \hspace{30pt} 
        \left. +36(A^{2})_{11}^{2}(A^{4})_{11} +36(A^{3})^{2}_{11}(A^{2})_{11} -30(A^{2})_{11}^{4} \right)
\end{align*}

\newpage
For $j=8$, the recursion (\ref{beta_jl_recursion_almost_regular_graph})
becomes%
\begin{align*}
    \beta_{8l} 
%%%==1
    &  =\frac{1}{x}\left(  \sum_{m=2}^{N} \beta_{7m} a_{lm}-\sum
    _{m=1}^{N}a_{1m}\sum_{k=1}^{6} \beta_{kl}\beta_{7-k,m}\right) \\
%%%%%==2
    &  =\frac{1}{x} \left(
        \sum_{m=2}^{N} \beta_{7m} a_{lm} -\beta_{1l}\sum_{m=1}^{N}a_{1m}\beta_{6m}
        -\beta_{2l}\sum_{m=1}^{N}a_{1m}\beta_{5m} \right.\\
    &   \hspace{30pt} \left.-\beta_{3l}\sum_{m=1}^{N}a_{1m}\beta_{4m}  
        -\beta_{4l}\sum_{m=1}^{N}a_{1m}\beta_{3m}
        -\beta_{5l}\sum_{m=1}^{N}a_{1m}\beta_{2m} -\beta_{6l}\sum_{m=1}^{N}a_{1m}\beta_{1m} \right) \\
%%%%%==3
    &  =\frac{1}{x^{8}}\left(
        \sum_{m=2}^{N} a_{lm} \left( 
        (A^{7})_{m1} -6(A^{5})_{m1}(A^{2})_{11}-5(A^{4})_{m1}(A^{3})_{11}
        -4(A^{3})_{m1}(A^{4})_{11}-3(A^{2})_{m1}(A^{5})_{11} \right. \right.\\
    &   \hspace{90pt}    
        +18(A^{3})_{1m}(A^{2})^{2}_{11}+27(A^{2})_{m1}(A^{2})_{11}(A^{3})_{11} \\
    &   \hspace{90pt}
        +  \left. \left[  -2(A^{6})_{11}+18(A^{2})_{11}(A^{4})_{11}+9(A^{3})_{11}^{2}
        -30(A^{2})_{11}^{3}\right]  a_{m1}  \right) \\
    &   \hspace{30pt}
        -a_{l1}\sum_{m=1}^{N}a_{1m}\left((A^{6})_{m1} -5(A^{4})_{m1}(A^{2})_{11}
            -4(A^{3})_{m1}(A^{3})_{11}-3(A^{2})_{m1}(A^{4})_{11} \right.\\
    &   \hspace{110pt}
        \left.   +12(A^{2})_{m1}(A^{2})^{2}_{11}-2(A^{5})_{11}a_{m1}+16(A^{2})_{11}           (A^{3})_{11}a_{m1}\right) \\
    &   \hspace{30pt}
        -(A^{2})_{1l}\sum_{m=1}^{N}a_{1m}\left[ 
        (A^{5})_{1m} -4(A^{3})_{1m}(A^{2})_{11}
        -3(A^{2})_{1m}(A^{3})_{11} -2(A^{4})_{11}a_{m1}+7(A^{2})^{2}%
        _{11}a_{m1} \right] \\
    &    \hspace{30pt}
        - \left[  (A^{3})_{1l} -2(A^{2})_{11}a_{l1} \right]  \sum_{m=1}^{N}%
        a_{1m}\left[  (A^{4})_{m1} -3(A^{2})_{1m}(A^{2})_{11}-2(A^{3})_{11}%
        a_{m1}\right] \\
    &   \hspace{30pt}
        -\left[  (A^{4})_{l1} -3(A^{2})_{1l}(A^{2})_{11}-2(A^{3})_{11}a_{l1}\right]
        \sum_{m=1}^{N}a_{1m} \left[  (A^{3})_{1m} -2(A^{2})_{11}a_{m1}\right] \\
    &   \hspace{30pt}
        -\left[ 
         (A^{5})_{1l} -4(A^{3})_{1l}(A^{2})_{11}
         -3(A^{2})_{1l}(A^{3})_{11}-2(A^{4})_{11}a_{l1}+7(A^{2})^{2}_{11}a_{l1}\right]  \sum_{m=1}^{N}a_{1m}(A^{2})_{1m}\\
    &    \hspace{30pt}
        -\left(
        (A^{6})_{l1} -5(A^{4})_{l1}(A^{2})_{11}-4(A^{3})_{l1}(A^{3})_{11}
        -3(A^{2})_{l1}(A^{4})_{11} \right.\\
    &    \hspace{50pt}
        \left. \left.  +12(A^{2})_{l1}(A^{2})^{2}_{11} + \left[  -2(A^{5})_{11}+16(A^{2})_{11}%
        (A^{3})_{11} \right]  a_{l1} \right)  \sum_{m=1}^{N}a_{1m}a_{m1} \right)
\end{align*}

\begin{align*} 
    &  =\frac{1}{x^{8}}\left(
            (A^{8})_{l1} -6(A^{6})_{l1}(A^{2})_{11}-5(A^{5})_{l1}(A^{3})_{11}
            -4(A^{4})_{l1}(A^{4})_{11}-3(A^{3})_{l1}(A^{5})_{11} \right.\\
    &   \hspace{30pt}     
        +18(A^{4})_{l1}(A^{2})^{2}_{11}+27(A^{3})_{l1}(A^{2})_{11}(A^{3})_{11}\\
    &   \hspace{30pt}
        + (A^{2})_{l1}\left[  -2(A^{6})_{11}+18(A^{2})_{11}(A^{4})_{11}+9(A^{3})_{11}^{2}-30(A^{2})_{11}^{3}\right] \\
    &   \hspace{30pt}
        +a_{l1}\left[  -(A^{7})_{11} +9(A^{5})_{11}(A^{2})_{11}+9(A^{4})_{11}%
        (A^{3})_{11} -45(A^{2})_{11}^{2}(A^{3})_{11}\right]\\
    &    \hspace{30pt}
        +a_{l1}\left[  -(A^{7})_{11} +7(A^{5})_{11}(A^{2})_{11}+7(A^{4})_{11}%
        (A^{3})_{11}- 28(A^{2})_{11}^{2}(A^{3})_{11} \right] \\
    &   \hspace{30pt}
        -(A^{2})_{1l}(A^{6})_{11}+6(A^{2})_{1l}(A^{2})_{11}(A^{4})_{11}+3(A^{2}%
        )_{1l}(A^{3})_{11}^{2}-7(A^{2})_{1l}(A^{2})_{11}^{3}\\
    &   \hspace{30pt}
        - (A^{3})_{1l}(A^{5})_{11}+5(A^{3})_{1l}(A^{2})_{11}(A^{3})_{11}+\left[
        2(A^{2})(A^{5})_{11}-10 (A^{2})^{2}_{11}(A^{3})_{11}\right]  a_{l1}\\
    &   \hspace{30pt}
        -(A^{4})_{l1}(A^{4})_{11} +2(A^{4})_{l1}(A^{2})_{11}^{2}
        +3(A^{2})_{1l}(A^{2})_{11}(A^{4})_{11}   -6(A^{2})_{1l}(A^{2})_{11}^{3}\\
    &   \hspace{30pt}      
        +\left[  2(A^{3})_{11}(A^{4})_{11}-4(A^{2})_{11}^{2}(A^{3})_{11}\right]a_{l1}
        -(A^{5})_{1l}(A^{3})_{11} +4(A^{3})_{1l}(A^{2})_{11}(A^{3})_{11}+3(A^{2})_{1l}(A^{3})_{11}^{2}\\
    &   \hspace{30pt}
        +2(A^{3})_{11}(A^{4})_{11}a_{l1}-7(A^{2})^{2}_{11}(A^{3})_{11}a_{l1}
        -(A^{6})_{l1}(A^{2})_{11} +5(A^{4})_{l1}(A^{2})_{11}^{2}+4(A^{3})_{l1}
        (A^{2})_{11}(A^{3})_{11} \\
    &   \hspace{30pt}
        +3(A^{2})_{l1}(A^{2})_{11}(A^{4})_{11}
        \left. -12(A^{2})_{l1}(A^{2})^{3}_{11} + \left[  2(A^{2})_{11}(A^{5})_{11}%
        -16(A^{2})_{11}^{2}(A^{3})_{11} \right]  a_{l1} \right) \\
%%%%%%%======
    &   =\frac{1}{x^{8}}\left(
        (A^{8})_{l1} -7(A^{6})_{l1}(A^{2})_{11}-6(A^{5})_{l1}(A^{3})_{11}
        -5(A^{4})_{l1}(A^{4})_{11} -4(A^{3})_{l1}(A^{5})_{11}-3(A^{2})_{l1}(A^{6})_{11}\right.\\
    &   \hspace{30pt}
        +25(A^{4})_{1l}(A^{2})^{2}_{11}+40(A^{3})_{l1}(A^{2})_{11}(A^{3})_{11}\\
    &   \hspace{30pt}
        +30(A^{2})_{l1}(A^{2})_{11}(A^{4})_{11}+15(A^{2})_{l1}(A^{3})_{11}%
        ^{2}-55(A^{2})_{l1}(A^{2})_{11}^{3}\\
    &   \hspace{30pt}
        \left. + \left[  -2(A^{7})_{11}+20(A^{2})_{11}(A^{5})_{11}+20(A^{3})_{11}(A^{4}%
        )_{11}-110(A^{2})_{11}^{2}(A^{3})_{11}\right]  a_{l1} \right)
\end{align*}

The coefficient $c_{9}\left(  q\right)  $ is
\begin{align*}
    c_{9}\left(  q \right)   
    &  =\sum_{l=1}^{N}\beta_{8l}a_{1l}\\
    &  =\frac{1}{x^{8}}\sum_{l=1}^{N}\left(
        (A^{8})_{l1} -7(A^{6})_{l1}(A^{2})_{11}-6(A^{5})_{l1}(A^{3})_{11}
        -5(A^{4})_{l1}(A^{4})_{11}-4(A^{3})_{l1}(A^{5})_{11} \right.\\
    &  \hspace{50pt}
        -3(A^{2})_{l1}(A^{6})_{11} +25(A^{4})_{1l}(A^{2})^{2}_{11}+40(A^{3})_{l1}(A^{2})_{11}(A^{3})_{11} \\
    &  \hspace{50pt}
        +30(A^{2})_{l1}(A^{2})_{11}(A^{4})_{11}+15(A^{2})_{l1}(A^{3})_{11}^{2}-55(A^{2})_{l1}(A^{2})_{11}^{3}\\
    &   \hspace{50pt}
        \left. + \left[ 
        -2(A^{7})_{11}+20(A^{2})_{11}(A^{5})_{11}
        +20(A^{3})_{11}(A^{4})_{11}-110(A^{2})_{11}^{2}(A^{3})_{11}
        \right]  a_{l1} \right)  a_{1l}\\
    &   =\frac{1}{x^{8}} \left(
        (A^{9})_{11}-9(A^{7})_{11}(A^{2})_{11} -9(A^{6})_{11}(A^{3})_{11}%
        -9(A^{5})_{11}(A^{4})_{11} \right.\\
    &  \hspace{30pt}
        \left.  +45(A^{2})_{11}^{2}(A^{5})_{11}+90(A^{2})_{11}(A^{3})_{11}(A^{4})_{11}
        +15(A^{3})_{11}^{3} -165(A^{2})_{11}^{3}(A^{3})_{11}\right)
\end{align*}

\newpage
For $j=9$, the recursion (\ref{beta_jl_recursion_almost_regular_graph})
becomes%
\begin{align*}
    \beta_{9l}  
    &  =\frac{1}{x}\left(  \sum_{m=2}^{N}\beta_{8m}a_{lm}-\sum
    _{m=1}^{N}a_{1m}\sum_{k=1}^{7}\beta_{kl}\beta_{8-k,m}\right) \\
    &  =\frac{1}{x}\left(
        \sum_{m=2}^{N}\beta_{8m}a_{lm}-\beta_{1l}\sum_{m=1}^{N}a_{1m}\beta_{7m}%
        -\beta_{2l}\sum_{m=1}^{N}a_{1m}\beta_{6m}
        -\beta_{3l}\sum_{m=1}^{N}a_{1m}\beta_{5m} \right.\\
    &   \hspace{30pt} 
        \left.-\beta_{4l}\sum_{m=1}^{N}a_{1m}\beta_{4m}-              \beta_{5l}\sum_{m=1}^{N}a_{1m}
        \beta_{3m}-\beta_{6l}\sum_{m=1}^{N}a_{1m}\beta_{2m}-\beta_{7l}\sum_{m=1}
        ^{N}a_{1m}\beta_{1m} \right) \\
%%%%%%====3
    &  =\frac{1}{x^{9}}\left(
        \sum_{m=2}^{N}a_{lm}\left(
            (A^{8})_{m1}-7(A^{6})_{m1}(A^{2})_{11}-6(A^{5})_{m1}(A^{3})_{11}
            -5(A^{4})_{m1}(A^{4})_{11}-4(A^{3})_{m1}(A^{5})_{11} \right. \right.\\
    &   -3(A^{2})_{m1}(A^{6})_{11} +25(A^{4})_{1m}(A^{2})_{11}^{2}+40(A^{3})_{m1}                (A^{2})_{11}(A^{3})_{11}+30(A^{2})_{m1}(A^{2})_{11}(A^{4})_{11}\\
    &  +15(A^{2})_{m1}(A^{3})_{11}^{2}-55(A^{2})_{m1}(A^{2})_{11}^{3}
    +\left[  
        -2(A^{7})_{11}+20(A^{2})_{11}(A^{5})_{11}
        +20(A^{3})_{11}(A^{4})_{11} \right.\\
    &  \left.\left. -110(A^{2})_{11}^{2}(A^{3})_{11}
        \right]  a_{m1} \right) -a_{l1}\sum_{m=1}^{N}a_{1m}\left(
        (A^{7})_{m1}-6(A^{5})_{m1}(A^{2})_{11}-5(A^{4})_{m1}(A^{3})_{11}
        -4(A^{3})_{m1}(A^{4})_{11} \right.\\
    &   -3(A^{2})_{m1}(A^{5})_{11} 
        +18(A^{3})_{1m}(A^{2})_{11}^{2}
        +27(A^{2})_{m1}(A^{2})_{11}(A^{3})_{11}
        + \left[  -2(A^{6})_{11}+18(A^{2})_{11}(A^{4})_{11} \right.\\
    &   \left. \left. +9(A^{3})_{11}^{2}-30(A^{2})_{11}^{3}\right]  a_{m1}\right) 
    -(A^{2})_{1l}\sum_{m=1}^{N}a_{1m}\left(
        (A^{6})_{m1}-5(A^{4})_{m1}(A^{2})_{11}-4(A^{3})_{m1}(A^{3})_{11}  \right.\\
    &  \left. -3(A^{2})_{m1}(A^{4})_{11}+12(A^{2})_{m1}(A^{2})_{11}^{2} 
        -2(A^{5})_{11}a_{m1}+16(A^{2})_{11}(A^{3}%
        )_{11}a_{m1} \right)
         -\left[  (A^{3})_{1l}-2(A^{2})_{11}a_{l1}\right]\\
    &     \sum_{m=1}^{N}a_{1m}\left(
        (A^{5})_{1m}-4(A^{3})_{1m}(A^{2})_{11} -3(A^{2})_{1m}                        (A^{3})_{11}-2(A^{4})_{11}a_{m1}+7(A^{2})_{11}^{2}a_{m1}\right)\\
    &    -\left[  (A^{4})_{l1}-3(A^{2})_{1l}(A^{2})_{11}-2(A^{3})_{11}a_{l1}\right]  \sum_{m=1}^{N}a_{1m} \left[  (A^{4})_{m1}-3(A^{2})_{1m}(A^{2})_{11}-2(A^{3})_{11}a_{m1}\right]\\
    & -\left( (A^{5})_{1l}-4(A^{3})_{1l}(A^{2})_{11}-3(A^{2})_{1l}(A^{3})_{11}        -2(A^{4})_{11}a_{l1} +7(A^{2})_{11}^{2}a_{l1}\right)  \sum_{m=1}^{N}a_{1m}\left[  (A^{3})_{1m} -2(A^{2})_{11}a_{m1}\right]\\
    &   -\left(
        (A^{6})_{l1}-5(A^{4})_{l1}(A^{2})_{11}-4(A^{3})_{l1}(A^{3})_{11}
        -3(A^{2})_{l1}(A^{4})_{11}  +12(A^{2})_{l1}(A^{2})_{11}^{2} \right.\\
    &   \left.   +\left[  -2(A^{5})_{11}+16(A^{2})_{11}
        (A^{3})_{11}\right]  a_{l1} \right)  \sum_{m=1}^{N}a_{1m}(A^{2})_{1m}
        -\left(
        (A^{7})_{l1}-6(A^{5})_{l1}(A^{2})_{11}-5(A^{4})_{l1}(A^{3})_{11} \right.\\
    &   -4(A^{3})_{l1}(A^{4})_{11}-3(A^{2})_{l1}(A^{5})_{11} 
        +18(A^{3})_{1l}(A^{2})_{11}^{2}+27(A^{2})_{l1}(A^{2})_{11}(A^{3})_{11} \\
    &   \left.  \left.   +\left[  -2(A^{6})_{11}+18(A^{2})_{11}(A^{4})_{11} 
            +9(A^{3})_{11}^{2}-30(A^{2})_{11}^{3}\right]  a_{l1}%
        \right)  \sum_{m=1}^{N}a_{1m}a_{m1} \right)
\end{align*}

\begin{align*}
    &  =\frac{1}{x^{9}}\left(
            (A^{9})_{l1}-7(A^{7})_{l1}(A^{2})_{11}-6(A^{6})_{l1}(A^{3})_{11}
            -5(A^{5})_{l1}(A^{4})_{11}-4(A^{4})_{l1}(A^{5})_{11}-3(A^{3})_{l1}(A^{6})_{11} \right. \\
    &   \hspace{30pt}     +25(A^{5})_{1l}(A^{2})_{11}^{2}+40(A^{4})_{l1}(A^{2})_{11}(A^{3})_{11}
            +30(A^{3})_{l1}(A^{2})_{11}(A^{4})_{11}+15(A^{3})_{l1}(A^{3})_{11}
            ^{2}\\
    &   \hspace{30pt}      -55(A^{3})_{l1}(A^{2})_{11}^{3}+(A^{2})_{1l}\left[  -2(A^{7})_{11}+20(A^{2})_{11}(A^{5})_{11} 
            +20(A^{3})_{11}(A^{4})_{11}-110(A^{2})_{11}^{2}(A^{3})_{11}\right] \\
    &  \hspace{30pt}       +a_{l1}\left(
            -(A^{8})_{11}+10(A^{2})_{11}(A^{6})_{11}+10(A^{3})_{11}(A^{5})_{11}
            +5(A^{4})_{11}^{2} \right.\\
    &   \hspace{70pt}    \left.   \left.  -55(A^{4})_{11}(A^{2})_{11}^{2}
            -55(A^{3})_{11}^{2}(A^{2})_{11}+55(A^{2})_{11}^{4}%
            \right) \right) \\
    &   \hspace{30pt}     +a_{l1}\left(
            -(A^{8})_{11}+8(A^{6})_{11}(A^{2})_{11}
            +8(A^{5})_{11}(A^{3})_{11}+4(A^{4})_{11}^{2} \right.\\
    &    \hspace{70pt}    \left.     -36(A^{2})_{11}^{2}(A^{4})_{11}-36(A^{2})_{11}                 (A^{3})_{11}^{2}+30(A^{2}%
            )_{11}^{4} \right) \\
    &   \hspace{30pt}  +(A^{2})_{1l}\left[  -(A^{7})_{11}+7(A^{2})_{11}(A^{5})_{11}  +7(A^{3})_{11}               (A^{4})_{11}-28(A^{2})_{11}^{2}(A^{3})_{11}\right] \\
    &   \hspace{30pt}  +(A^{3})_{1l}\left[  -(A^{6})_{11}+6(A^{2})_{11}(A^{4})_{11}+3(A^{3})_{11}%
            ^{2}-7(A^{2})_{11}^{3}\right] \\
    &   \hspace{30pt}      +\left[  2(A^{2})_{11}(A^{6})_{11}-12(A^{2})_{11}^{2}(A^{4})_{11}-6(A^{2})_{11}        (A^{3})_{11}^{2}+14(A^{2})_{11}^{4}\right]  a_{l1} \\
    &   \hspace{30pt}  -(A^{4})_{l1}(A^{5})_{11}+5(A^{4})_{l1}(A^{2})_{11}(A^{3})_{11}
            +3(A^{2})_{l1}(A^{2})_{11}(A^{5})_{11}
            -15(A^{2})_{1l}(A^{2})_{11}^{2}(A^{3})_{11}\\
    &   \hspace{30pt}     
            +\left[  2(A^{3})_{11}(A^{5})_{11}-10(A^{2})_{11}(A^{3})_{11}^{2}\right]
            a_{l1} \\
    &   \hspace{30pt}    -(A^{5})_{1l}(A^{4})_{11}+2(A^{5})_{1l}(A^{2})_{11}^{2}+4(A^{3})_{1l}
            (A^{2})_{11}(A^{4})_{11}+3(A^{2})_{1l}(A^{3})_{11}(A^{4})_{11}\\
    &   \hspace{30pt}      -8(A^{3})_{1l}(A^{2})_{11}^{3}-6(A^{2})_{1l}(A^{2})_{11}^{2}(A^{3}%
            )_{11}
            +a_{l1}\left[  2(A^{4})_{11}^{2}-11(A^{2})_{11}^{2}(A^{4})_{11}%
            +14(A^{2})_{11}^{4}\right]\\
    &   \hspace{30pt}     -(A^{6})_{l1}(A^{3})_{11}+5(A^{4})_{l1}(A^{2})_{11}(A^{3})_{11}+4(A^{3}%
            )_{l1}(A^{3})_{11}^{2}
            +3(A^{2})_{l1}(A^{3})_{11}(A^{4})_{11}
            \\
    &   \hspace{30pt}      -12(A^{2})_{l1}(A^{2})_{11}^{2}(A^{3})_{11}+\left[  2(A^{3})_{11}(A^{5}%
            )_{11}-16(A^{2})_{11}(A^{3})_{11}^{2}\right]  a_{l1} \\
    &   \hspace{30pt}   -(A^{7})_{l1}(A^{2})_{11}+6(A^{5})_{l1}(A^{2})_{11}^{2}+5(A^{4})_{l1}%
            (A^{3})_{11}(A^{2})_{11}
            +4(A^{3})_{l1}(A^{2})_{11}(A^{4})_{11}\\
    &   \hspace{30pt}     +3(A^{2})_{l1}(A^{2})_{11}(A^{5})_{11} -18(A^{3})_{1l}(A^{2})_{11}^{3}-27(A^{2})_{l1}(A^{2})_{11}^{2}(A^{3})_{11} \\
    &   \hspace{30pt}      +\left[  2(A^{2})_{11}(A^{6})_{11}-18(A^{2})_{11}^{2}(A^{4})_{11}%
            -9(A^{2})_{11}(A^{3})_{11}^{2}+30(A^{2})_{11}^{4}\right]  a_{l1}  \\
    &  =\frac{1}{x^{9}}\left(
        (A^{9})_{l1}-8(A^{7})_{l1}(A^{2})_{11}-7(A^{6})_{l1}(A^{3})_{11}-6(A^{5}%
        )_{l1}(A^{4})_{11}-5(A^{4})_{l1}(A^{5})_{11} -4(A^{3})_{l1}(A^{6})_{11}\right.\\
    &   \hspace{30pt}  -3(A^{2})_{l1}(A^{7})_{11}
        +33(A^{5})_{1l}(A^{2})_{11}^{2}+55(A^{4})_{l1}(A^{2})_{11}(A^{3}%
        )_{11}+44(A^{3})_{l1}(A^{2})_{11}(A^{4})_{11}\\
    &   \hspace{30pt}  +22(A^{3})_{l1}(A^{3})_{11}^{2}+33(A^{2})_{l1}(A^{2})_{11}(A^{5}%
        )_{11}+33(A^{2})_{l1}(A^{3})_{11}(A^{4})_{11}
        -88(A^{3})_{l1}(A^{2})_{11}^{3}\\
    &   \hspace{30pt}  -198(A^{2})_{l1}(A^{2})_{11}^{2}(A^{3})_{11}+\left(
            -2(A^{8})_{11}+22(A^{6})_{11}(A^{2})_{11}+22(A^{5})_{11}(A^{3})_{11}%
            +11(A^{4})_{11}^{2} \right.\\
    & \hspace{30pt} \left. \left.    -132(A^{4})_{11}(A^{2})_{11}^{2}-132(A^{3})_{11}^{2}(A^{2})_{11}%
            +143(A^{2})_{11}^{4}
        \right)  a_{l1}  \right)
\end{align*}

\newpage
The coefficient $c_{10}\left(  q\right)  $ is
\begin{align*}
    c_{10}  &  =\sum_{l=1}^{N}\beta_{9l}a_{1l}\\
    &  =\frac{1}{x^{9}}\sum_{l=1}^{N}\left(
        (A^{9})_{l1}-8(A^{7})_{l1}(A^{2})_{11}-7(A^{6})_{l1}(A^{3})_{11}
        -6(A^{5})_{l1}(A^{4})_{11}-5(A^{4})_{l1}(A^{5})_{11} \right.\\
    &   \hspace{50pt} -4(A^{3})_{l1}(A^{6})_{11}-3(A^{2})_{l1}(A^{7})_{11}
        +33(A^{5})_{1l}(A^{2})_{11}^{2}
        +55(A^{4})_{l1}(A^{2})_{11}(A^{3})_{11}\\
    &    \hspace{50pt}
    +44(A^{3})_{l1}(A^{2})_{11}(A^{4})_{11} +22(A^{3})_{l1}(A^{3})_{11}^{2}+33(A^{2})_{l1}(A^{2})_{11}(A^{5})_{11}+33(A^{2})_{l1}(A^{3})_{11}(A^{4})_{11}\\
    &  \hspace{50pt}
    -88(A^{3})_{l1}(A^{2})_{11}^{3}-198(A^{2})_{l1}(A^{2})_{11}^{2}(A^{3})_{11}+
    \left(-2(A^{8})_{11}+22(A^{6})_{11}(A^{2})_{11} \right.\\
    &   \hspace{50pt} \left.\left.
    +22(A^{5})_{11}(A^{3})_{11}
            +11(A^{4})_{11}^{2}
            -132(A^{4})_{11}(A^{2})_{11}^{2}-132(A^{3})_{11}^{2}(A^{2})_{11}%
            +143(A^{2})_{11}^{4}%
        \right)  a_{l1} \right)  a_{1l}\\
    &  =\frac{1}{x^{9}}\left(
        (A^{10})_{11}-10(A^{8})_{11}(A^{2})_{11}-10(A^{7})_{11}(A^{3})_{11}%
        -10(A^{6})_{11}(A^{4})_{11}
        -5(A^{5})_{11}^{2}
        +55(A^{6})_{11}(A^{2})_{11}^{2} \right.\\
    &  \hspace{30pt} +110(A^{5})_{11}(A^{2})_{11}(A^{3})_{11}
        +55(A^{4})_{11}^{2}(A^{2})_{11}+55(A^{4})_{11}(A^{3})_{11}^{2}\\
    &   \hspace{30pt}\left.  -220(A^{4})_{11}(A^{2})_{11}^{3}-330(A^{3})_{11}^{2}(A^{2})_{11}^{2}%
        +143(A^{2})_{11}^{5}\right)
\end{align*}

Summary.
\begin{align*}
\beta_{1l}  &  =\frac{a_{l1}}{x}\\
\beta_{2l}  &  =\frac{(A^{2})_{1l}}{x^{2}}\\
\beta_{3l}  &  =\frac{(A^{3})_{1l}-2(A^{2})_{11}a_{l1}}{x^{3}}\\
\beta_{4l}  &  =\frac{(A^{4})_{l1}-3(A^{2})_{1l}(A^{2})_{11}-2(A^{3}%
)_{11}a_{l1}}{x^{4}}\\
\beta_{5l}  &  =\frac{1}{x^{5}}
\left((A^{5})_{1l}-4(A^{3})_{1l}(A^{2})_{11}-3(A^{2})_{1l}(A^{3})_{11}-2(A^{4})_{11}a_{l1}+7(A^{2})_{11}^{2}a_{l1}\right)\\
\beta_{6l}  &  =\frac{1}{x^6} \left(
(A^{6})_{l1}-5(A^{4})_{l1}(A^{2})_{11}-4(A^{3})_{l1}(A^{3})_{11}-3(A^{2})_{l1}(A^{4})_{11}+12(A^{2})_{l1}(A^{2})_{11}^{2} \right.\\
& \hspace{30pt} + \left. \left[  -2(A^{5})_{11}+16(A^{3})_{11}(A^{2})_{11}\right]  a_{l1} \right)\\
\beta_{7l}  &  =\frac{1}{x^{7}}\left(
        (A^{7})_{l1}-6(A^{5})_{l1}(A^{2})_{11}-5(A^{4})_{l1}(A^{3})_{11}-4(A^{3}%
        )_{l1}(A^{4})_{11}-3(A^{2})_{l1}(A^{5})_{11}
        +18(A^{3})_{1l}(A^{2})_{11} ^{2} \right. \\
&   \hspace{30pt} \left.  +27(A^{2})_{l1}(A^{2})_{11}(A^{3})_{11}+\left[  -2(A^{6})_{11}+18(A^{4})_{11}(A^{2})_{11}+9(A^{3})_{11}^{2}-30(A^{2})_{11}^{3}\right]  a_{l1}\right) \\
\beta_{8l}  &  =\frac{1}{x^{8}}\left(
        (A^{8})_{l1}-7(A^{6})_{l1}(A^{2})_{11}-6(A^{5})_{l1}(A^{3})_{11}-5(A^{4}%
        )_{l1}(A^{4})_{11}-4(A^{3})_{l1}(A^{5})_{11}-3(A^{2})_{l1}(A^{6})_{11} \right.\\
    &   +25(A^{4})_{1l}(A^{2})_{11}^{2}+40(A^{3})_{l1}(A^{2})_{11}(A^{3})_{11}+30(A^{2})_{l1}(A^{2})_{11}(A^{4})_{11} +15(A^{2})_{l1}(A^{3})_{11}^{2}-55(A^{2})_{l1}(A^{2})_{11}^{3}\\
    &    +\left. \left[ -2(A^{7})_{11}+20(A^{5})_{11}(A^{2})_{11}+20(A^{4})_{11}(A^{3})_{11}-110(A^{3})_{11}(A^{2})_{11}^{2}\right]  a_{l1}    \right) \\
\beta_{9l}  &  =\frac{1}{x^{9}}\left(
        (A^{9})_{l1}-8(A^{7})_{l1}(A^{2})_{11}-7(A^{6})_{l1}(A^{3})_{11}-6(A^{5}%
        )_{l1}(A^{4})_{11}-5(A^{4})_{l1}(A^{5})_{11}
        -4(A^{3})_{l1}(A^{6})_{11} \right.\\
    &   \hspace{30pt} -3(A^{2})_{l1}(A^{7})_{11}+33(A^{5})_{1l}(A^{2})_{11}^{2}+55(A^{4})_{l1}(A^{2})_{11}(A^{3}%
        )_{11}+44(A^{3})_{l1}(A^{2})_{11}(A^{4})_{11}\\
    & \hspace{30pt}
    +22(A^{3})_{l1}(A^{3})_{11}^{2}
        +33(A^{2})_{l1}(A^{2})_{11}(A^{5})_{11}+33(A^{2})_{l1}(A^{3})_{11}(A^{4}%
        )_{11}\\
    & \hspace{30pt}-88(A^{3})_{l1}(A^{2})_{11}^{3}-198(A^{2})_{l1}(A^{2})_{11}^{2}(A^{3})_{11}\\
    &   \hspace{30pt}
    +\left[
            -2(A^{8})_{11}+22(A^{6})_{11}(A^{2})_{11}+22(A^{5})_{11}(A^{3} )_{11}%
            +11(A^{4})_{11}^{2} \right.\\
    &   \hspace{50pt} \left. \left.    
    -132(A^{4})_{11}(A^{2})_{11}^{2}-132(A^{3})_{11}^{2}(A^{2})_{11} +143(A^{2})_{11}^{4}\right]  a_{l1}
    \right)
\end{align*}

\begin{align*}
c_{2}  &  =\frac{1}{x}\left(  A^{2}\right)  _{11}\\
c_{3}  &  =\frac{1}{x^{2}}\left(  A^{3}\right)  _{11}\\
c_{4}  &  =\frac{1}{x^{3}}\left(  \left(  A^{4}\right)  _{11}-2\left(
A^{2}\right)  _{11}^{2}\right) \\
c_{5}  &  =\frac{1}{x^{4}}\left(  \left(  A^{5}\right)  _{11}-5\left(
A^{2}\right)  _{11}\left(  A^{3}\right)  _{11}\right) \\
c_{6}  &  =\frac{1}{x^{5}}\left(  \left(  A^{6}\right)  _{11}-6\left(
A^{2}\right)  _{11}\left(  A^{4}\right)  _{11}-3\left(  A^{3}\right)
_{11}^{2}+7\left(  A^{2}\right)  _{11}^{3}\right) \\
c_{7}  &  =\frac{1}{x^{6}}\left(  (A^{7})_{11}-7(A^{5})_{11}(A^{2}%
)_{11}-7(A^{4})_{11}(A^{3})_{11}+28(A^{3})_{11}(A^{2})_{11}^{2}\right) \\
c_{8}  &  =\frac{1}{x^{7}}\left(  (A^{8})_{11}-8(A^{6})_{11}(A^{2})_{11}-8(A^{5})_{11}(A^{3})_{11}-4(A^{4})_{11}^{2}+36(A^{4})_{11}(A^{2}
)_{11}^{2} \right.\\
& \hspace{30pt} \left. +36(A^{3})_{11}^{2}(A^{2})_{11}-30(A^{2})_{11}^{4}\right) \\
c_{9}  &  =\frac{1}{x^{8}}\left(
        (A^{9})_{11}-9(A^{7})_{11}(A^{2})_{11}-9(A^{6})_{11}(A^{3})_{11}-9(A^{5}%
        )_{11}(A^{4})_{11}+ \right.\\
    &\hspace{30pt} \left.    45(A^{5})_{11}(A^{2})_{11}^{2}+90(A^{4})_{11}(A^{2})_{11}(A^{3})_{11}+15(A^{3})_{11}^{3}-165(A^{3})_{11}(A^{2})_{11}^{3}\right) \\
c_{10}  &  =\frac{1}{x^{9}}\left(
        (A^{10})_{11}-10(A^{8})_{11}(A^{2})_{11}-10(A^{7})_{11}(A^{3})_{11}%
        -10(A^{6})_{11}(A^{4})_{11}-5(A^{5})_{11}^{2} \right.\\
    &  \hspace{30pt}  +55(A^{6})_{11}(A^{2})_{11}^{2}+110(A^{5})_{11}(A^{2})_{11}(A^{3}%
        )_{11}+55(A^{4})_{11}^{2}(A^{2})_{11}+55(A^{4})_{11}(A^{3})_{11}^{2}\\
    &   \hspace{30pt} \left. -220(A^{4})_{11}(A^{2})_{11}^{3}-330(A^{3})_{11}^{2}(A^{2})_{11}^{2}%
        +143(A^{2})_{11}^{5} \right)
\end{align*}

\section{The recursion for the characteristic coefficient $\mathcal{A}%
[k,m]$}
\label{sec_ recursion_characteristic_coefficient_A}

The definition of the characteristic coefficient is
\begin{align*}
\mathcal{A}[k,m]  &  =\sum_{\sum_{i=1}^{k}j_{i}=m;j_{i}>0}\prod_{i=1}%
^{k}\left(  A^{j_{i}}\right)  _{11}\\
&  =\sum_{\sum_{i=1}^{k-1}j_{i}+j_{k}=m;j_{i}>0}\left(  A^{j_{k}}\right)
_{11}\prod_{i=1}^{k-1}\left(  A^{j_{i}}\right)  _{11}%
\end{align*}
Since $1\leq j_{k}\leq m-\sum_{i=1}^{k-1}j_{i}$ and since $j_{1}\geq1$, it
holds that $1\leq j_{k}\leq m-\left(  k-1\right)  $ and we can write the sum,
after replacing $j_{k}$ by $j$, as%
\[
\mathcal{A}[k,m]=\sum_{j=1}^{m-k+1}\left(  A^{j}\right)  _{11}\sum_{\sum
_{i=1}^{k-1}j_{i}=m-j;j_{i}>0}\prod_{i=1}^{k-1}\left(  A^{j_{i}}\right)  _{11}%
\]
and we obtain the characteristic coefficient recursion, for $k>1$,%
\begin{equation}
\mathcal{A}[k,m]=\sum_{j=1}^{m-k+1}\left(  A^{j}\right)  _{11}\mathcal{A}%
[k-1,m-j] \label{chc_A_recursion}%
\end{equation}
with initial starting value $\mathcal{A}[1,m]=\left(  A^{m}\right)  _{11}$.
The recursion (\ref{chc_A_recursion}) enables to compute the characteristic
coefficient $\mathcal{A}[k,m]$ for any adjacency matrix $A$.

The characteristic coefficients of a complex function $f(z)$ with Taylor series
$f(z) = \sum_{k=0}^{\infty}f_{k}(z_{0})(z-z_{0})^{k}$ are defined as \cite{PVM_charcoef}
\begin{equation*}
\left.  s[k,m]\right\vert _{ f\left(  z\right)  }(z_{0})=\sum_{\sum_{i=1}^{k}j_{i}=m;j_{i}>0}\prod_{i=1}%
^{k} f_{j_{i}}(z_{0})
\end{equation*}
and obey $\left. s[k,m]\right\vert _{f\left(  z\right)  }(z_{0})= 0$ if $k<0$ and $k>m$.
Moreover, $\left. s[k,m]\right\vert _{f\left(  z\right)  }(z_{0})$ possesses a recursion form and has been derived from function theory in \cite[Eq. (3), in Section
2]{PVM_ASYM}
\begin{align}
\left.  s[1,m]\right\vert _{ f\left(  z\right)  }(z_{0}  &  =f_{m}(z_{0})\nonumber\\
\left.  s[k,m]\right\vert _{ f\left(  z\right)  }(z_{0} &  =\sum_{j=1}^{m-k+1}f_{j}\;s[k-1,m-j]\hspace{2cm}(k>1)
\label{s_recursive}%
\end{align}
While (\ref{chc_A_recursion}) presents the combinational derivation.
The differentiation formula of the characteristic coefficient $\left. s[k,m]\right\vert _{f\left(  z\right)  }(z_{0})$ in \cite{PVM_charcoef} is
\begin{equation}
\left.  s[k,m]\right\vert _{f\left(  z\right)  }(z_{0})=\frac{1}{2\pi i}%
\int_{C(z_{0})}\frac{[f(z)-f(z_{0})]^{k}}{(z-z_{0})^{m+1}}\,dz=\frac{1}%
{m!}\,\left.  \frac{d^{m}}{dz^{m}}[f(z)-f(z_{0})]^{k}\right\vert _{z=z_{0}}.
\label{s_general}%
\end{equation}
The analogy with chc in (\ref{chc_A_recursion}) suggests that $f_{k}(z_{0})=\left(  A^{k}\right) _{11}$ and the corresponding generating function is 
$f_{G}\left(  z\right)  =\sum_{m=0}^{\infty}\left(  A^{m}\right)  _{11}z^{m}$
for a graph $G$ and $z_{0}=0$.

\section{The characteristic coefficient $\mathcal{A}[k,m]$ of the
complete graph}

The maximum value of $\left(  A^{m}\right)  _{11}$ is attained in the complete
graph $K_{N}$ with adjacency matrix $A=J-I$. The maximum number of closed
walks of lenght $m$ starting and ending at node $j$ is $\left(  J-I\right)
_{jj}^{m}=\frac{1}{N}\left(  \left(  N-1\right)  ^{m}-\left(  -1\right)
^{m}\right)  +\left(  -1\right)  ^{m}$, as shown in \cite[eq. (2.18), p.
30]{PVM_graphspectra_second_edition}. The corresponding generating function
for the graph $G=$ $K_{N}$ is, for $\left\vert z\right\vert <\frac{1}{N-1}$,%
\begin{align*}
f_{K_{N}}\left(  z\right)   &  =\sum_{m=0}^{\infty}\left(  J-I\right)
_{jj}^{m}z^{m}=\frac{1}{N}\sum_{m=0}^{\infty}\left(  N-1\right)  ^{m}%
z^{m}+\left(  1-\frac{1}{N}\right)  \sum_{m=0}^{\infty}\left(  -1\right)
^{m}z^{m}\\
&  =\frac{1}{N}\frac{1}{1-\left(  N-1\right)  z}+\left(  1-\frac{1}{N}\right)
\frac{1}{1+z}%
\end{align*}
from which $f_{K_{N}}\left(  0\right)  =1=\left(  A^{0}\right)  _{11}$,
$\left(  A^{1}\right)  _{11}=0$, $\left(  A^{2}\right)  _{11}=N-1$, etc. The
characteristic coefficient $\mathcal{A}[k,m]$ for the complete graph $K_{N}$
around $z_{0}=0$ can be computed from the differentiation formula in (\ref{s_general})
with
\begin{align*}
f_{K_{N}}(z)-f_{K_{N}}(z_{0})  &  =\frac{1}{N}\frac{1}{1-\left(  N-1\right)  z}+\frac{1}%
{N}\frac{N-1}{1+z}-1\\
&  =\frac{\left(  N-1\right)  z^{2}}{\left(  1-\left(  N-1\right)  z\right)
\left(  1+z\right)  }%
\end{align*}
The contour integral in (\ref{s_general}) becomes%
\begin{align*}
\mathcal{A}[k,m]  &  =\frac{\left(  N-1\right)  ^{k}}{2\pi i}\int_{C(0)}%
\frac{z^{2k}}{z^{m+1}\left(  1-\left(  N-1\right)  z\right)  ^{k}\left(
1+z\right)  ^{k}}\,dz\\
&  =\frac{\left(  N-1\right)  ^{k}}{2\pi i}\int_{C(0)}\frac{\left(  1-\left(
N-1\right)  z\right)  ^{-k}\left(  1+z\right)  ^{-k}dz}{z^{m-2k+1}}\\
&  =\frac{\left(  N-1\right)  ^{k}}{\left(  m-2k\right)  !}\frac{d^{m-2k}%
}{dz^{m-2k}}\,\left.  \left(  \left(  1-\left(  N-1\right)  z\right)
^{-k}\left(  1+z\right)  ^{-k}\right)  \right\vert _{z=0}%
\end{align*}
The Taylor series around $z_{0}=0$ of $\left(  1-\left(  N-1\right)  z\right)
^{-k}=\sum_{n=0}^{\infty}\binom{-k}{n}\left(  -\left(  N-1\right)  \right)
^{n}z^{n}$ and of $\left(  1+z\right)  ^{-k}=\sum_{n=0}^{\infty}\binom{-k}%
{n}z^{n}$ indicates that the Cauchy product is%
\begin{align*}
\left(  1-\left(  N-1\right)  z\right)  ^{-k}\left(  1+z\right)  ^{-k}  &
=\sum_{n=0}^{\infty}\binom{-k}{n}\left(  -\left(  N-1\right)  \right)
^{n}z^{n}\sum_{n=0}^{\infty}\binom{-k}{n}z^{n}\\
&  =\sum_{n=0}^{\infty}\left(  \sum_{r=0}^{n}\binom{-k}{r}\binom{-k}%
{n-r}\left(  -\left(  N-1\right)  \right)  ^{r}\right)  z^{n}%
\end{align*}
from which
\[
\mathcal{A}[k,m]=\left(  N-1\right)  ^{k}\sum_{r=0}^{m-2k}\binom{-k}{r}%
\binom{-k}{m-2k-r}\left(  -1\right)  ^{r}\left(  N-1\right)  ^{r}%
\]
Using $\binom{s+j-1}{j}=(-1)^{j}{\binom{-s}{j}}$ in
(\ref{binomial_min_z_in_positive_z}),
we arrive at%
\begin{equation}
\mathcal{A}[k,m]=(-1)^{m}\left(  N-1\right)  ^{k}\sum_{r=0}^{m-2k}%
\binom{k+r-1}{r}\binom{m-k-r-1}{m-2k-r}\left(  -1\right)  ^{r}\left(
N-1\right)  ^{r} \label{chc_complete_graph}%
\end{equation}
which is an alternating sum in $d_{\max}=N-1$ and, unfortunately, difficult to
simplify further. From $s\left[  1,m\right]  =f_{m}$, we have
\begin{align*}
\mathcal{A}[1,m]  &  =\left(  -1\right)  ^{m}\left(  1-\frac{1-\left(
-d_{\max}\right)  ^{m}}{1-\left(  -d_{\max}\right)  }\right)  =\left(
-1\right)  ^{m}\left(  1-\sum_{j=0}^{m-1}\left(  -d_{\max}\right)  ^{j}\right)
\\
&  =\left(  -1\right)  ^{m-1}\sum_{j=1}^{m-1}\left(  -d_{\max}\right)
^{j}=\sum_{j=1}^{m-1}\left(  -1\right)  ^{m-1-j}d_{\max}^{j}%
\end{align*}
We also observe that $\mathcal{A}[k,m]=0$ if $k>\frac{m}{2}$.

\section{ An upper bound for the characteristic coefficient
$\mathcal{A}[k,m]$ of the complete graph}

An upper bound for the characteristic coefficient $\mathcal{A}[k,m]$ in
(\ref{chc_adjacency_closed_walks_node_q}) can be deduced from the contour
integral in (\ref{s_general}). We evaluate the integral along a circle
$z=re^{i\theta}$ around the origin with $r<\frac{1}{N-1}$ (to avoid encircling
the pole of $\frac{1}{1-\left(  N-1\right)  z}$ at $z=\frac{1}{N-1}$)%
\begin{align*}
\mathcal{A}[k,m]  &  =\frac{\left(  N-1\right)  ^{k}}{2\pi i}\int_{C(0)}%
\frac{z^{2k-m-1}}{\left(  1-\left(  N-1\right)  z\right)  ^{k}\left(
1+z\right)  ^{k}}\,dz\\
&  =\frac{\left(  N-1\right)  ^{k}}{2\pi}\int_{0}^{2\pi}\frac{\left(
re^{i\theta}\right)  ^{2k-m}}{\left(  1-\left(  N-1\right)  re^{i\theta
}\right)  ^{k}\left(  1+re^{i\theta}\right)  ^{k}}\,d\theta
\end{align*}
and%
\[
\mathcal{A}[k,m]=\frac{r^{2k-m}\left(  N-1\right)  ^{k}}{2\pi}\int_{0}^{2\pi
}\frac{e^{i\theta\left(  2k-m\right)  }}{\left(  1-\left(  N-1\right)
re^{i\theta}\right)  ^{k}\left(  1+re^{i\theta}\right)  ^{k}}\,d\theta
\]
We bound the integral as%
\begin{align*}
\left\vert \mathcal{A}[k,m]\right\vert  &  \leq\frac{r^{2k-m}\left(
N-1\right)  ^{k}}{2\pi}\int_{0}^{2\pi}\frac{d\theta}{\left\vert 1-\left(
N-1\right)  re^{i\theta}\right\vert ^{k}\left\vert 1+re^{i\theta}\right\vert
^{k}}\,\\
&  \leq\frac{r^{2k-m}\left(  N-1\right)  ^{k}}{2\pi}\int_{0}^{2\pi}%
\frac{d\theta}{\left(  1-\left(  N-1\right)  r\right)  ^{k}\left(  1-r\right)
^{k}}\,
\end{align*}
and obtain%
\[
\left\vert \mathcal{A}[k,m]\right\vert \leq\frac{r^{2k-m}\left(  N-1\right)
^{k}}{\left(  1-\left(  N-1\right)  r\right)  ^{k}\left(  1-r\right)  ^{k}%
}=\frac{1}{r^{m}}\left(  \frac{r^{2}\left(  N-1\right)  }{\left(  1-\left(
N-1\right)  r\right)  \left(  1-r\right)  }\right)  ^{k}%
\]
which holds for any $r<\frac{1}{N-1}$. Let $r=\frac{x}{N-1}$ for $0<x<1$, then
the upper bound becomes%
\[
\left\vert \mathcal{A}[k,m]\right\vert \leq\frac{\left(  N-1\right)  ^{m}%
}{x^{m}}\left(  \frac{x^{2}\left(  N-1\right)  ^{-1}}{\left(  1-x\right)
\left(  1-\frac{x}{N-1}\right)  }\right)  ^{k}=\frac{\left(  N-1\right)  ^{m}%
}{\left(  N-1-x\right)  ^{k}}\frac{x^{2k-m}}{\left(  1-x\right)  ^{k}}%
\]
For example, for $x=\frac{1}{2}$,%
\[
\left\vert \mathcal{A}[k,m]\right\vert \leq2^{m-k}\frac{\left(  N-1\right)
^{m}}{\left(  N-\frac{3}{2}\right)  ^{k}}%
\]
The minimizer of $\frac{x^{2k-m}}{\left(  1-x\right)  ^{k}}$ occurs at the
root of $\frac{d}{dx}\left(  \frac{x^{2k-m}}{\left(  1-x\right)  ^{k}}\right)
=0$, which is $x=\frac{m-2k}{m-k}$. Hence, the minimum is $\frac{\left(
\frac{m-2k}{m-k}\right)  ^{2k-m}}{\left(  \frac{k}{m-k}\right)  ^{k}}%
=\frac{\left(  m-k\right)  ^{m-k}}{\left(  m-2k\right)  ^{m-2k}k^{k}}$ and we
find%
\[
\left\vert \mathcal{A}[k,m]\right\vert \leq\frac{\left(  N-1\right)  ^{m}%
}{\left(  N-1-\frac{m-2k}{m-k}\right)  ^{k}}\frac{\left(  m-k\right)  ^{m-k}%
}{\left(  m-2k\right)  ^{m-2k}k^{k}}%
\]
which\footnote{Stirling's approximation $\Gamma\left(  z+1\right)
\approx\sqrt{2\pi z}z^{z}e^{-z}$ yields%
\[
\frac{\left(  m-k\right)  !}{k!\left(  m-2k\right)  !}=\sqrt{\frac{m-k}{2\pi
k\left(  m-2k\right)  }}\frac{\left(  m-k\right)  ^{\left(  m-k\right)  }%
}{k^{k}\left(  m-2k\right)  ^{\left(  m-2k\right)  }}%
\]
and suggests%
\[
\left\vert \mathcal{A}[k,m]\right\vert \leq\binom{m-k}{k}\frac{\left(
N-1\right)  ^{m}}{\left(  N-1-\frac{m-2k}{m-k}\right)  ^{k}}%
\]
which, is, however worse.} turns out to be exact for $k=\frac{m}{2}$ and $m$
odd, but which is considerably higher for $\left(  k,m\right)  $-tuples.

\end{document}